%% file: paper.tex
\documentclass[onefignum,onetabnum,onealgnum,oneeqnum]{siamart190516}
\usepackage{amssymb}
\usepackage{amsmath}
\usepackage{enumitem}
\usepackage{booktabs}
\usepackage{multirow}
\usepackage{microtype}
\usepackage{bm}
\usepackage{algorithm}
\usepackage[noend]{algpseudocode}
\makeatletter
\algrenewcommand\ALG@beginalgorithmic{\small}
\renewcommand{\ALG@name}{\small Algorithm}
\makeatother
\setlist{leftmargin=1.5em}

\DeclareMathOperator*{\argmin}{arg\,min}
\newcommand{\ldbrack}{[\![}
\newcommand{\rdbrack}{]\!]}
\newcommand{\jump}[1]{\ldbrack #1 \rdbrack}
\newcommand{\trans}{\mathsf{T}}
\newcommand{\vtau}{\bm{\tau}}
\newcommand{\vsigma}{\bm{\sigma}}
\newcommand{\vu}{\mathbf u}
\newcommand{\vn}{\mathbf n}

\crefname{section}{\S\!}{\S\!}
\crefname{subsection}{\S\!}{\S\!}
\crefname{figure}{Fig.\!}{Figs.\!}

\newcommand{\vomega}{\bm{\omega}}
\def\trianglesymbol{\raisebox{-0.1em}{\rotatebox{90}{\scalebox{0.95}[0.75]{$\blacktriangle$}}}\hspace{0.08em}}
\def\squaresymbol{\raisebox{-0.05em}{\rotatebox{45}{\scalebox{0.55}{$\blacksquare$}}}}

\def\pentagonsymbol{\raisebox{0.05em}{\includegraphics[width=0.45em]{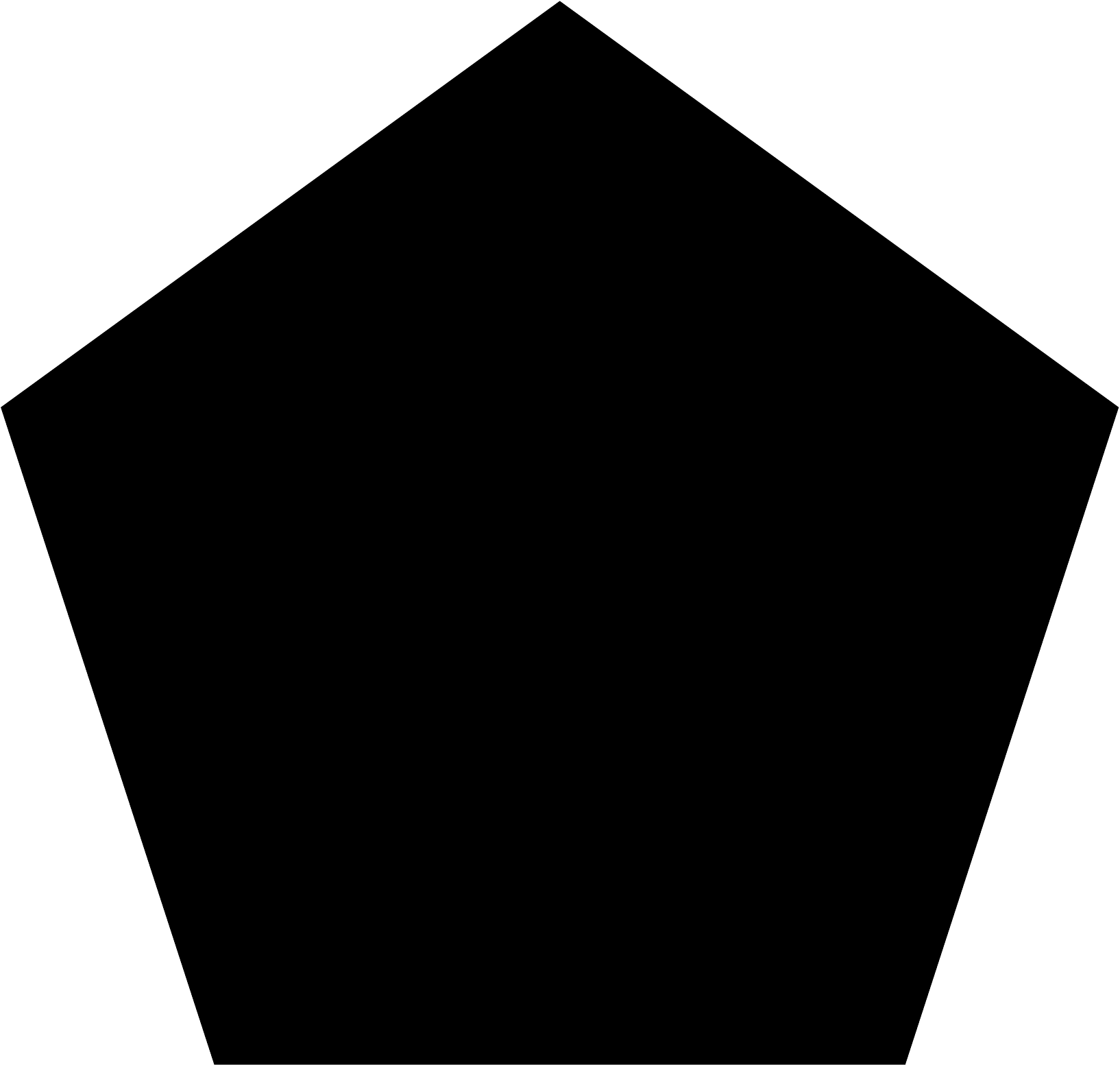}}}
\newcommand{\bmg}{\bm\eta}
\newcommand{\bmgg}{\eta}

\usepackage{etoolbox}
\makeatletter
\newif\if@algorithm
\patchcmd{\@float@Hx}% <cmd>
  {\@nodocument}% <search>
  {\@nodocument%
   \ifnum\pdfstrcmp{#1}{algorithm}=0 % Check whether using algorithm float
     \@algorithmtrue% Inside algorithm environment
     \setlength{\columnwidth}{\algorithmwidth}% Update column width
   \fi}% <replace>
  {}{}% <success><failure>
\renewcommand\float@makebox[1]{%
  \hbox{%
    \if@algorithm\hspace*{-.5\dimexpr\algorithmwidth-\textwidth}\fi%
    \vbox{\hsize=#1 \@parboxrestore
      \@fs@pre\@fs@iftopcapt
        \ifvoid\@floatcapt\else\unvbox\@floatcapt\par\@fs@mid\fi
        \unvbox\@currbox
      \else\unvbox\@currbox
        \ifvoid\@floatcapt\else\par\@fs@mid\unvbox\@floatcapt\fi
      \fi\par\@fs@post\vskip\z@}}}
\makeatother
\newlength{\algorithmwidth}
\AtBeginDocument{\setlength{\algorithmwidth}{0.8\textwidth}}

\title{Fast multigrid solution of high-order accurate multi-phase Stokes problems}

\author{Robert Saye\thanks{Mathematics Group, Lawrence Berkeley National Laboratory, Berkeley, CA 94720 (\texttt{rsaye@lbl.gov})}}
\date{\today}

\begin{document}

\maketitle

\begin{abstract}
A fast multigrid solver is presented for high-order accurate Stokes problems discretised by local discontinuous Galerkin (LDG) methods. The multigrid algorithm consists of a simple V-cycle, using an element-wise block Gauss-Seidel smoother. The efficacy of this approach depends on the LDG pressure penalty stabilisation parameter---provided the parameter is suitably chosen, numerical experiment shows that: (i) for steady-state Stokes problems, the convergence rate of the multigrid solver can match that of classical geometric multigrid methods for Poisson problems; (ii) for unsteady Stokes problems, the convergence rate further accelerates as the effective Reynolds number is increased. An extensive range of two- and three-dimensional test problems demonstrates the solver performance as well as high-order accuracy---these include cases with periodic, Dirichlet, and stress boundary conditions; variable-viscosity and multi-phase embedded interface problems containing density and viscosity discontinuities several orders in magnitude; and test cases with curved geometries using semi-unstructured meshes. 
\end{abstract}



\section{Introduction}

Stokes flow describes the motion of an incompressible viscous fluid at slow speeds, or small scales, and can be used to model a wide range of intricate phenomena, including mantle dynamics, the swimming of micro-organisms, the sedimentation of particulates, and the flotation of water droplets in clouds. In the steady-state case, the corresponding governing equations of motion are given by the Stokes equations, which generally take on one of two forms: either 
\begin{equation} \label{eq:stokes1} \begin{aligned} -\mu \nabla^2 \vu + \nabla p &= {\mathbf f} \\ \nabla \cdot \vu &= 0, \end{aligned} \end{equation}
or, alternatively,
\begin{equation} \label{eq:stokes2} \begin{aligned} -\nabla \cdot \bigl( \mu (\nabla \vu + \nabla \vu^\trans) \bigr) + \nabla p &= {\mathbf f} \\ \nabla \cdot \vu &= 0, \end{aligned} \end{equation}
where $\mu$ specifies the viscosity of the fluid, $\vu$ and $p$ describes its velocity and pressure fields, and $\mathbf f$ specifies the net external forces acting on the fluid. The form given in \eqref{eq:stokes1} is here referred to as the \textit{standard form} of the Stokes equations, and \eqref{eq:stokes2} as the \textit{viscous-stress form}, with the particular choice depending on the end application. For example, the viscous-stress form is generally applicable when $\mu$ is variable or if boundary conditions on stress are imposed.

Our motivation in this work is to develop fast multigrid solvers for computing high-order accurate solutions of the Stokes systems \eqref{eq:stokes1} or \eqref{eq:stokes2}, with extension also to multi-phase variants involving interfacial jump conditions in velocity and stress, as well as to time-dependent problems. In particular, we consider a framework based on local discontinuous Galerkin (LDG) methods \cite{ldg}, and build on prior work developing efficient multigrid algorithms for LDG discretisations of elliptic interface (Poisson-like) problems \cite{dgmg,fluxx}. We show that standard and simple-to-implement geometric multigrid algorithms can be applied to the resulting multi-phase Stokes problems---in the steady-state case, results show that the solver can match the speeds of fast geometric multigrid methods for Poisson problems; in the time-dependent case, convergence rates further accelerate as the effective Reynolds number is increased.

In particular, the presented multigrid algorithm consists of a standard V-cycle using an element-wise block Gauss-Seidel smoother---individual blocks correspond to individual mesh elements, such that the elemental degrees of freedom of both velocity and pressure are collected into the same block. 
Key to the rapid convergence of this approach is a suitable choice of the pressure penalty stabilisation parameter underlying the LDG framework---if the parameter is chosen well, then a highly efficient multigrid algorithm is obtained. 
We discuss how to choose this parameter for steady-state Stokes problems, and develop a simple strategy for generalising this choice to time-dependent Stokes problems.
Extensive tests of the multigrid methods are presented in this paper, including problems which impose Dirichlet or stress boundary conditions, variable-viscosity problems, test cases with curved geometry using semi-unstructured meshes, and multi-phase embedded interface problems with viscosity and density coefficients exhibiting discontinuities several orders in magnitude.

\subsection*{Previous work}

A vast amount of work in computational science and engineering has been devoted to the efficient solution of Stokes systems and saddle point problems in general; for an in-depth review of the correspondingly wide array of different approaches and their applications, see Benzi, Golub, and Liesen \cite{BenziGolubLiesen2005}. These approaches include, among others: block preconditioner methods, which operate on the viscosity, gradient, and divergence operator block structure of the Stokes equations; Schur complement methods, which manipulate, and usually approximate, the Schur complement of the saddle point system; and stationary iterative methods, such as the well-known Uzawa method, which alternates between updates of velocity and pressure, holding the other fixed, through the two governing equations in \eqref{eq:stokes1}. Here, we briefly review work on multigrid-style methods, particularly schemes based on the coupled solution of velocity and pressure, as is relevant to the present work; see also the reviews \cite{Wesseling2001,WesselingOosterlee2001,OosterleeGaspar2008}. 

As mentioned, the multigrid algorithm developed here uses a block Gauss-Seidel relaxation method, with each block collecting the velocity and pressure degrees of freedom on each mesh element. This approach is similar in essence to the ``symmetric coupled Gauss-Seidel method'' of Vanka \cite{Vanka1986,Vanka1986b} and can also be considered as a kind of ``box relaxation'' scheme \cite{BrandtLivne}. Vanka-type smoothers, originally devised for staggered-grid finite difference methods, visit each grid cell, solve for the velocity and pressure unknowns simultaneously via a local Stokes-like problem, and then move onto the next cell. Typically, a damping/under-relaxation parameter is needed to ensure convergence. In the original Vanka method, the Stokes system is restricted to the local variables in each grid cell and off-diagonal entries of the viscous operator are zeroed-out to facilitate a simpler update for the unknowns \cite{Wesseling2001,OlshanskiiTyrtyshnikov}. Variations have led to schemes which include the off-diagonal terms \cite{ThompsonFerziger1989,BorzacchielloLericheBlottiereGuillet2017}, line-based sweeping methods \cite{ThompsonFerziger1989,OosterleeWesseling1993b,OosterleeWesseling1993}, and have been examined with local mode analyses \cite{Sivaloganathan1991}. Vanka-type smoothers may also be considered as iterative Schwarz solvers, whereby the subdomains of the Schwarz method corresponds to the collection of degrees of freedom in each grid cell. Schemes building on this idea have since been developed for finite volume and finite element methods with much of the attention devoted to the choice of relaxation parameters, the choice of subdomains (e.g., whether to use one cell, or patches of cells), and on theoretical proofs of convergence in a multigrid setting, see, e.g., \cite{SchoberlZulehner2003,Manservisi2006,OosterleeGaspar2008,GasparNotayOosterleeRodrigo2014,Chen2015,HeMacLachlan2017,ColeyBenzakenEvans2017}; they have also found application in variable-viscosity Stokes problems \cite{BorzacchielloLericheBlottiereGuillet2017} and in computational solid mechanics \cite{WobkerTurek2009,HongKrausXuZikatanov2016}. As an example, in very recent work, Farrel, He, and MacLachlan \cite{FarrellHeMacLachlan2019} demonstrated the application of local Fourier analysis on these smoothers and found that smaller patches result in better convergence per floating point operation. Meanwhile, solvers specific to discontinuous Galerkin methods of the Stokes equations have also been devised; here, one possible approach is to exactly enforce the divergence constraint across the multigrid hierarchy through manipulation of the DG spaces, e.g., through $H(\text{div},\Omega)$-conforming discretisations. In the associated multigrid solvers, the divergence constraint is built into the coarse and fine mesh approximations, see, e.g., \cite{Chen2015,KanschatMao2015,HongKrausXuZikatanov2016,AdlerBensonMacLachlan2016,ColeyBenzakenEvans2017}. 

In many of these works, satisfactory multigrid convergence rates are reported, but they generally do not match the speeds of an efficient geometric multigrid method designed for scalar elliptic equations. In some cases, performance degrades as the mesh is refined, or as viscosity ratios increase, or in the case of time-dependent Stokes problems, as the Reynolds number changes. In contrast, for the LDG schemes devised here, we found that a simple block Gauss-Seidel relaxation method, which does not use any under- or over-relaxation parameters, can result in rapid multigrid convergence across a variety of challenging Stokes problems. Another advantage to a block Gauss-Seidel method is its simple implementation and the possibility of parallelism; for example, Gmeiner \textit{et al} \cite{GmeinerHuberJohnRudeWohlmuth2016} and Bauer \textit{et al} \cite{BauerKlementOberhuberZabka2016} demonstrate massively parallel and GPU implementations, respectively. For example, some of the three-dimensional tests in this paper used half a billion degrees of freedom and scaled to several hundred computing nodes, though we do not report on scaling performance here.

\subsection*{Outline}

In the main article, the central ideas and results are presented, while nonessential details of the LDG discretisation, grid convergence analyses, and implementation possibilities are deferred to the \iffalse Supplementary Material\fi Appendices. \iffalse The remainder of this article is organised as follows. \fi First, we outline the essential components of the LDG framework for the multi-phase Stokes equations. Second, the design of a standard multigrid V-cycle is outlined, afterwhich the role of pressure penalty stabilisation on multigrid efficiency is examined. In the remaining two sections, results are presented for a variety of problems for the steady-state and time-dependent Stokes equations, respectively. We then conclude, summarising the key observations made in this work along with a discussion of future research avenues.

\section{Local Discontinuous Galerkin Methods for Multi-Phase Stokes Problems}

In this work, we build on the LDG schemes developed by Cockburn \textit{et al} \cite{CockburnKanschatSchotzauSchwab2002} and extend them to the variable-viscosity multi-phase Stokes problem. The governing equations are written as follows: we seek to determine a velocity field $\vu : \Omega \to {\mathbb R}^d$ and pressure field $p : \Omega \to {\mathbb R}$ such that
\begin{equation} \label{eq:govern1} \left. \begin{aligned} -\nabla \cdot \bigl(\mu_i (\nabla \vu + \gamma \nabla \vu^\trans) \bigr) + \nabla p &= \mathbf f \\ -\nabla \cdot \mathbf \vu &= f \end{aligned} \right\} \text{ in $\Omega_i$,} \end{equation}
subject to the interfacial jump conditions,
\begin{equation} \label{eq:govern2} \left. \begin{aligned} \jump{\vu} &= {\mathbf g}_{ij} \\ \jump{\mu (\nabla \vu + \gamma \nabla \vu^\trans) \vn - p \vn} &= {\mathbf h}_{ij} \end{aligned} \right\} \text{ on $\Gamma_{ij}$,} \end{equation}
and boundary conditions,
\begin{equation} \label{eq:govern3} \begin{aligned} \vu &= {\mathbf g}_{\partial} && \text{on } \Gamma_D, \\ \mu (\nabla \vu + \gamma \nabla \vu^\trans) \vn - p \vn &= {\mathbf h}_{\partial} && \text{on } \Gamma_N, \end{aligned} \end{equation}
where $\Omega$ is a domain in $\mathbb R^d$ divided into one or more subdomains $\Omega_i$ (denoted ``phases''), $\Gamma_{ij} := \partial \Omega_i \cap \partial \Omega_j$ is the interface between phase $i$ and $j$, and $\Gamma_D$ and $\Gamma_N$ denote the parts of $\partial \Omega$ on which velocity Dirichlet or stress boundary conditions are imposed, respectively. Here, either $\gamma = 0$ or $\gamma = 1$ depending on whether the Stokes equations are in standard form or viscous-stress form, respectively. The operator $\jump{\cdot}$ denotes the jump in a quantity across an interface and $\vn$ is to be understood from context---on $\partial \Omega$, $\vn$ denotes the outward unit normal to the domain boundary, whereas for an interface $\Gamma_{ij}$, $\vn$ denotes the unit normal to $\Gamma_{ij}$, oriented consistently with the definition of the jump operator. Finally, $\mu_i$ is a phase-dependent viscosity coefficient, while $\mathbf f$, $f$, $\mathbf g$, and $\mathbf h$ provide the data to the multi-phase Stokes problem and are given functions defined on $\Omega$, its boundary, and internal interfaces.

Here, we mainly consider meshes arising from Cartesian grids along with semi-unstructured quadtree/octree-based implicitly defined meshes of more complex curved domains. In this setting, it is natural to adopt a tensor-product piecewise polynomial space. Let $\mathcal E = \bigcup_i E_i$ denote the set of mesh elements, let $p \geq 1$ be an integer,\footnote{The meaning of $p$, whether as pressure or polynomial degree, should be clear from context.} and define ${\mathcal Q}_p(E)$ as the space of tensor-product polynomials of (one-dimensional) degree $p$ on element $E$. For example, ${\mathcal Q}_2$ is the space of biquadratic or triquadratic polynomials, with dimension $9$ or $27$ in 2D or 3D, respectively. Define the corresponding space of discontinuous piecewise polynomial functions as
\[ V_h = \bigl\{ u : \Omega \to \mathbb R \enskip \bigl| \enskip u|_E \in \mathcal Q_p(E) \text{ for every } E \in \mathcal E \bigr\},\]
with analogous definitions for the space of piecewise polynomial vector-valued fields, $V_h^d$, and the space of matrix-valued fields, $V_h^{d \times d}$. As discussed in ref.~\cite{CockburnKanschatSchotzauSchwab2002}, it is possible to build LDG methods for the Stokes equations wherein the discrete pressure field has either the same polynomial degree as the discrete velocity field or is in a lower degree space. In this work, we focus on the case the two have the same degree, i.e., we seek a discrete solution such that $\vu_h \in V_h^d$ and $p_h \in V_h$. 

In one possible construction of the LDG framework, the governing set of equations \eqref{eq:govern1}--\eqref{eq:govern3} can be discretised in a three-step process: (i) define a discrete stress tensor $\vtau_h \in V_h^{d \times d}$ equal to the discretisation of $\nabla \vu_h + \gamma \nabla \vu_h^\trans$, taking into account Dirichlet source data $\mathbf g$; (ii) define $\vsigma_h \in V_h^{d \times d}$ as the viscous stress $\mu \vtau_h - p_h {\mathbb I}$ via an $L^2$ projection of $\mu \vtau_h$ onto $V_h^{d \times d}$; (iii) compute a discrete divergence of $\vsigma_h$, taking into account Neumann-like data $\mathbf h$, and add penalty stabilisation parameters for both velocity and pressure, setting the result equal to the $L^2$ projection of the given right hand side $\mathbf f$. Details of this construction, along with the associated treatment of the divergence constraint, are provided in \iffalse the Supplementary Material \fi \cref{sec:ldg}; here, we summarise the main outcomes of essential relevance. \iffalse to the remainder of this the main article.\fi The LDG discretisation results in a symmetric linear system for $(\vu_h,p_h)$ having the form
\begin{equation} \label{eq:blockform} \begin{pmatrix} A & M {\mathcal G} \\ {\mathcal G}^\trans M & -E \end{pmatrix} \begin{pmatrix} \vu_h \\ p_h \end{pmatrix} = \begin{pmatrix} {\mathbf b}_{\vu} \\ b_p \end{pmatrix} \end{equation}
where $M$ is the block-diagonal mass matrix, $\mathcal G$ is a discrete gradient operator, and $({\mathbf b}_{\vu}, b_p)$ collects the entire influence of the source data $\mathbf f$, $f$, $\mathbf g$, and $\mathbf h$ onto the right hand side. Here, $A$ implements the viscous part of the Stokes momentum equations, and can be written in $d \times d$ block form corresponding to its action on the $d$ components of $\vu_h$, with the $(i,j)$th block given by
\[ A_{ij} = \delta_{ij} \bigl( \textstyle{\sum_{k=1}^d} G_k^\trans M_\mu G_k \bigr) + \gamma\, G_j^\trans M_\mu G_i + \delta_{ij} \tilde{E}, \]
where $M_\mu$ is a $\mu$-weighted mass matrix, $G$ is a second discrete gradient operator closely related to the adjoint of ${\mathcal G}$, and $\tilde E$ is the operator associated with velocity penalty stabilisation. Note that if $\gamma = 0$, then $A$ is block diagonal with identical blocks corresponding to a discretisation of the Laplacian operator $-\nabla \cdot (\mu \nabla)$. Meanwhile, noting that the adjoint of $\mathcal G$ is given by $M^{-1} {\mathcal G}^\trans M$, one may observe that the divergence constraint of the Stokes equations is implemented in the $(p,\vu)$ block of \eqref{eq:blockform} through an effective discrete divergence operator which is equal to the negative adjoint of $\mathcal G$. 

There is one last operator to define in \eqref{eq:blockform}, whose presence is of key importance to multigrid efficiency: the pressure stabilisation operator $E$, which weakly enforces continuity of the pressure field. The symmetric positive semidefinite matrix $E$ is defined such that\footnote{In a convenient abuse of notation, a piecewise polynomial function (e.g., in $V_h$) may carry the same notation as its corresponding coefficient vector in the basis of $V_h$, with the precise meaning understood from context. For example, in the identity $\smash{u^\trans M v = \scalebox{0.8}{$\int_\Omega$} u\,v}$, the left hand side employs vectors and matrices relative to the chosen basis of $V_h$, whereas the right hand side employs the functional form.}
\begin{equation} \label{eq:eop} u^\trans E v = \int_{\Gamma_0} \tau_p \,\jump{u} \jump{v} \end{equation}
holds for every $u, v \in V_h$. Here, the integral is taken over the union of every non-interfacial interior mesh face, $\jump{\cdot}$ denotes the jump across the face, and $\tau_p$ is a \textit{pressure stabilisation penalty parameter} which scales proportional to the element size $h$ and inversely proportional to the (local) viscosity coefficient:
\begin{equation} \label{eq:taueqn} \tau_p = \tau\,h / \mu, \end{equation}
where $\tau$ is a user-defined constant prefactor. Provided $\tau$ is positive, Cockburn \textit{et al} \cite{CockburnKanschatSchotzauSchwab2002} (see also extensions \cite{CockburnKanschatSchotzau2003,CockburnKanschatSchotzau2004}) prove the well-posedness of the single-phase symmetric saddle point problem in \eqref{eq:blockform}, including satisfaction of the inf-sup conditions. 
In these cited works, however, the particular choice of $\tau$ in \eqref{eq:taueqn} is not extensively discussed. One of the main results in this work is to demonstrate that $\tau$ can be chosen so as to achieve excellent multigrid solver efficiency when computing solutions to the Stokes equations.

In the remainder of this paper, we will refer to the Stokes problem both in its operator block form \eqref{eq:blockform} and through the more succinct notation 
\[ {\mathcal A}_h x_h = b_h, \]
where ${\mathcal A}_h$ is the symmetric saddle point operator and $x_h$ collects $\vu_h$ and $p_h$ into one set of unknowns.

\section{Multigrid Design}

Prior work on designing geometric multigrid methods for LDG discretisations of \iffalse elliptic \fi Poisson-like equations \cite{ImplicitMeshPartOne,ImplicitMeshPartTwo,dgmg,fluxx} shows that one can build an efficient solver through standard multigrid concepts: a V-cycle applied to a mesh hierarchy using straightforward interpolation and restriction operators, together with standard relaxation methods, such as block Gauss-Seidel in which each block corresponds to the collective set of unknowns on each mesh element. Here, we show the same can be done for the discretised Stokes problem \eqref{eq:blockform}. (In the following, it is assumed the reader is familiar with the general design of multigrid methods; see, e.g., the books \cite{Briggs_00_01,BrandtLivne,Wesseling2001,OlshanskiiTyrtyshnikov} for reviews and applications.)

The multigrid methods designed here may be considered as a ``purely-geometric'' approach, wherein the Stokes problem is discretised on each level of the mesh hierarchy. (A convenient strategy for constructing the coarse-mesh problems---without having to explicitly form the coarse meshes themselves---is discussed shortly.) Three preliminary ingredients are needed to specify its design:
\begin{itemize}
\item \textit{Mesh hierarchy.} In this work, quadtrees and octrees are used to define the finest mesh. The tree structure naturally defines a hierarchical procedure by which to agglomerate elements to create a nested mesh hierarchy,  coarsening by a factor of two in each dimension down each level. Regarding the multi-phase case, element agglomeration is permitted only between elements of the same phase---as such, the interface is sharply preserved throughout the entire multigrid hierarchy.
\item \textit{Interpolation operator.} Owing to the presence of a nested mesh hierarchy, the interpolation operator $I_{2h}^h$, which transfers coarse mesh corrections to a fine mesh, is naturally defined via injection. In particular, we define $(I_{2h}^h u)|_{E_f} = u|_{E_c}$, where $E_f$ is a fine mesh element and $E_c \supseteq E_f$ is its corresponding coarse mesh element.
\item \textit{Restriction operator.} The restriction operator $R_h^{2h}$, which transfers the residual of a fine-mesh problem to the coarse mesh, is defined as the $L^2$ projection onto the coarse mesh (or, equivalently, as the adjoint of the interpolation operator). It is related to the interpolation operator via $R_h^{2h} = M_{2h}^{-1} (I_{2h}^h)^\trans M_h$, where $M_h$ and $M_{2h}$ are the mass matrices of the two meshes.
\end{itemize}
The last essential multigrid ingredient, and perhaps most important, is the relaxation method. As mentioned earlier, we have used a simple block Gauss-Seidel method, where each block corresponds to the collective set of degrees of freedom (i.e., velocity and pressure combined) on each mesh element. Specifically, consider a repartitioning of $\mathcal A x = b$ according to these blocks, such that $x_i$ denotes the set of velocity and pressure values on element $i$, and ${\mathcal A}_{ij}$ denotes the $(i,j)$th block of ${\mathcal A}$, whence $b_i = \sum_j {\mathcal A}_{ij} x_j$. Then, the block Gauss-Seidel method simply sweeps over the elements, in some particular order, replacing $x_i \leftarrow {\mathcal A}_{ii}^{-1}(b_i - \sum_{j \neq i} {\mathcal A}_{ij} x_j)$. Here, ${\mathcal A}_{ii}$ is the $i$th diagonal block of ${\mathcal A}$ and takes on the form of a miniature Stokes operator; referring to \eqref{eq:blockform}, we have
\[ {\mathcal A}_{ii} = \begin{pmatrix} A_{ii} & M_{ii} {\mathcal G}_{ii} \\ {\mathcal G}_{ii}^\trans M_{ii} & -E_{ii} \end{pmatrix}. \]
Note that ${\mathcal A}_{ii}$ needs to be inverted in the Gauss-Seidel update of element $i$. Assuming that the global Stokes saddle-point problem ${\mathcal A} x = b$ satisfies the inf-sup conditions, it is straightforward to show that so too does ${\mathcal A}_{ii}$, and hence the local element-wise problem is well-posed; this has also been confirmed through numerous and extensive numerical tests. In our specific implementation, we precompute a symmetric indefinite factorisation of ${\mathcal A}_{ii}$ for every $i$, and use this factorisation as a direct solver for each of these mini-Stokes problems in the Gauss-Seidel sweep.
Regarding the element ordering, we have opted for a multi-coloured Gauss-Seidel method. The primary reason for this choice is that a multi-coloured sweep affords a simpler parallel implementation of the method, both in terms of multi-threading and in a distributed environment (e.g., through standard domain decomposition methods using MPI). 

\begin{algorithm}[t]
	\caption{\sffamily\small Multigrid V-cycle $V({\mathcal E}_h, x_h, b_h)$ with $\nu_1$ pre- and $\nu_2$ post-smoothing steps on mesh ${\mathcal E}_h$ of the hierarchy}
	\begin{algorithmic}[1]
		\If{${\mathcal E}_h$ is the bottom level}
			\State Solve ${\mathcal A}_h x_h = b_h$ with bottom solver
		\Else
			\State Apply block Gauss-Seidel relaxation $\nu_1$ times
			\State $r_{2h} := (I_{2h}^h)^\trans (b_h - {\mathcal A}_h x_h)$
			\State $x_{2h} := V({\mathcal E}_{2h}, 0, r_{2h})$
			\State $x_h \leftarrow x_h + I_{2h}^h x_{2h}$
			\State Apply block Gauss-Seidel relaxation $\nu_2$ times
		\EndIf
		\State{\textbf{return} $x_h$}
	\end{algorithmic}%
	\label{algo:vcycle}%
\end{algorithm}%

Using the defined interpolation and restriction operators and the block Gauss-Seidel relaxation method, the construction of a multigrid V-cycle is relatively standard and is outlined in \cref{algo:vcycle}.\footnote{Note that $(I_{2h}^h)^\trans$ appears on line 5, rather than the restriction operator $R_h^{2h}$; this follows from a convenient simplification common to many finite element methods: briefly, viewed as an operator which maps $\smash{V_h^d \otimes V_h}$ to $\smash{V_h^d \otimes V_h}$, the discrete Stokes operator is given by $\smash{M^{-1} {\mathcal A}}$. Therefore, the residual of the fine mesh problem, as a piecewise polynomial function, is $\smash{M_h^{-1} (b_h - {\mathcal A}_h x_h)}$. This residual is then multiplied by $\smash{R_h^{2h}}$ to define the source data for the coarse mesh problem $\smash{M_{2h}^{-1} {\mathcal A}_{2h} = R_h^{2h} M_h^{-1} (b_h - A_h x_h)}$. Rearranging, one obtains line 5.} In this algorithm, $\mathcal A_h$ is assumed to be pre-computed on every level of the mesh; a particularly convenient method for doing so---without having to explicitly mesh each level, build quadrature rules for coarse mesh elements, or build LDG operators via coarse-mesh numerical fluxes, etc.---uses the operator-coarsening ideas of ref.~\cite{dgmg}. In this technique, the discrete gradient and penalty operators underlying \eqref{eq:blockform} are coarsened solely based on the interpolation operator hierarchy, using simple block-sparse linear algebra; these methods are further described in \iffalse the Supplementary Material \fi \cref{sec:op}. Returning to \cref{algo:vcycle}, note that the V-cycle computes coarse-grid corrections for both the velocity and pressure, simultaneously, and there is no need to strictly enforce the divergence constraint on any level of the hierarchy. Regarding the bottom solver, in this work a direct solver using a symmetric indefinite factorisation of ${\mathcal A}_h$ on the coarsest mesh is used, together with an appropriate treatment of its associated trivial kernel.\footnote{With periodic boundary conditions, the Stokes problem has a trivial kernel of dimension $d+1$, spanned by constant velocity and pressure fields; with velocity Dirichlet boundary conditions, the kernel is one-dimensional, spanned by constant pressure fields; with stress boundary conditions in viscous-stress form, the kernel is spanned by constant velocity fields as well as less trivial velocity modes such as, e.g., the velocity field $(x,y) \mapsto (-y,x)$ in 2D. The bottom solver robustly treats these modes through a simple least squares approach which (pre)computes the symmetric eigendecomposition of $\mathcal A_h$, ``snapping'' any nearly-zero eigenvalues to exactly zero.} 

As is typical, applying more and more pre- and post-smoothing steps increases the convergence rate of the multigrid solver, but at greater computational cost. According to a variety of numerical experiments, a general observation made in this work is that a V-cycle with three pre- and post-smoothing steps is a good all-rounder, based on the metric of fastest computation time in reducing solution error by a given factor. On occasion, four pre-smoothing steps and two post-smoothing steps, or vice versa, performs marginally better, but on a problem-specific basis. Naturally, the optimal choice of multigrid design parameters is implementation- and problem-dependent, influenced by a wide variety of aspects, e.g., the relative computational costs of interpolation, restriction, and relaxation operators, or computing hardware characteristics, such as shared memory or distributed memory architectures and their associated memory communication costs. Further comments on V-cycle design, or counterparts such as W-cycles, are provided in the concluding remarks.

To complete the description of the multigrid method, we note that although the V-cycle can be used as a standalone iterative solver, solver efficiency can be further accelerated by using it as a preconditioner of a Krylov method \cite{BenziGolubLiesen2005}. In this work, we have used  a single V-cycle (with $\nu_1=\nu_2=3$ pre- and post-smoothing steps) as a left-preconditioner of the GMRES method. Specifically, applying the V-cycle to an initial guess of zero on the fine mesh results in a linear operator, denoted in the remainder of this paper as $V$; the preconditioned system is then $V \mathcal A$. For simplicity, we do not consider restarted variants of GMRES here, in part because experiments show that convergence is generally attained in as few as 5-15 steps for a ten-fold reduction in the order of magnitude of the residual. Moreover, convergence behaviour is generally smooth during the iterations, such that the residual reduces in norm by a nearly constant factor each iteration of the GMRES method.

\section{Influence of Pressure Penalty Stabilisation on Multigrid Efficiency}
\label{sec:tauinfluence}

As outlined, the Stokes multigrid solver consists of a standard V-cycle, using a block Gauss-Seidel relaxation method in which each block corresponds to the collective set of degrees of freedom, of both velocity and pressure, on each mesh element. Key to rapid multigrid convergence is an apt choice of the user-defined pressure penalty stabilisation prefactor parameter $\tau$ in \eqref{eq:taueqn}. In general terms, if $\tau$ is below some positive threshold, the V-cycle fails to converge; above this threshold, there is a range of values for which convergence rates can match that of fast geometric multigrid methods for scalar Poisson problems; and, beyond this range, multigrid efficiency will degrade.

\begin{figure}%[tbhp]
\centering
\includegraphics[scale=0.91]{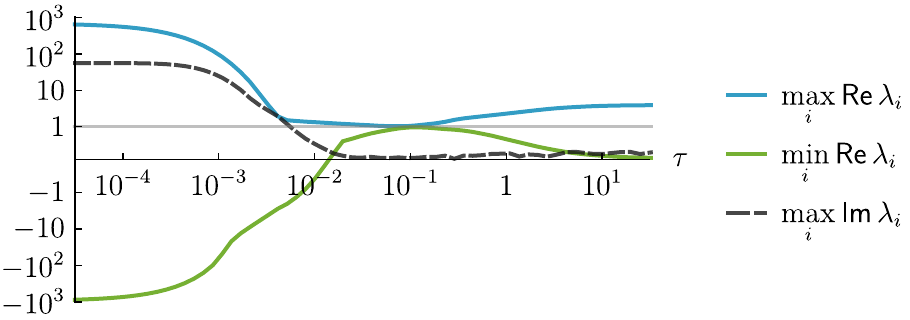}
\caption{Spectral properties of the multigrid-preconditioned system $V{\mathcal A}$ as a function of pressure penalty stabilisation parameter; note the quasi-logarithmic axis. Results correspond to a two-dimensional Stokes problem in standard form with $p = 2$; similar characteristics are obtained with other $p$, or in 3D, or with the viscous-stress form of the Stokes equations.}
\label{fig:tauspectral}
\end{figure}

To illustrate this behaviour, as well as the suitability of the V-cycle as a preconditioner for the Stokes system ${\mathcal A}$, we examine the spectral properties of the preconditioned system $V {\mathcal A}$. In this particular example, we consider the standard form of the Stokes equations ($\gamma = 0$), with periodic boundary conditions on a $128 \times 128$ Cartesian grid using $p = 2$ biquadratic elements. As a function of $\tau$, \cref{fig:tauspectral} plots three quantities concerning the spectrum\footnote{The spectral analysis deliberately excludes the zero eigenvalues associated with the trivial kernel of the Stokes operator. In addition, regarding the results plotted in \cref{fig:tauspectral}, the spectral extremes have been estimated via the GMRES method, through computation of the eigenvalues of the upper Hessenberg matrix of the corresponding Arnoldi iteration. Although only an approximation to the true spectrum of $V{\mathcal A}$, estimation via GMRES is significantly more efficient than forming and computing the spectrum of the dense matrix $V{\mathcal A}$; furthermore, the computed quantities are found in practice to be sufficiently accurate for the present purpose.} of $V{\mathcal A}$: (i) the real part of the right-most eigenvalue, $\max_i \textsf{Re}\,\lambda_i$; (ii) the real part of the left-most eigenvalue, $\min_i \textsf{Re}\,\lambda_i$; and (iii) the greatest imaginary part, $\max_i \textsf{Im}\,\lambda_i$. Ideally, the eigenvalues of $V{\mathcal A}$ should be clustered around $1$, and we observe this is the case when $\tau \approx 0.1$; furthermore, near this value, the eigenvalues are nearly real. However, if $\tau$ is too small, then the left-most eigenvalue of $V{\mathcal A}$ crosses the imaginary axis; in \cref{fig:tauspectral} this occurs when $\tau \lessapprox 10^{-2}$. This represents a breakdown of the V-cycle, as then the (non-trivial) eigenvalues of $V{\mathcal A}$ cease to be bounded away from zero. (In fact, numerical experiments examining the efficacy of the V-cycle as a standalone iterative method show that the V-cycle ceases to have spectral radius less than one at this same point.) Meanwhile, if $\tau$ is too large, e.g., $\tau \gtrapprox 1$, then the eigenvalues of $V{\mathcal A}$ begin to diverge away from $1$. 

\begin{figure}%[tbhp]
\centering
\includegraphics[scale=0.91]{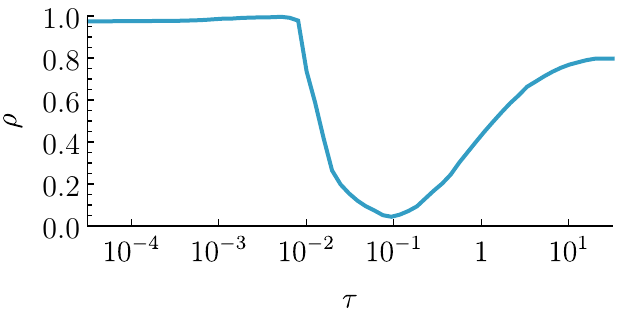}
\caption{Multigrid efficiency as a function of pressure penalty stabilisation parameter. Here, the convergence rate $\rho$, defined by eq.~\ref{eq:rho}, measures the average reduction factor per iteration in the residual of the V-cycle preconditioned GMRES method.}
\label{fig:taurho}
\end{figure}

In this work, the primary metric used to assess multigrid efficiency is the convergence rate $\rho$ of the multigrid-preconditioned GMRES method. This metric correlates with the spectral characteristics described above, and also to the performance of the V-cycle as a standalone iterative solver. Here, $\rho$ is defined as the average residual reduction factor per iteration of the (left) preconditioned GMRES method,
\begin{equation} \label{eq:rho} \rho = \exp \Bigl( \frac1n \log \frac{\|V{\mathcal A} x_n - V b\|_2}{\|V{\mathcal A} x_0 - V b\|_2} \Bigr), \end{equation}
where $n$ is the number of iterations required to reduce the residual by a factor of $10^8$ from its starting value. In particular, a right hand side of $b = 0$ is used, with initial guess $x_0$ given by a randomly generated $(d+1)$-dimensional vector field. With high probability, the randomly generated field contains modes which are damped slowest by the multigrid method, and thus $\rho$ represents a typical ``worst case'' convergence rate. \cref{fig:taurho} shows the convergence rate $\rho$ as a function of $\tau$ for the same example considered in \cref{fig:tauspectral}. Optimal convergence is attained when $\tau \approx 0.1$, precisely when the spectrum of $V{\mathcal A}$ is tightly clustered around 1; meanwhile, when $\tau \lessapprox 10^{-2}$, $\rho$ is approximately one, representing the fact that GMRES is unable to effectively reduce the residual owing to a breakdown of the preconditioner. 

Similar behaviour to that seen in \cref{fig:tauspectral} and \cref{fig:taurho} is observed for other choices of polynomial degree $p$, for different grid sizes, for the viscous-stress form of the Stokes equations, and in 3D as well as 2D. In all cases, the convergence rate $\rho$ exhibits a well-defined valley as a function of $\tau$. At present, a formula for the corresponding optimal value of $\tau$ is not known. In this work, a simple one-dimensional parameter sweep was used to find the optimal value of $\tau$ on successive grid sizes $n \times n \,(\,\times n)$ for $n = 4,8,16,\ldots$ up to $n = 256$ in 2D and $n = 128$ in 3D; experiments indicated that the optimum essentially converges around $n = 64$ or $128$, beyond which $\argmin_{\tau} \rho$ is relatively insensitive to the grid size. \cref{tab:optimaltau} contains the results of the search for a variety of $p$ in 2D and 3D. In addition to the optimal value, a ``window'' of acceptable $\tau$ may also be computed: one can search for all $\tau$ such that $\rho(\tau) \leq (\min \rho)^{1 - \epsilon}$, where $0 \leq \epsilon < 1$ represents a user-defined threshold for which the number of multigrid iterations increases by a factor of $1/(1-\epsilon)$ above the optimal minimum. 
For example, with $\epsilon = 1/9$, the corresponding range of $\tau$ values yields at most 12.5\% more iterations than optimal; \iffalse Table SM1 in the Supplementary Material \fi \cref{tab:tauwindow} contains the corresponding ranges, and shows that, even if $\tau$ is not chosen exactly at the optimum, there is a relatively wide range of values for $\tau$ that will nevertheless attain good multigrid efficiency. Finally, numerical experiments indicate that the pressure penalty parameter has very little influence on the velocity or pressure discretisation error (see, for example, the demonstration given in \cref{app:tauaccuracy}); this is ideal, as it allows us to concentrate mainly on the impact of $\tau$ on multigrid performance. Further comments concerning the selection of $\tau$ are given in the concluding remarks.

\begin{table}%[tbhp]
\centering%
\sffamily\footnotesize%
\caption{Optimal values of pressure penalty stabilisation parameter, attaining minimal multigrid iteration count.}
\begin{tabular}{lcccccc}
\multirow{2}{*}{Stokes form} & \multirow{2}{*}{$d$} & \multicolumn{5}{c}{Polynomial degree $p$} \\
 & & 1 & 2 & 3 & 4 & 5 \\
\midrule
\multirow{2}{*}{Standard} & 2D & 0.19 & 0.10 & 0.086 & 0.019 & 0.031 \\ 
& 3D & 0.12 & 0.088 & 0.084 & -- & -- \\
\midrule
\multirow{2}{*}{Viscous stress} & 2D & 0.14 & 0.046 & 0.034 & 0.0095 & 0.011 \\
& 3D & 0.12 & 0.039 & 0.040 & -- & -- \\
\bottomrule
\end{tabular}
\label{tab:optimaltau}
\end{table}

\begin{table}%[tbhp]
\centering%
\caption{Range of pressure penalty stabilisation parameters for which the number of multigrid iterations is at most 12.5\% more than optimal.}%
\sffamily%
\footnotesize%
\newcommand{\ww}{0.83em}%
\begin{tabular}{l@{\hspace{\ww}}c@{\hspace{\ww}}c@{\hspace{\ww}}c@{\hspace{\ww}}c@{\hspace{\ww}}c@{\hspace{\ww}}c}
\multirow{2}{*}{Stokes form} & \multirow{2}{*}{$d$} & \multicolumn{5}{c}{Polynomial degree $p$} \\
 & & 1 & 2 & 3 & 4 & 5 \\
\midrule
\multirow{2}{*}{Standard} & 2D & {(0.15, 0.28)} & {(0.082, 0.12)} & {(0.067, 0.11)} & {(0.013, 0.029)} & {(0.021, 0.041)} \\ 
& 3D & {(0.061, 0.30)} & {(0.070, 0.11)} & {(0.064, 0.12)} & -- & -- \\
\midrule
\multirow{2}{*}{Viscous stress} & 2D & {(0.091, 0.18)} & {(0.040, 0.056)} & {(0.027, 0.043)} & {(0.0058, 0.021)} & {(0.0072, 0.020)} \\
& 3D & {(0.025, 0.19)} & {(0.031, 0.066)} & {(0.021, 0.059)} & -- & -- \\
\bottomrule
\end{tabular}
\label{tab:tauwindow}
\end{table}

\section{Multigrid Efficiency for the Time-Independent Stokes Equations} 
\label{sec:timeindep}

In the next two sections, multigrid performance is examined for a variety of Stokes problems, in both standard and viscous-stress form, and with different types of boundary conditions. We also consider test cases with variable viscosity, multi-phase problems exhibiting large discontinuities in $\mu$ across an embedded interface, and curved geometry problems which use semi-unstructured meshes. Our primary focus is on demonstrating effective multigrid performance under the action of mesh refinement. The order of accuracy in the velocity and pressure, in the $L^2$ and maximum error norms, is also measured and reported. In two dimensions, grid sizes typically range from $4 \times 4$ up to $1024 \times 1024$, with polynomial degrees $p = 1, 2, \ldots, 5$; in three dimensions, owing to limited computing resources, only $p = 1, 2,$ and $3$ is considered (i.e., trilinear, triquadratic, and tricubic polynomials) on grid sizes up to $128 \times 128 \times 128$ (the largest of these problems has half a billion degrees of freedom and requires 1TB of memory to store just the block-diagonal component of the block-sparse matrix ${\mathcal A}$). In the remainder of this paper, for every test case, the pressure penalty stabilisation parameter is chosen equal to the values reported in \cref{tab:optimaltau}. Multigrid convergence rates are assessed using the average reduction factor in residual per iteration of the GMRES method, i.e., using \eqref{eq:rho}, on a test problem with right hand side $b = 0$, initial guess defined by a randomly generated $(d+1)$-dimensional vector field, over as many iterations as necessary to reduce the initial residual by a factor of $10^8$.

\subsection{Periodic boundary conditions}
\label{prob:periodic}

\begin{figure}%[tbhp]
\centering%
\includegraphics[scale=0.91]{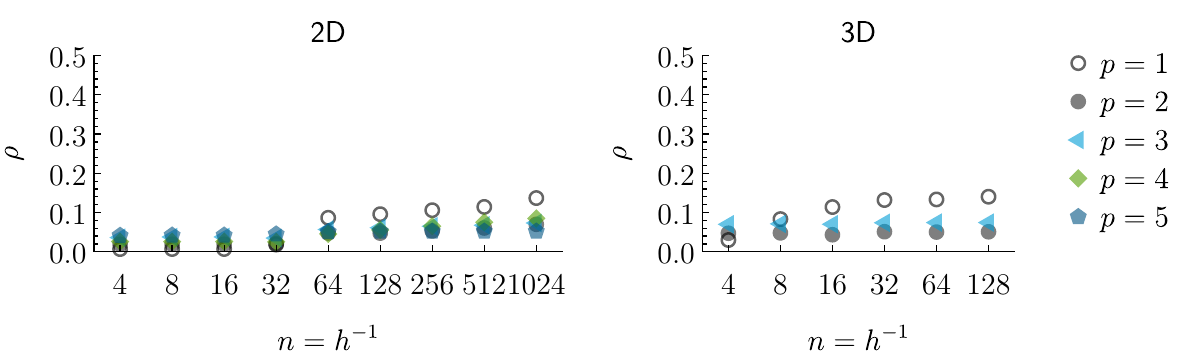}%
\caption{Measured multigrid convergence rates when solving the Stokes equations in standard form, with $\mu = 1$ and periodic boundary conditions.}%
\label{fig:periodic-stokes}%
\end{figure}

We begin with perhaps the simplest Stokes problem, i.e., the Stokes equations in standard form, with $\mu = 1$ and periodic boundary conditions, on the unit square/cube domain $\Omega = (0,1)^d$. \Cref{fig:periodic-stokes} plots the measured convergence rate $\rho$ as a function of grid size $n \times n\,(\times\, n)$. With the exception of $p = 1$, ideal convergence rates are attained. For example, $\rho \approx 0.05$ corresponds to needing only seven iterations to reduce the residual by a factor of $10^9$. When $p = 1$, however, slower convergence is seen---compared to the higher-degree cases, it appears that the an element-wise block Gauss-Seidel relaxation method is less effective for a bilinear and trilinear LDG discretisation of the Stokes equations. This behaviour is consistently observed across all of the presented tests; see also the next two examples. Regarding the order of accuracy, numerical experiments show that both velocity and pressure attain order $p + 1$ accuracy, in both the $L^2$ and maximum error norms, for all $p$ considered, in 2D and 3D; see \iffalse Fig.~SM3 in the Supplementary Material. \fi \cref{sec:gridconv}.

\subsection{Dirichlet boundary conditions}
\label{prob:dirichlet}

\begin{figure}%[tbhp]
\centering%
\includegraphics[scale=0.91]{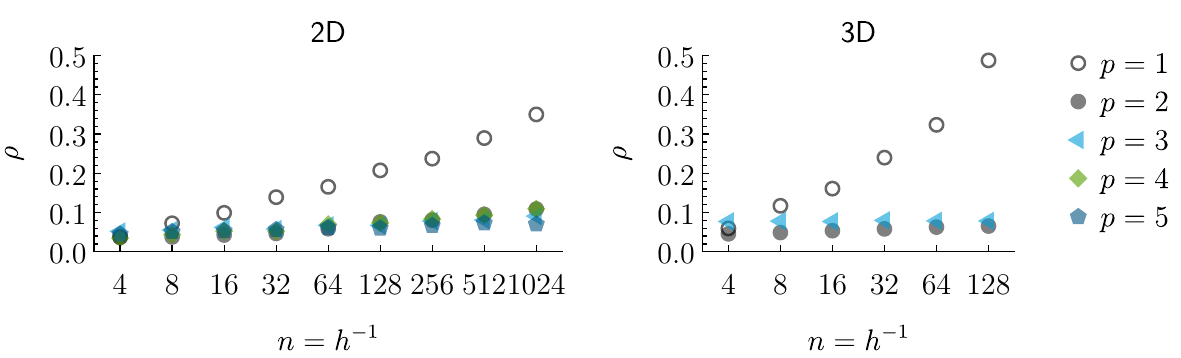}%
\caption{Measured multigrid convergence rates when solving the Stokes equations in standard form, with $\mu = 1$ and velocity Dirichlet boundary conditions.}%
\label{fig:dirichlet-stokes}%
\end{figure}

\begin{figure}
\centering
\newcommand{\ww}{1.5in}
\begin{tabular}{ccc}
{\setlength{\fboxsep}{0pt}\setlength{\fboxrule}{0.5pt}\fbox{\includegraphics[width=\ww]{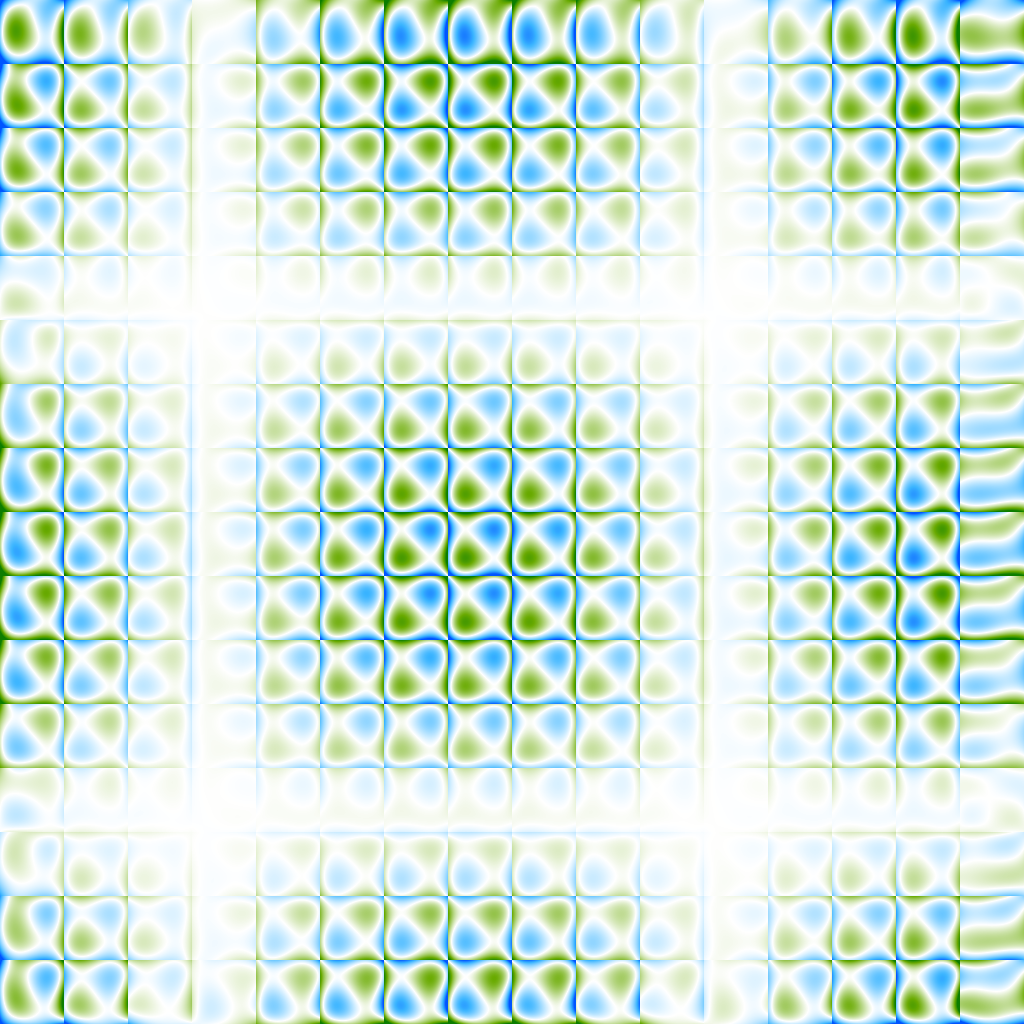}}}
 & {\setlength{\fboxsep}{0pt}\setlength{\fboxrule}{0.5pt}\fbox{\includegraphics[width=\ww]{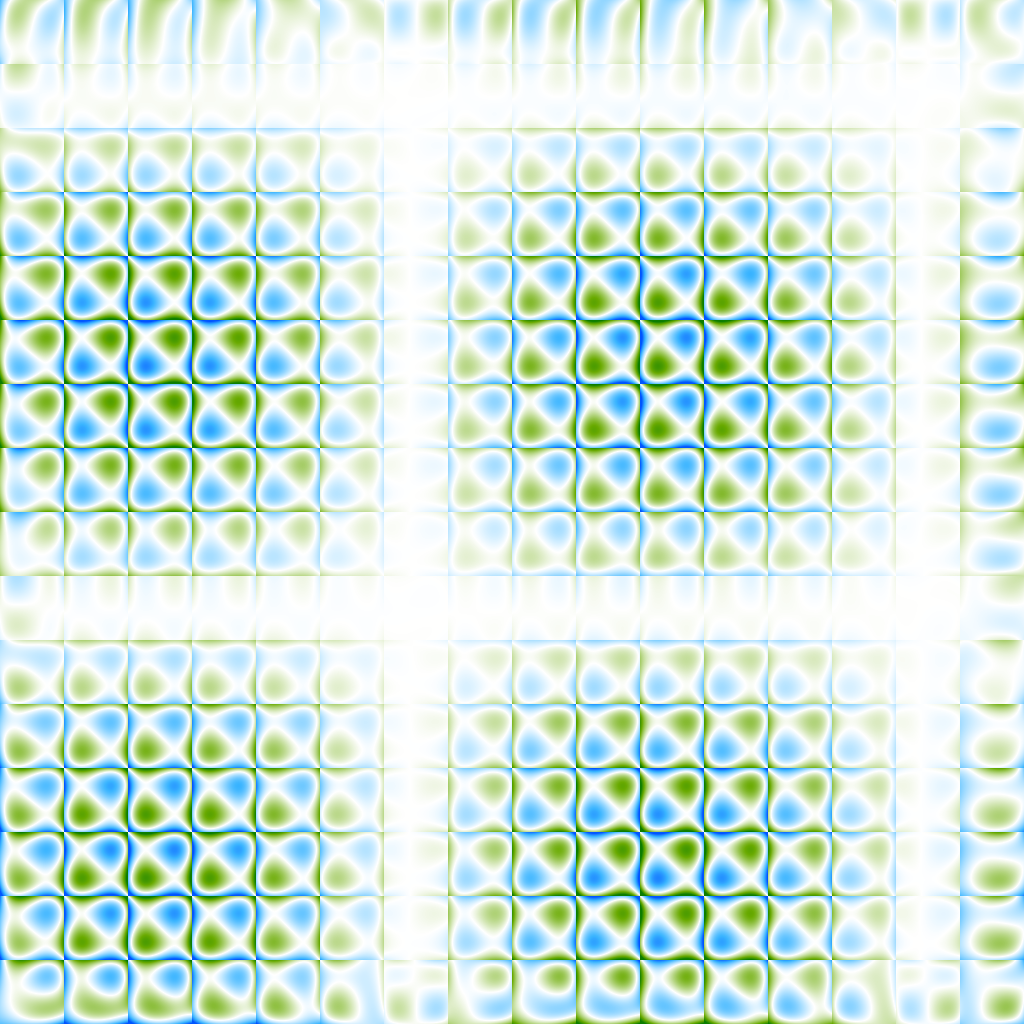}}}
 & {\setlength{\fboxsep}{0pt}\setlength{\fboxrule}{0.5pt}\fbox{\includegraphics[width=\ww]{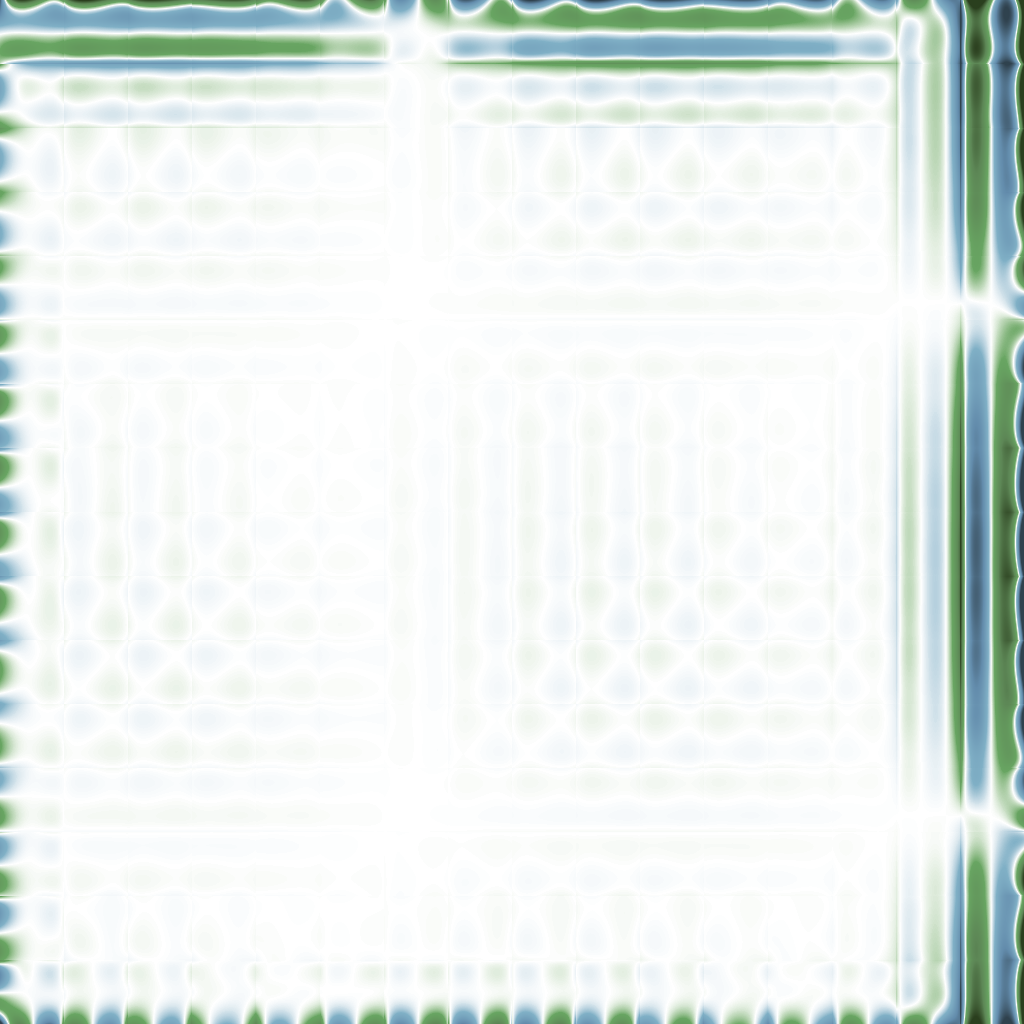}}} \\[1em]
\makebox[2cm][c]{{\begingroup%
\setlength{\unitlength}{2cm}%
\def\ticwidth{0.03}%
\def\ticsep{-0.02}%
\setlength\fboxsep{0pt}%
\begin{picture}(1,0.1)%
	\put(0,\ticwidth){\fbox{\includegraphics[width=0.995\unitlength,height=0.05\unitlength]{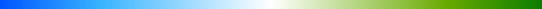}}}%
	\put(0.0035,0){\line(0,1){\ticwidth}}%
	\put(0.25,0){\line(0,1){\ticwidth}}%
	\put(0.5,0){\line(0,1){\ticwidth}}%
	\put(0.75,0){\line(0,1){\ticwidth}}%
	\put(1.0055,0){\line(0,1){\ticwidth}}%
	\put(0.0,\ticsep){\makebox(0,0)[ct]{\rotatebox{0}{\scalebox{0.6}{$-0.000015\hphantom{-}$}}}}%
	\put(0.5,\ticsep){\makebox(0,0)[ct]{\rotatebox{0}{\scalebox{0.6}{$0$}}}}%
	\put(1,\ticsep){\makebox(0,0)[ct]{\rotatebox{0}{\scalebox{0.6}{$0.000015$}}}}%
	\put(0.5,0.23){\makebox(0,0)[ct]{\rotatebox{0}{\scalebox{0.75}{$u-u_h$}}}}%
\end{picture}%
\global\let\ticwidth\undefined%
\global\let\ticsep\undefined%
\endgroup}}
&
\makebox[2cm][c]{{\begingroup%
\setlength{\unitlength}{2cm}%
\def\ticwidth{0.03}%
\def\ticsep{-0.02}%
\setlength\fboxsep{0pt}%
\begin{picture}(1,0.1)%
	\put(0,\ticwidth){\fbox{\includegraphics[width=0.995\unitlength,height=0.05\unitlength]{bluegreen}}}%
	\put(0.0035,0){\line(0,1){\ticwidth}}%
	\put(0.25,0){\line(0,1){\ticwidth}}%
	\put(0.5,0){\line(0,1){\ticwidth}}%
	\put(0.75,0){\line(0,1){\ticwidth}}%
	\put(1.0055,0){\line(0,1){\ticwidth}}%
	\put(0.0,\ticsep){\makebox(0,0)[ct]{\rotatebox{0}{\scalebox{0.6}{$-0.000015\hphantom{-}$}}}}%
	\put(0.5,\ticsep){\makebox(0,0)[ct]{\rotatebox{0}{\scalebox{0.6}{$0$}}}}%
	\put(1,\ticsep){\makebox(0,0)[ct]{\rotatebox{0}{\scalebox{0.6}{$0.000015$}}}}%
	\put(0.5,0.23){\makebox(0,0)[ct]{\rotatebox{0}{\scalebox{0.75}{$v-v_h$}}}}%
\end{picture}%
\global\let\ticwidth\undefined%
\global\let\ticsep\undefined%
\endgroup}}
&
\makebox[2cm][c]{{\begingroup%
\setlength{\unitlength}{2cm}%
\def\ticwidth{0.03}%
\def\ticsep{-0.02}%
\setlength\fboxsep{0pt}%
\begin{picture}(1,0.1)%
	\put(0,\ticwidth){\fbox{\includegraphics[width=0.995\unitlength,height=0.05\unitlength]{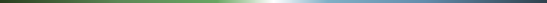}}}%
	\put(0.0035,0){\line(0,1){\ticwidth}}%
	\put(0.25,0){\line(0,1){\ticwidth}}%
	\put(0.5,0){\line(0,1){\ticwidth}}%
	\put(0.75,0){\line(0,1){\ticwidth}}%
	\put(1.0055,0){\line(0,1){\ticwidth}}%
	\put(0.0,\ticsep){\makebox(0,0)[ct]{\rotatebox{0}{\scalebox{0.6}{$-0.004\hphantom{-}$}}}}%
	\put(0.5,\ticsep){\makebox(0,0)[ct]{\rotatebox{0}{\scalebox{0.6}{$0$}}}}%
	\put(1,\ticsep){\makebox(0,0)[ct]{\rotatebox{0}{\scalebox{0.6}{$0.004$}}}}%
	\put(0.5,0.23){\makebox(0,0)[ct]{\rotatebox{0}{\scalebox{0.75}{$p-p_h$}}}}%
\end{picture}%
\global\let\ticwidth\undefined%
\global\let\ticsep\undefined%
\endgroup}}
\end{tabular}
\caption{Illustration of the discrete error for the test case considered in \S\ref{prob:dirichlet}, corresponding to a single-phase Stokes problem in standard form with velocity Dirichlet boundary conditions. The error in velocity $\vu = (u,v)$ and pressure $p$ is shown in the case of a $16 \times 16$ Cartesian mesh, for $p = 3$ bicubic polynomials. Note the numerical boundary layer in pressure, which according to grid convergence analyses, does not impact the maximum norm optimal order accuracy of the velocity field.}
\label{fig:2dsquare}
\end{figure}

The next test problem is identical to the previous, but with velocity Dirichlet boundary conditions imposed on $\partial \Omega$. \Cref{fig:dirichlet-stokes} plots the measured multigrid convergence rates and shows that, with the exception of $p = 1$, excellent multigrid performance is attained, similar to that of the periodic case in \cref{fig:periodic-stokes}. When $p = 1$, we observe a stronger failure of multigrid efficiency, with $\rho$ diverging toward one as the mesh is refined. Another difference compared to the periodic case concerns the order of accuracy: numerical experiments show that the velocity attains order $p + 1$, in both the $L^2$ and maximum error norms, whereas the pressure field attains order $p+\tfrac12$ in the $L^2$ norm, and order $p$ in the maximum norm; see \iffalse Fig.~SM4 in the Supplementary Material.\fi \cref{sec:gridconv}. The order reduction in the computed pressure field, as compared to the order $p + 1$ observed in the periodic case, is due to a numerical boundary layer \iffalse (see, e.g., Fig.~SM2 in the SM); \fi (see, e.g., \cref{fig:2dsquare}); it is important to note, however, that this numerical boundary layer does not impact the optimal order accuracy of the computed velocity. See also the discussion of Cockburn \textit{et al} \cite{CockburnKanschatSchotzauSchwab2002}, wherein a priori estimates shows that order $p$ accuracy in pressure is to be expected for this LDG discretisation of the Stokes equations. 

\subsection{Stress boundary conditions}
\label{prob:stress}

\begin{figure}%[tbhp]
\centering%
\includegraphics[scale=0.91]{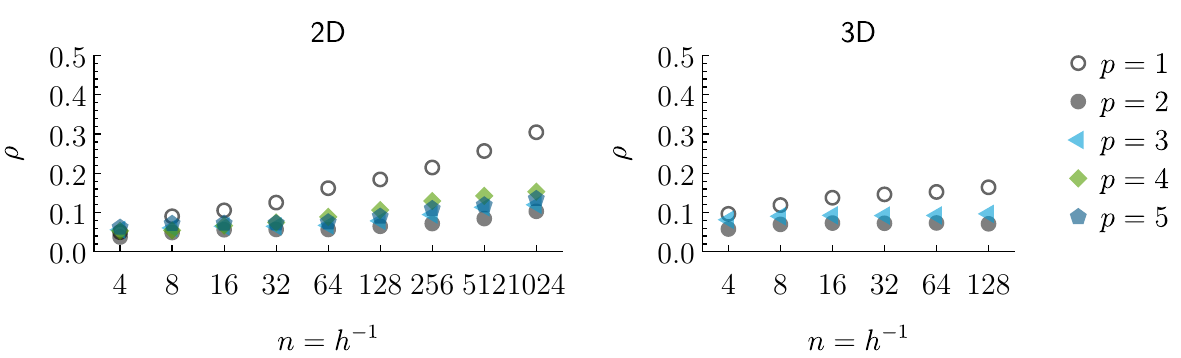}%
\caption{Measured multigrid convergence rates when solving the Stokes equations in viscous-stress form, with $\mu = 1$ and stress boundary conditions.}%
\label{fig:neumann-stress}%
\end{figure}

We next test performance of the Stokes multigrid solver when stress boundary conditions are imposed, in which case the pertinent form of the Stokes equations is the viscous-stress form ($\gamma = 1$). Similar to the previous test problems, eq.~\eqref{eq:stokes2} with $\mu = 1$ is solved on a $n \times n\,(\times\,n)$ mesh of a unit square/cube domain $\Omega = (0,1)^d$, with stress boundary conditions ${\bm \sigma} \cdot \vn = {\mathbf h}_{\partial}$ imposed on $\partial \Omega$. \Cref{fig:neumann-stress} plots the measured multigrid convergence rates. Compared to the previous two test problems (which employed the standard form of the Stokes equations), we observe a mild increase in $\rho$ for the cases in which $p > 1$, most visible in 2D. The slight increase in $\rho$ is mainly attributed to the viscous-stress form of the Stokes equations, and not solely to the imposition of stress boundary conditions. (Indeed, if one imposes Neumann-like boundary conditions for the Stokes equations in standard form, convergence rates similar to those in \cref{fig:periodic-stokes} and \cref{fig:dirichlet-stokes} are obtained.) In \cref{fig:neumann-stress}, for $p = 1$, we once again see less than ideal multigrid efficiency, though with marginal improvements compared to the case of velocity Dirichlet boundary conditions. Since the case of $p = 1$ is generally not of significant practical interest in the context of high-order accurate DG methods, we will focus on degrees $p > 1$ in the remainder of the presented results. Meanwhile, for the current test problem, numerical experiments show that the velocity field attains order $p+1$ accuracy in the maximum norm, while the pressure field attains order $p + \tfrac12$ in the $L^2$ norm and order $p$ in the maximum norm, see \cref{sec:gridconv}; as a general rule, our results indicate that whenever boundary conditions (or interfacial jump conditions) are imposed, the pressure field loses one order of accuracy near the boundary (or interface) owing to a numerical boundary layer, but this never affects the optimal order accuracy of the computed velocity field.

\subsection{Variable viscosity}
\label{prob:variableviscosity}

\begin{figure}%[tbhp]
\centering%
\includegraphics[scale=0.91]{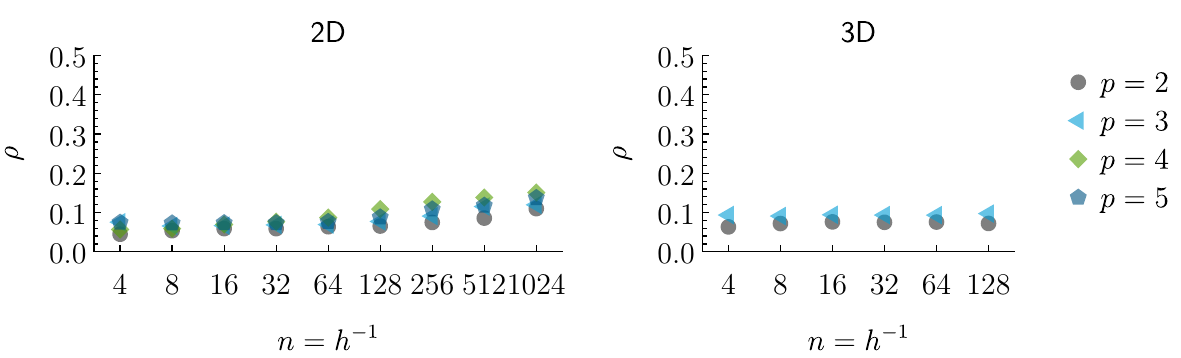}%
\caption{Measured multigrid convergence rates when solving the Stokes equations in viscous-stress form with a non-constant viscosity function $\mu: \Omega \to {\mathbb R}^+$ given by \eqref{eq:mu}, and stress boundary conditions.}%
\label{fig:variable-mu-stress}%
\end{figure}

In the next example, we consider the possibility of a non-constant viscosity function $\mu : \Omega \to {\mathbb R}^+$, varying throughout the domain. This problem serves three main purposes: (i) to test the application of the local inverse scaling by $\mu$ of the pressure penalty stabilisation parameter suggested in \eqref{eq:taueqn}; (ii) to demonstrate whether or not multigrid efficiency is impacted by variable ellipticity coefficient; and (iii) to examine the order of accuracy of the discrete solution in the variable coefficient case. Specifically, we consider the Stokes equations in viscous-stress form on the unit square/cube domain $\Omega = (0,1)^d$, where $\mu : \Omega \to {\mathbb R}^+$ is given by
\begin{equation} \label{eq:mu} \mu = \begin{cases} 1 + \tfrac12 \sin 4\pi x \sin 4 \pi y & \text{in 2D,} \\ 1 + \tfrac12 \sin 4\pi x \sin 4\pi y \sin 4 \pi z & \text{in 3D,} \end{cases} \end{equation}
together with stress boundary conditions on $\partial\Omega$. Measured multigrid convergence rates, plotted in \cref{fig:variable-mu-stress}, show very similar behaviour to the constant-viscosity test problem of the previous example (\cref{fig:neumann-stress}). Meanwhile, numerical experiments examining the order of accuracy (see \cref{sec:gridconv}) confirm that the velocity field attains order $p+1$ in the maximum norm, while pressure attains order $p+\tfrac12$ in the $L^2$ norm and order $p$ in the maximum norm.

\subsection{Curved domain geometry}
\label{prob:sphere}

\begin{figure}%[tbhp]
\centering\sffamily\footnotesize%
\begin{minipage}{0.3\textwidth}\centering%
\includegraphics[width=\textwidth]{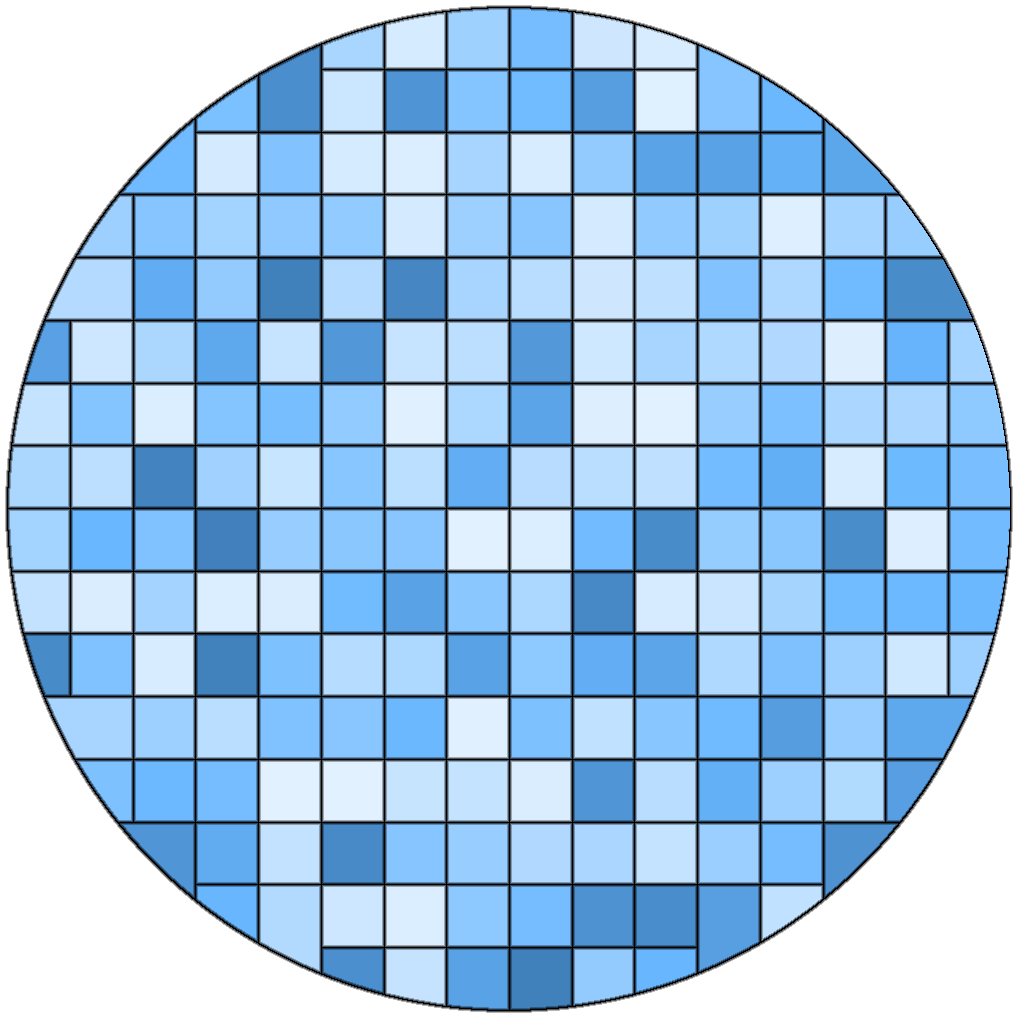}\\[0.5em]
(a) Circular domain
\end{minipage}
\qquad\qquad
\begin{minipage}{0.3\textwidth}\centering%
\includegraphics[width=\textwidth]{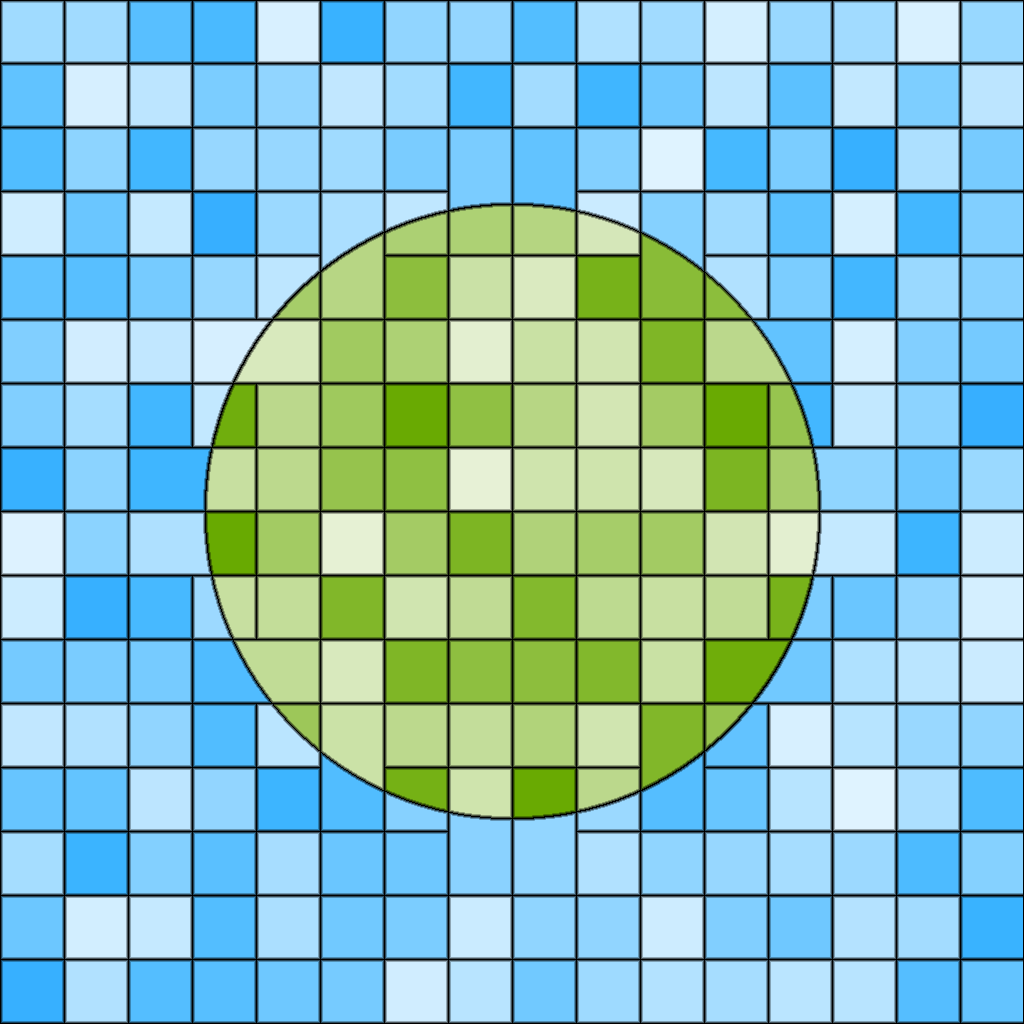}\\[0.5em]
(b) Circular interface
\end{minipage}%
\caption{Examples of implicitly defined meshes generated with a background $16 \times 16$ Cartesian grid. (a) Circular domain, used in the single-phase Stokes problem examined in \cref{fig:sphere}. (b) Square domain with an embedded circular interface, used in the two-phase time-dependent Stokes problem considered later in this article, in \cref{fig:water-bubble}.}%
\label{fig:circle}%
\end{figure}

In the last single-phase example of this section, we consider a Stokes problem in a curved domain. Here, and in all others involving curved geometry, we make use of implicitly defined meshes, building on the implicit mesh DG framework developed in prior work \cite{ImplicitMeshPartOne,ImplicitMeshPartTwo}. Briefly, an implicitly defined mesh uses one or more level set functions, describing the domain geometry and/or embedded interfaces, to cut through the cells of a background quadtree or octree; tiny cut cells are then merged with neighbouring cells to create a mesh such that elements adjacent to curved geometry have their shape defined implicitly by the level set functions. In particular, the resulting mesh is interface- and boundary-conforming, thereby sharply representing the implicitly-defined geometry. To use an implicitly defined mesh with a DG method, the main task is to compute quadrature schemes for the elements and faces whose geometry is implicitly defined; these are computed using the high-order accurate algorithms detailed in refs.~\cite{HighorderImplicitQuad,algoim}, and then used when computing mass matrices, discrete gradient operators, $L^2$ projections, etc. For details on the implicit mesh DG framework, see refs.~\cite{ImplicitMeshPartOne,ImplicitMeshPartTwo}, and for illustrations demonstrating the associated (implicitly-formed) mesh hierarchy underlying the multigrid method, see ref.~\cite{dgmg}.

\begin{figure}%[tbhp]
\centering%
\includegraphics[scale=0.91]{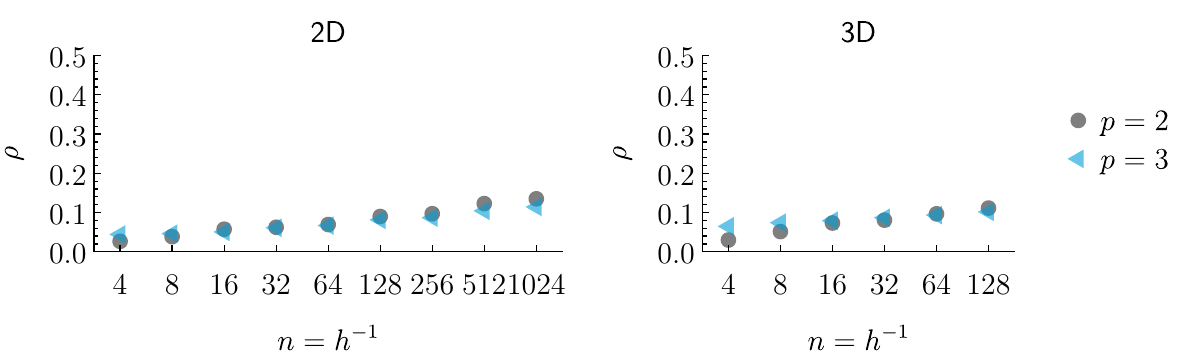}%
\caption{Measured multigrid convergence rates when solving the Stokes equations in standard form, with $\mu = 1$ in a unit diameter spherical domain using implicitly defined meshes together with velocity Dirichlet boundary conditions. Here, $n$ denotes the number of cells in the background uniform Cartesian $n\times n\,(\,\times n)$ grid used to build the corresponding implicitly defined mesh; see \cref{fig:circle}(a).}%
\label{fig:sphere}%
\end{figure}

In the current test problem, we consider solving the Stokes equations, in standard form, on a unit diameter circle or sphere, $\Omega = \{ x \in {\mathbb R}^d : |x| < \tfrac12 \}$, with $\mu = 1$ and velocity Dirichlet boundary conditions. An example of the corresponding implicitly defined mesh generated by a background $16 \times 16$ Cartesian grid is shown in \cref{fig:circle}(a). Measured multigrid convergence rates for this test case are shown in \cref{fig:sphere} for polynomial degrees $p = 2$ and $p = 3$. Overall, efficient multigrid convergence rates are witnessed, however once $n \geq 256$ we observe that $\rho$ increases in value compared to previous test cases; we attribute this to the worsening conditioning of the discrete Stokes operator. In particular, inspection of numerical results shows that the conditioning of ${\mathcal A}$ exceeds $10^8$ in these cases. Since the V-cycle is built by coarsening the components of the fine mesh operator ${\mathcal A}$, and since $\rho$ is measured according to how many GMRES iterations it takes to reduce the (preconditioned) residual by a factor $10^8$, it follows that $\rho$ may be impacted by this conditioning.  
Regarding the grid convergence analysis (see \cref{sec:gridconv}), numerical results once again confirm that the velocity field attains order $p+1$ in the maximum norm, whereas the pressure attains order $p+\tfrac12$ in the $L^2$ norm and order $p$ in the maximum norm.

\subsection{Multi-phase Stokes equations with interfacial jump conditions}
\label{prob:square}

In this last example for the steady-state Stokes equations, we consider a challenging situation in which the viscosity coefficient $\mu$ exhibits jumps several orders in magnitude across an embedded interface. In particular, let $\Omega = (0,1)^d$ be divided into an interior square or cubic phase, $\Omega_1 = (\tfrac14,\tfrac34)^d$, together with an exterior phase, $\Omega_2 = \Omega \setminus \overline{\Omega_1}$. Four different ratios of viscosity jump across the interface are considered: in every case, the exterior phase has unit viscosity, $\mu_2 = 1$, whereas the interior phase will have one of four values, $\mu_1 \in \{10^{-6}, 10^{-3}, 10^3, 10^6\}$. We consider the multi-phase viscous-stress form of the Stokes equations, \eqref{eq:govern1}--\eqref{eq:govern3} with $\gamma = 1$, for which interfacial jump conditions in velocity and the stress tensor are imposed on $\Gamma = \partial \Omega_1 \cup \partial \Omega_2$.
To tackle this problem, we apply two distinct but complementary strategies: viscosity-upwinded weighted fluxes on interfacial mesh faces, and a diagonal scaling to improve the conditioning of the Stokes operator; these are discussed next.

Prior work on LDG methods for elliptic interface problems \cite{fluxx} shows that, in order to obtain ideal multigrid efficiency and solution accuracy, one should apply a biasing strategy for the LDG numerical fluxes on interfacial faces, wherein $\hat{\vu}$ and $\hat{\vsigma}$ (see \iffalse the Supplementary Material \fi \cref{sec:ldg}) are biased toward one phase or the other, depending on the local viscosity coefficients. Here, we build on these ideas and extend it the multi-phase Stokes case.

The strategy of biasing can be intuitively motivated as follows. Suppose that the interior phase has a vastly smaller viscosity than the exterior phase, i.e., $\mu_1 = \epsilon$ with $0 < \epsilon \ll 1$. Suppose also that the velocity $\vu$ and its gradient is unit order in magnitude near the interface. From a physical standpoint, the pressure scales as $p \sim \mu U/L$ (where $U$ is a typical velocity scale and $L$ is a typical length scale), and so, suppose also that pressure and its gradient is unit order magnitude in the exterior phase, and scales with $\epsilon$ in the interior phase. Then, in the limit of vanishingly small $\epsilon$, the stress jump condition in \eqref{eq:govern2} approximately reduces to phase $\Omega_2$ having a stress-like boundary condition, $\mu_2 (\nabla \mathbf u + \nabla \mathbf u^T) \mathbf n - p \mathbf n \approx {\mathbf h}_{12}$ on $\Gamma$. Except for the modes associated with the trivial kernel of the corresponding Stokes operator, this is enough to determine the solution $(\mathbf u_2,p_2)$ in $\Omega_2$ and thus the phase $\Omega_2$ (nearly) decouples from phase $\Omega_1$. Once $(\mathbf u_2,p_2)$ is found, the remaining jump condition in \eqref{eq:govern2} essentially reduces to a Dirichlet boundary condition for the Stokes problem in phase $\Omega_1$, i.e., $\mathbf u_1|_\Gamma \approx  {\mathbf u}_2|_\Gamma - {\mathbf g}_{12}$, which is enough to determine $(\mathbf u_1,p_1)$ up to the modes associated with its kernel. Therefore, in the limit of vanishingly small $\epsilon$, the highly viscous phase essentially ``sees'' a stress boundary condition on $\Gamma$ (whose data is nearly independent of the other phase), and the nearly inviscid phase ``sees'' a Dirichlet boundary condition on $\Gamma$ (whose data depends on the solution across the interface). It follows, therefore, that it may be advantageous to bias the numerical fluxes of the LDG formulation in the same way, to reflect the physical nature of the interfacial jump conditions in \eqref{eq:govern2}. 

For elliptic interface problems, a common approach is to bias the fluxes according to a kind of harmonic weighting; in ref.~\cite{fluxx}, a stronger kind of biasing is advocated, denoted viscosity-upwinded weighting. Here, we adopt the same strategy and refer the reader to ref.~\cite{fluxx} for discussion, much of which is directly analogous to the multi-phase Stokes case, and to \iffalse the Supplementary Material \fi \cref{sec:ldg} with details on how the multi-valued numerical flux functions $\hat\vu$ and $\hat\vsigma$ are chosen for interfacial mesh faces. 

In addition to the application of viscosity-upwinded numerical fluxes, one other numerical technique is employed to improve the conditioning of the multi-phase Stokes equations. For a single-phase, constant-coefficient Stokes problem with viscosity $\mu$, i.e., $\{-\mu \nabla^2 \mathbf u + \nabla p = \mathbf f, -\nabla \cdot \mathbf u = f\}$, the largest positive eigenvalue of the discretised Stokes operator scales as $\mu/h^2$, whereas all negative eigenvalues scale inversely-proportional to $\mu$ (and independently of $h$). Thus, unlike a Poisson problem, whose conditioning is independent of an arbitrary multiple of its ellipticity coefficient, the conditioning of the Stokes problem worsens quadratically in $\mu$ as $\mu$ is made arbitrarily large. In essence, the two operators of the momentum equations, $\mu \nabla^2 \mathbf u$ and $\nabla p$, are not on equal footing if $\mu$ is large and $\mathbf u$ and $p$ are treated as independent variables. However, as noted earlier, from a physical point of view, the pressure scales as $p \sim \mu U/L$, and so the magnitude of $p$ depends on $\mu$. This apparent ill-conditioning of the Stokes operator can be easily remedied by rescaling one of the variables $\mathbf u$ or $p$ by $\mu$ or $1/\mu$, respectively, effectively recasting the Stokes equations into a unit-viscosity form, $-\nabla^2 \tilde{\mathbf u} + \nabla \tilde p = \tilde{\mathbf f}$. 

This issue of ill-conditioning is exaggerated for the multi-phase Stokes equations, in which the largest positive eigenvalue scales as $\max_i \mu_i /h^2$. However, as motivated by the single-phase case, a simple rescaling of the solution variables $\mathbf u$ and $p$ can be used to mitigate the issue (see also refs.~\cite{FuruichiMayTackley2011,AksoyluUnlu2014,BorzacchielloLericheBlottiereGuillet2017}). In this work, we achieve this through a diagonal pre- and post-scaling of the Stokes operator, i.e., replace the Stokes operator in \eqref{eq:blockform} with
\[ \begin{pmatrix} \alpha & 0 \\ 0 & \beta \end{pmatrix} \begin{pmatrix} A & M {\mathcal G} \\ {\mathcal G}^\trans M & -E \end{pmatrix} \begin{pmatrix} \alpha & 0 \\ 0 & \beta \end{pmatrix} \]
where $\alpha$ and $\beta$ are diagonal matrices whose entries equal $1/\sqrt{\mu}$ and $\sqrt{\mu}$, respectively. (More precisely, for every velocity field $\mathbf u$ and pressure field $p$, we have that $\smash{(\alpha \mathbf u)|_E = \mu_E^{-1/2} {\mathbf u}|_E}$ and $\smash{(\beta p)|_E = \mu_E^{1/2} p|_E}$ for every element $E \in {\mathcal E}$, where $\mu_E$ is the viscosity on element $E$.) This pre- and post-scaling is built into the multigrid schemes by replacing the Stokes operator with the scaled version on all levels of the hierarchy. Note that it is straightforward to adapt the original linear system to the scaled approach: instead of solving ${\mathcal A}x = b$, one instead solves for $\tilde {\mathcal A} \tilde x = \tilde b$, where $\tilde {\mathcal A}$ is the scaled Stokes operator and $\tilde b = \operatorname{diag}(\alpha,\beta) b$, whereupon solving for $\tilde x$ gives the original unscaled solution as $x = \operatorname{diag}(\alpha,\beta) \tilde x$.

\begin{figure}[t]%[tbhp]
\centering\sffamily\footnotesize%
(a) $\tfrac{\mu_1}{\mu_2} = 10^{-3}$\\[0em]%
\includegraphics[scale=0.91]{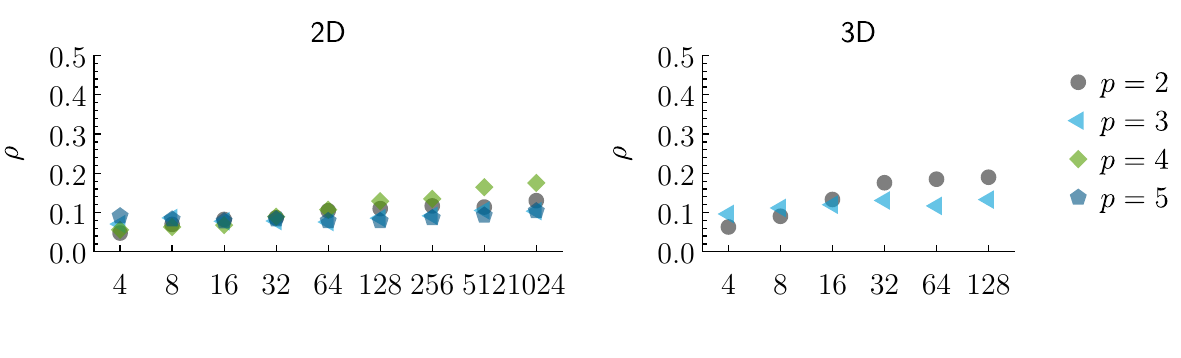}\\[-1em]%
(b) $\tfrac{\mu_1}{\mu_2} = 10^{+3}$\\[-1em]%
\includegraphics[scale=0.91]{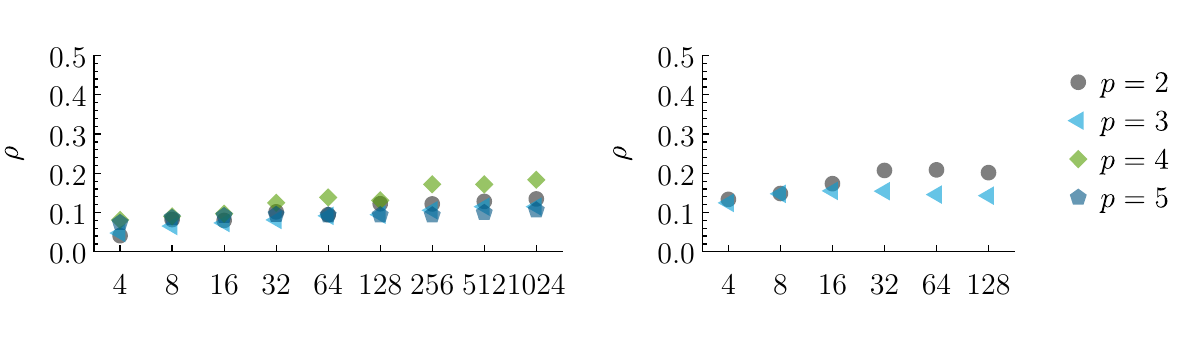}\\[-1em]%
(c) $\tfrac{\mu_1}{\mu_2} = 10^{-6}$\\[-1em]%
\includegraphics[scale=0.91]{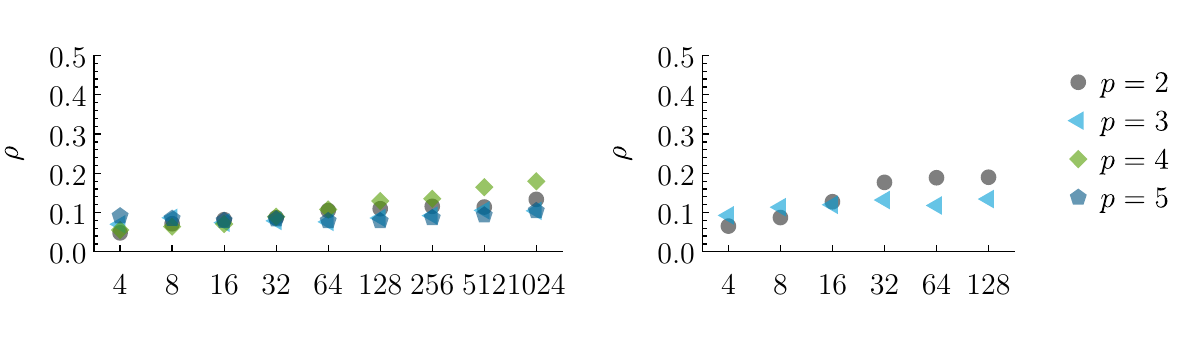}\\[-1em]%
(d) $\tfrac{\mu_1}{\mu_2} = 10^{+6}$\\[-1em]%
\includegraphics[scale=0.91]{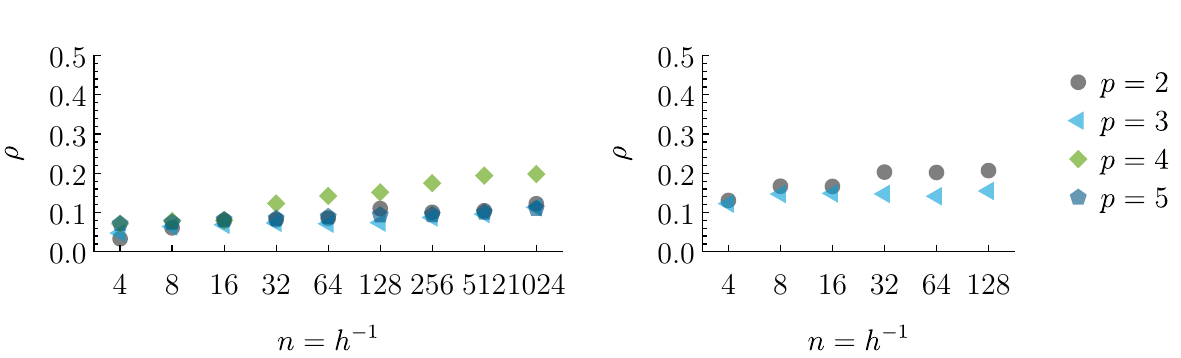}%
\caption{Measured multigrid convergence rates when solving the multi-phase Stokes equations in viscous-stress form, wherein $\Omega_1 = (\tfrac14,\tfrac34)^d$ has viscosity $\mu_1$ and $\Omega_2 = (0,1)^d \setminus \overline{\Omega_1}$ has viscosity $\mu_2$, with the viscosity ratio as indicated, and periodic boundary conditions. In each case, 2D results are plotted in the left column, and 3D results in the right column.}%
\label{fig:square}%
\end{figure}

With the application of viscosity-upwinded numerical fluxes and the simple diagonal scaling to improve conditioning, \cref{fig:square} plots the measured multigrid convergence rates when solving the multi-phase Stokes equation in viscous-stress form with periodic boundary conditions, for the four different viscosity ratios considered, $\tfrac{\mu_1}{\mu_2} \in \{10^{-6}, 10^{-3}, 10^3, 10^6\}$. Compared to previous test problems, we observe marginally slower convergence rates, representative of the challenging Stokes problems at hand. Nevertheless, good convergence rates are obtained across significant viscosity ratios, i.e., typically 7--10 iterations of multigrid-preconditioned GMRES for a $10^8$ reduction in the residual norm.  Numerical experiments examining the order of accuracy (see \iffalse the Supplementary Material \fi \cref{sec:gridconv}) show that the velocity attains order $p+1$ in the maximum norm, while pressure attains order $p+\tfrac12$ in the $L^2$ norm and order $p$ in the maximum norm. As in other test problems, a numerical boundary layer exists in the pressure field: in the present case, the layer is adjacent to the interface, but, as in other test problems, does not impact the optimal order accuracy of the velocity.

\section{Multigrid Efficiency for the Time-Dependent Stokes Equations}

So far, we have focused on the efficacy of the multigrid solver when applied to the steady-state Stokes equations. The time-dependent Stokes equations, however, may pose additional challenges, owing to the competing effects of the temporal derivative and viscous stress operator \cite{BenziGolubLiesen2005}. To illustrate, consider the time-dependent Stokes equations in the following form:
\begin{equation} \label{eq:timedep1} \left. \begin{aligned} \frac{\rho_i}{\delta} \vu -\nabla \cdot \bigl(\mu_i (\nabla \vu + \gamma \nabla \vu^\trans) \bigr) + \nabla p &= \mathbf f \\ -\nabla \cdot \mathbf \vu &= f \end{aligned} \right\} \text{ in $\Omega_i$,} \end{equation}
where $\rho_i$ is a phase-dependent density and $\delta$ is a parameter proportional to the time step $\Delta t$ of a temporal integration method. Its LDG discretisation is a straightforward modification to the corresponding steady-state Stokes equations, and reads
\begin{equation} \label{eq:timedep2} \begin{pmatrix} \frac{1}{\delta} M_\rho + A & M {\mathcal G} \\ {\mathcal G}^\trans M & -E \end{pmatrix} \begin{pmatrix} \vu_h \\ p_h \end{pmatrix} = \begin{pmatrix} {\mathbf b}_{\vu} \\ b_p \end{pmatrix}, \end{equation}
where $M_\rho$ is a $\rho$-weighted mass matrix, with all other operators unchanged.

When the viscous operator dominates (i.e., a small effective Reynolds number such that $\mu \delta/\rho$ is sufficiently large, depending on mesh resolution), the dominant operator is the steady-state Stokes equations, with a small $\rho/\delta$-weighted identity shift added to the viscous operator; in this case, one may expect a good steady-state Stokes solver to be effective. On the other hand, in the case when viscous effects are weak (i.e., a large effective Reynolds number such that $\mu \delta/\rho$ is small or the mesh is unable to resolve viscous effects), then \eqref{eq:timedep1} essentially reduces to a Helmholtz-Hodge projection problem (having strong connections to Chorin's projection method for solving the incompressible Navier-Stokes equations \cite{Chorin,GuermondMinevShen,sdpc}). In the latter case, a solver designed specifically for the steady-state Stokes equations may deteriorate. 

In this work, a simple strategy is employed to automatically account for these two competing effects, resulting in fast multigrid solvers across a full range of Reynolds numbers. In essence, the strategy chooses the pressure penalty stabilisation parameter $\tau_p$ in \eqref{eq:eop} according to the expected scaling of the maximal eigenvalue of the operator $\tfrac{\rho}{\delta} {\mathbb I} - \nabla \cdot (\mu(\nabla + \gamma \nabla^\trans))$. 
In the steady-state case, e.g., $\rho = 0$, the penalty parameter should scale as $\tau_p \sim h/\mu$. However, in the degenerate time-dependent Stokes case, in which $\mu = 0$, \eqref{eq:timedep2} reduces to an LDG method, written in flux form, for computing the solution of a Poisson problem with operator $\smash{\nabla \cdot (\tfrac{\delta}{\rho} \nabla p)}$; the appropriate scaling of the pressure penalty stabilisation parameter in this case is then $\tau_p \sim \delta/(\rho h)$ \cite{ldg,fluxx}. Both of these scaling statements may be summarised as follows: the penalty parameter should scale such that $\tau \sim (h \Lambda)^{-1}$, where $\Lambda$ is the maximal eigenvalue of either the discretised operator $-\nabla \cdot (\mu (\nabla + \gamma \nabla^\trans))$ (in the former case) or $\tfrac{\rho}{\delta} {\mathbb I}$ (in the latter case). This leads to a simple, but effective, idea to treat the general case: let $\tau \sim (h \Lambda)^{-1}$ where $\Lambda$ is the sum of the expected scalings of the maximal eigenvalues of the two operators. Using this idea, the definition of the pressure penalty parameter in \eqref{eq:eop}--\eqref{eq:taueqn} is replaced with the following relation:
\begin{equation} \label{eq:taugeneral} \tau_p = \Bigl( \frac{h \rho}{\tau_0 \delta} + \frac{\mu}{\tau h} \Bigr)^{-1}. \end{equation}
Note that \eqref{eq:taugeneral} reduces to the correct penalty parameter for the steady-state Stokes case when $\rho = 0$, i.e., $\tau_p = \tau h/\mu$, where $\tau$ is the multigrid-optimal parameter given in \cref{tab:optimaltau}; in the other extreme, when $\mu = 0$ (or $\delta$ is vanishingly small), \eqref{eq:taugeneral} reduces to $\tau = \tau_0 \delta / (h \rho)$, appropriate for a scalar Poisson problem with ellipticity coefficient $\tfrac{\delta}{\rho}$, where $\tau_0$ is a user-defined constant prefactor. In the general situation, i.e., between these two extremes, it is important to note the scaling of \eqref{eq:taugeneral} may change across the multigrid hierarchy---on a very fine mesh, $h$ is small so that the second term may dominate; on a coarse mesh, however, $h$ is large and the first term may dominate. This should be taken into account when building the multigrid method---in a purely-geometric approach, wherein $\mathcal A_h$ is explicitly built on every level, \eqref{eq:taugeneral} can be utilised directly; for the operator coarsening strategy, 
a simple modification to existing schemes is discussed in \iffalse the Supplementary Material. \fi \cref{sec:op}. 

The penalisation choice in \eqref{eq:taugeneral} blends across the two extremes of the time-dependent Stokes equations: a steady-state Stokes problem at one extreme, and a pure Poisson problem (written in flux form) at the other. Across a full range of Reynolds numbers, convergence results confirm optimal order of accuracy and that this choice of $\tau$ neither saturates the discretisation error nor underpenalises. As in the case of steady-state Stokes problems, we found that a V-cycle with three pre- and post-smoothing steps (used throughout this work) resulted in approximately optimal convergence speed, independent of the effective Reynolds number; occasionally a different combination of pre- and post-smoothing steps (at least two and at most four) can be slightly faster, but as mentioned earlier, the precise optimal value depends on a host of implementation, hardware, and problem-dependent factors.

\subsection{\hspace{-0.04em}Single-phase time-dependent Stokes equations}
\label{prob:timedep1}

\begin{figure}%[tbhp]
\centering\sffamily\footnotesize%
(a) $\mu = 10^{-2}$ corresponding to $\textsf{Re} \approx 100$\\[0.2em]
\includegraphics[scale=0.91]{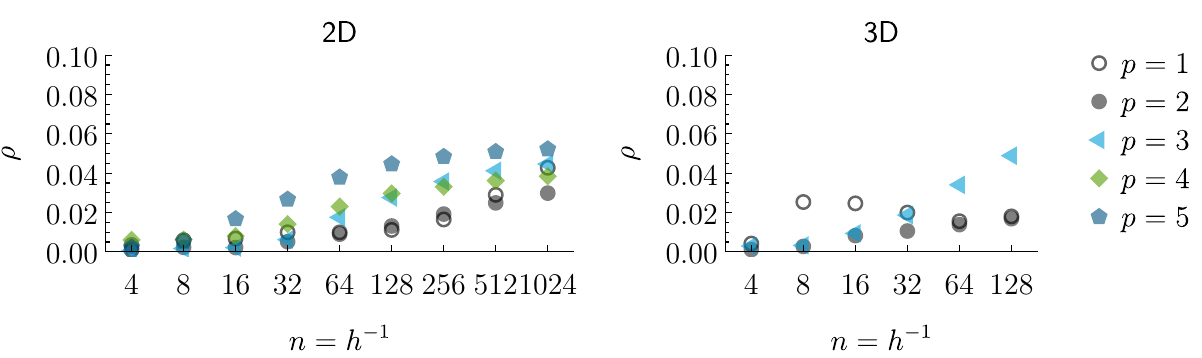}\\[1em]%
(b) $\mu = 10^{-4}$ corresponding to $\textsf{Re} \approx 10,\!000$\\[0.2em]
\includegraphics[scale=0.91]{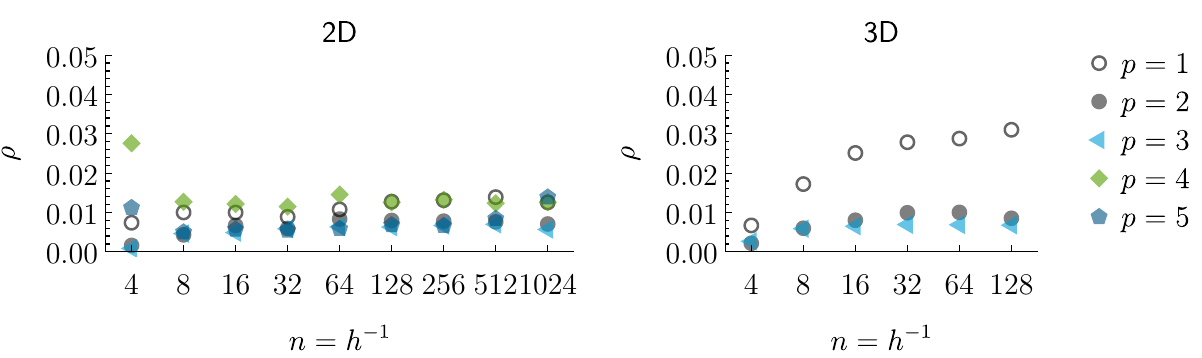}%
\caption{Measured multigrid convergence rates when solving the time-dependent, single-phase Stokes equations in standard form with $\mu$ as indicated and $\delta = 0.1 h$, together with Neumann boundary conditions.}%
\label{fig:time-single-phase-24}%
\end{figure}

To demonstrate multigrid performance on time-dependent single-phase Stokes problems, we consider two effective Reynolds numbers $\textsf{Re} = \rho U L/\mu$, one in which $\textsf{Re} \approx 100$, representing a viscous-dominated case (but where the time-derivative operator nevertheless influences performance characteristics) and the other with $\textsf{Re} \approx 10,\!000$ (wherein the time-derivative operator definitively dominates). In both cases, the velocity and length scales are unitary, $U = 1$, $L = 1$; density is set to $\rho = 1$, while viscosity satisfies $\mu = 10^{-2}$ in the former case and $\mu = 10^{-4}$ in the latter. In addition, we set $\delta = 0.1 h$, where $h$ is the element size on the finest-level mesh, representing a typical scenario of applying the time-dependent Stokes equations in a temporal integration method with CFL about $0.1$. As before, multigrid efficiency is quantified through the average convergence rate $\rho$; \cref{fig:time-single-phase-24} plots the results (note the magnified vertical axis).\footnote{In the results of this section, we have reincluded the case of $p = 1$, mainly to serve as a point of interest: specifically, for steady-state Stokes problems, we noted in \cref{sec:timeindep} that sub-optimal multigrid efficiency may occur when $p = 1$; for unsteady Stokes problems, and depending on the effective Reynolds number, ideal multigrid convergence can be restored in this case, as seen in \cref{fig:time-single-phase-24,fig:water-bubble}.} In \cref{fig:time-single-phase-24}(a), we observe an upwards trend in $\rho$ as the mesh is refined; this corresponds to the fact that as the mesh is refined, eventually the viscous operator will definitively dominate and multigrid convergence rates similar to the steady-state Stokes equations will be attained. In the weakly viscous case with $\textsf{Re} \approx 10,\!000$, \cref{fig:time-single-phase-24}(b) shows exceptionally fast multigrid convergence rates, with $\rho \approx 0.01$, corresponding to needing only 4 GMRES iterations to achieve a factor $10^8$ reduction in the initial residual. In the context of solving the incompressible Navier-Stokes equations, these results suggest that a fast non-stationary Stokes solver may outperform a fast projection method solver; further remarks on this topic are given in the conclusions. Meanwhile, grid convergence experiments examining the order of accuracy (see \cref{sec:gridconv}) show a departure from the typical results seen elsewhere in this work. In the weakly-viscous case with $\textsf{Re} \approx 10,\!000$, the velocity is order $p+1$ in the maximum norm, and so is the pressure, despite the presence of boundary conditions. This apparent ``super convergence'' in pressure is attributed to the property that, for very large Reynolds numbers, the time-dependent Stokes equations nearly reduce to a Helmholtz-Hodge projection, where one may naturally expect to attain optimal order accuracy in the pressure field, see, e.g., ref.~\cite{ImplicitMeshPartOne}. On the other hand, in the viscous-dominated case with $\textsf{Re} \approx 100$, the pressure reduces to order $p+\tfrac12$ in the $L^2$ norm and order $p$ in the maximum norm, while the velocity maintains order $p+1$ in the maximum norm. This is perhaps expected, based on the order of accuracy results reported elsewhere in this work, together with the property that the viscous-dominated case essentially represents a mildly perturbed stationary Stokes problem.

\subsection{Multi-phase time-dependent Stokes equations}
\label{prob:timedep2}

In our last two examples, we combine the challenging aspects of a multi-phase Stokes problem together with non-stationary effects and consider a problem in which both the density $\rho$ and viscosity $\mu$ have discontinuities several orders in magnitude across an embedded interface. The specific parameters considered correspond to two scenarios: a water bubble surrounded by gas, and a gas bubble surrounded by water. Specifically, $\rho_{\text{water}} = 1$, $\rho_{\text{gas}} = 0.001$, $\mu_{\text{water}} = 1$, and $\mu_{\text{gas}} = 0.0002$ (approximately accurate values for water and air at ambient temperature in CGS units). The radius of the bubble is $0.3$, centred in a unit square/cubic domain, $\Omega = (0,1)^d$; equations \eqref{eq:timedep1} with $\gamma = 1$ and \eqref{eq:govern2}--\eqref{eq:govern3} are solved with velocity Dirichlet boundary conditions and stress jump conditions across the gas-water interface. As in the test case on curved domain geometry, a level set function describing the interface geometry is used together with an implicit mesh DG framework to create a semi-unstructured, interface-conforming mesh; an example is shown in \cref{fig:circle}(b). As before, the time step is set equal to $\delta = 0.1 h$, where $h$ is the typical element size on the finest mesh, representing the application of a time stepping method with CFL number about $0.1$.

\begin{figure}%[tbhp]
\centering\sffamily\footnotesize%
(a) Water bubble surrounded by gas\\[0.2em]
\includegraphics[scale=0.91]{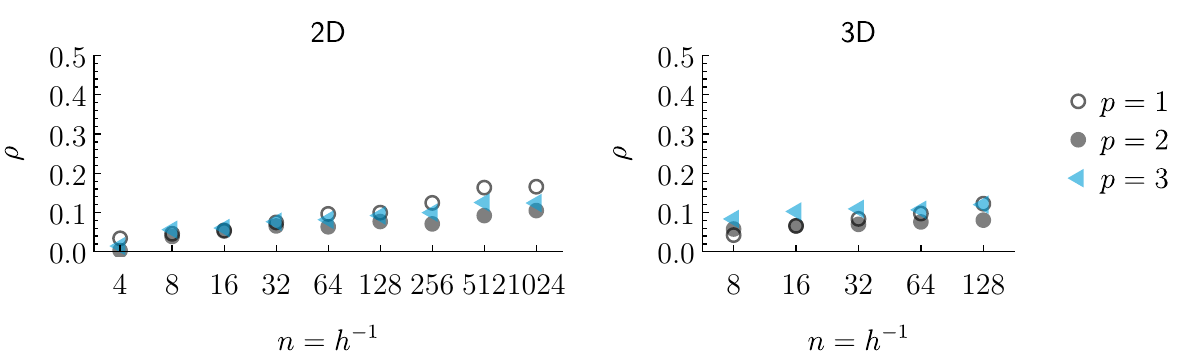}\\[1em]%
(b) Gas bubble surrounded by water\\[0.2em]
\includegraphics[scale=0.91]{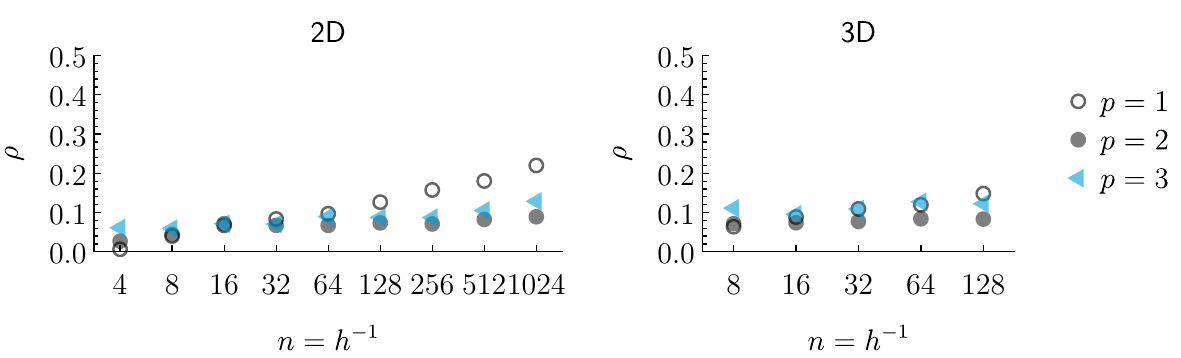}%
\caption{Measured multigrid convergence rates when solving the time-dependent multi-phase Stokes equations in viscous-stress form, with $\rho_{\text{water}} = 1$, $\rho_{\text{gas}} = 0.001$, $\mu_{\text{water}} = 1$, $\mu_{\text{gas}} = 0.0002$, and $\delta = 0.1 h$, together with velocity Dirichlet boundary conditions. Here, $n$ denotes the number of cells in the background uniform Cartesian $n\times n\,(\,\times n)$ grid used to build the corresponding implicitly defined mesh; see \cref{fig:circle}(b).}%
\label{fig:water-bubble}%
\end{figure}

This example combines three distinct but complementary strategies developed in prior examples: (i) viscosity-upwinded numerical fluxes, to robustly and accurately handle the large jump in viscosity across the interface; (ii) diagonal pre- and post-scaling of the Stokes operator, to remove the unnecessary ill-conditioning caused by viscosity coefficients differing several orders in magnitude; and (iii) a pressure penalty stabilisation parameter controlled by \eqref{eq:taugeneral}, to automatically adjust penalisation behaviour between the two extremes of the time-dependent Stokes equations. As such, this example serves to demonstrate a variety of subtleties in efficiently solving the multi-phase time-dependent Stokes equations, but as shown in the results plotted in \cref{fig:water-bubble}, combined together, one can attain highly efficient multigrid solvers. For the particular scales chosen in this problem, the water-bubble and gas-bubble problems correspond to a unit-order Reynolds number. As such, grid convergence analyses confirms the expectation that the velocity attains order $p+1$ in the maximum norm, and pressure attains order $p+\tfrac12$ in the $L^2$ norm and order $p$ in the maximum norm. 

\section{Concluding Remarks}

In this paper, we devised efficient geometric multigrid solvers for the Stokes equations discretised by local discontinuous Galerkin (LDG) methods. The approach follows standard geometric multigrid concepts, utilising a V-cycle and a simple block Gauss-Seidel relaxation method free from under-relaxation parameters. With a suitably-chosen pressure penalty stabilisation parameter, a wide array of tests showed that the Stokes multigrid solver can match the speed of classical geometric multigrid methods for Poisson problems \cite{dgmg,fluxx}. For example, typical convergence rates ranged from $\rho \approx 0.05$ to $0.1$, corresponding to needing about six to eight iterations for a $10^8$ reduction of the residual; fewer iterations are required for unsteady time-dependent Stokes equations. 
To implement the multigrid algorithm, one possibility is to explicitly build the mesh and associated LDG operators on every level of the hierarchy; an alternative method, not requiring the explicit formation of coarse meshes, can be implemented based on the operator coarsening algorithms detailed in \iffalse the Supplementary Material.\fi \cref{sec:op}. In addition, we also extended the LDG methods of ref.~\cite{CockburnKanschatSchotzauSchwab2002} to variable-viscosity and multi-phase problems exhibiting interfacial stress jump conditions; across all test problems, grid convergence analyses demonstrated order $p+1$ accuracy in the maximum norm for the computed velocity field, and at least order $p$ accuracy in the maximum norm for pressure.

Key findings of this work include the following aspects:
\begin{itemize}
\item Multigrid efficiency depends on an appropriate choice of the pressure penalty stabilisation parameter $\tau_p$. In general, for the steady-state Stokes equations, $\tau_p$ scales linearly with the local element size $h$ and inversely-proportional to the local viscosity, leaving the end-user to define the prefactor $\tau$ in the formula $\tau_p = \tau\, h/\mu$. \Cref{tab:optimaltau} provides values of $\tau$ resulting in optimal multigrid convergence rates, showing that $\tau$ depends on the polynomial degree of the DG space, the spatial dimension, and whether the standard form or viscous-stress form of the Stokes equations are being solved. \iffalse Table SM1 in the Supplementary Material \fi \Cref{tab:tauwindow} shows that an associated range of values exist for $\tau$ within which multigrid performance is very close to optimal. (See also Kanschat~\cite{Kanschat2005}, in which the performance of an approximate Schur form block preconditioner for LDG was also noted to depend crucially on the pressure penalty parameter, although the scaling of $\tau$ determined therein is very different to the results of the present work.)

\item In the case of the time-dependent Stokes equations, \eqref{eq:taugeneral} implements an effective strategy to automatically adjust the penalty parameter based on the competing effects of the viscous and time-derivative operators.
\item For multi-phase problems, in which the viscosity may exhibit jumps several orders in magnitude across an embedded interface, one should bias the numerical fluxes (here implemented with viscosity-upwinded weighting) to reflect the physical nature of the interfacial stress and velocity jump conditions. In addition, it can be advantageous to employ a diagonal preconditioning of the Stokes operator to correct the potential adverse affects large viscosities may have on the spectral characteristics of the saddle point problem.  
\end{itemize}

A formula for the optimal pressure penalty parameter, $\tau$, is not known at this time. In this work, a simple one-dimensional parameter sweep was used to find $\tau$ for particular choices of $p$, but a theoretical result would be ideal. Here, it may be possible to use analytical tools such as convergence criteria for block Gauss-Seidel methods \cite{ElsnerMhermann1991}, local mode analyses \cite{Sivaloganathan1991,GasparNotayOosterleeRodrigo2014} or local Fourier analyses \cite{HeMacLachlan2017,FarrellHeMacLachlan2019} to help determine a formula. On a similar note, it may also be possible to derive analytical proofs of the convergence of the overall multigrid method; see also, e.g., \cite{SchoberlZulehner2003,Manservisi2006,DrzisgaJohnRudeWohlmuthZulehner2018,KanschatMao2015,HongKrausXuZikatanov2016}.

Other areas of study include the following. In this work, we mainly considered structured meshes (such as Cartesian grids) and semi-unstructured, non-conforming implicitly defined meshes resulting from cell merging procedures. Efficacy of the multigrid Stokes solver on fully unstructured meshes is also worthy of examination--for highly anisotropic meshes, one may need to group elements into clusters for the block Gauss-Seidel method to be effective \cite{ThompsonFerziger1989,OosterleeWesseling1993,WobkerTurek2009}. Another possibility is to accelerate the Gauss-Seidel method through low-degree Chebyshev iterations, as was noted by Farrel \textit{et al} \cite{FarrellHeMacLachlan2019} for Vanka-type smoothers. Meanwhile, a wide variety of work on developing multigrid methods for the Stokes equations report that W-cycles can be more effective than V-cycles, or even variable V-cycles which change the smoothing counts between levels of the hierarchy; see, e.g., \cite{GasparNotayOosterleeRodrigo2014,DrzisgaJohnRudeWohlmuthZulehner2018,FarrellHeMacLachlan2019,KanschatMao2015,AdlerBensonMacLachlan2016}. Although a V-cycle with fixed pre-smoothing and post-smoothing steps was found highly effective in this work, these alternative strategies could prove useful for different applications. In the presented results, we also mentioned that the block Gauss-Seidel method is less effective when $p = 1$, i.e., for a bilinear or trilinear LDG discretisation of the Stokes equations; although $p = 1$ is rarely of interest for DG methods, one possibility here is to use a $p > 1$ method as a preconditioner for the $p = 1$ system, or perhaps cluster elements into larger blocks for the Gauss-Seidel method. Meanwhile, the results of this paper could also be used to inform the design of algebraic multigrid methods \cite{Notay2017}.

Finally, we remark that the measured multigrid convergence rates for the time-dependent Stokes equations are particularly encouraging, indicating a strong potential for developing fast solvers for the general incompressible Navier-Stokes equations. In particular, preliminary work indicates that a fast Stokes solver could outperform the well-known and widely-applied projection method of Chorin \cite{Chorin,GuermondMinevShen}, and could be integrated into arbitrary-order time-stepping methods for Navier-Stokes \cite{sdpc}; this will be further investigated in future work.

\subsection*{Acknowledgements} This research was supported by the Applied Mathematics Program of the U.S.~DOE Office of Advanced Scientific Computing Research under contract number DE-AC02-05CH11231. Some computations used resources of the National Energy Research Scientific Computing Center (NERSC), a U.S.~DOE Office of Science User Facility operated under Contract No.~DE-AC02-05CH11231.

\clearpage

\appendix

\section{Local discontinuous Galerkin methods for the multi-phase Stokes equations}
\label{sec:ldg}

In this section, an LDG framework is derived for the variable-viscosity multi-phase Stokes equations. The construction partly follows the schemes set out by Cockburn \textit{et al}~\cite{CockburnKanschatSchotzauSchwab2002} but with some differences, including: (i) the formulation is derived in a way which makes the role of the discrete gradient operator and its adjoint more visible (relevant to the operator-coarsening multigrid schemes presented in \S\ref{sec:op}); (ii) the penalty stabilisation operators are separated out from the numerical fluxes; and (iii) the formulation is extended to treat variable viscosity and multi-phase Stokes problems. For reference, the governing equations are repeated here: we seek to compute a velocity field $\vu : \Omega \to {\mathbb R}^d$ and pressure field $p : \Omega \to {\mathbb R}$ such that
\begin{equation} \label{eq:appgovern1} \left. \begin{aligned} -\nabla \cdot \bigl(\mu_i (\nabla \vu + \gamma \nabla \vu^\trans) \bigr) + \nabla p &= \mathbf f \\ -\nabla \cdot \mathbf \vu &= f \end{aligned} \right\} \text{ in $\Omega_i$,} \end{equation}
subject to the interfacial jump conditions,
\begin{equation} \label{eq:appgovern2} \left. \begin{aligned} \jump{\vu} &= {\mathbf g}_{ij} \\ \jump{\mu (\nabla \vu + \gamma \nabla \vu^\trans) \vn - p \vn} &= {\mathbf h}_{ij} \end{aligned} \right\} \text{ on $\Gamma_{ij}$,} \end{equation}
and boundary conditions
\begin{equation} \label{eq:appgovern3} \begin{aligned} \vu &= {\mathbf g}_{\partial} && \text{on } \Gamma_D, \\ \mu (\nabla \vu + \gamma \nabla \vu^\trans) \vn - p \vn &= {\mathbf h}_{\partial} && \text{on } \Gamma_N. \end{aligned} \end{equation}

We begin with some preliminary set up and notation. As stated in the main article, we consider meshes arising from Cartesian grids or, for domains or interfaces with curved geometry, semi-unstructured quadtree/octree-based implicitly defined meshes. In this setting, it is natural to adopt a tensor-product piecewise polynomial space. Let $\mathcal E = \bigcup_i E_i$ denote the set of mesh elements, let $p \geq 1$ be an integer, and define ${\mathcal Q}_p(E)$ as the space of tensor-product polynomials of (one-dimensional) degree $p$ on element $E$. We assume in this work that the mesh is interface-conforming, i.e., if there is an interface separating the domain into two or more phases, then the interface does not cut through any element. Then, regarding the faces of the mesh, we denote \textit{intraphase faces} as those shared by two elements in the same phase; \textit{interphase faces} as those shared by two elements in differing phases (therefore situated on $\Gamma_{ij}$ for some $i,j$); and \textit{boundary faces} as those situated on $\partial \Omega$. Each face has an associated normal vector $\vn$; on intraphase faces, which are always flat and lie in a particular coordinate plan, $\vn$ is defined to point ``left-to-right'', e.g., for vertical faces in 2D, $\vn = \hat{\mathbf x}$. Interphase faces adopt the same normal vector as the interface $\Gamma_{ij}$ on which they coincide, defined to point from the phase with smallest phase index, $i$, into the phase with largest index, $j > i$. Boundary faces adopt the natural outwards-pointing normal to the domain boundary. The notation $\jump{\cdot}$ denotes the jump of a quantity across an interface or face and is defined consistent with its orientation; in particular, $\jump{u} := u^- - u^+$ where $u^\pm(x) = \lim_{\epsilon \to 0^+} u(x \pm \epsilon \vn)$ denotes the left and right trace values, $u^-$ and $u^+$, respectively. In addition, define $\Gamma_0$ as the set of all points belonging to intraphase faces, and for an element $E \in \mathcal E$, define $\chi(E)$ to be the phase of that element, such that $E \subseteq \Omega_{\chi(E)}$.

In the first step of the LDG formulation, a discrete approximation to $\nabla \vu$ is defined through a ``strong-weak'' form.\footnote{The \textit{strong-weak form} states that $\bmg$ must satisfy \[ \int_E \bmg : \vomega = \int_E \nabla \vu : \vomega + \int_{\partial E} ({\hat \vu} - \vu) \cdot \vomega \cdot \vn \] for all test functions $\vomega$, whereas the \textit{weak form} states that $\bmg$ must satisfy \[ \int_E \bmg : \vomega = -\int_E \vu \cdot (\nabla \cdot \vomega) + \int_{\partial E} {\hat \vu} \cdot \vomega \cdot \vn.\] The two forms are equivalent whenever the associated quadrature scheme exactly preserves the identity of integration by parts. This is often the case for many implementations, including on quadrilateral, prismatic or simplicial elements. However, this property may not hold when approximate quadrature schemes are used, e.g., for implicitly defined meshes which use high-order accurate quadrature schemes (wherein integration by parts only holds up to a high-order truncation error). In the latter situation, to ensure symmetry of the final discretised Stokes operator, it is necessary to use the strong-weak form to define $\nabla \vu$ and the weak form to defined $\nabla \cdot \vsigma$ (see \eqref{eq:wdef}), or vice versa. For further discussion as it relates to the analogous case of elliptic interface problems, see ref.~\cite{ImplicitMeshPartOne}.} Given $\vu \in V_h$, $\bmg \in V_h^{d \times d}$ is defined such that
\begin{equation} \label{eq:geqn} \int_E \bmg : \bm \vomega = \int_E \nabla \vu : \bm \vomega + \int_{\partial E} ({\hat \vu}_{\chi(E)} - \vu) \cdot \bm \omega \cdot \bm n \end{equation}
holds for every element $E \in \mathcal E$ and every test function $\vomega \in V_h^{d \times d}$, where $\hat \vu_\chi$ is a numerical flux function defined as follows:
\begin{equation} \hat \vu_\chi := \begin{cases} \vu^- & \text{on any intraphase face}, \\ \lambda \vu^-  + (1-\lambda) (\vu^+ + {\mathbf g}_{\chi i}) & \text{on $\Gamma_{\chi i}$ if $\chi < i$}, \\ \lambda (\vu^- - {\mathbf g}_{i \chi}) + (1-\lambda) \vu^+ & \text{on $\Gamma_{i \chi}$ if $\chi > i$}, \\ \vu^- & \text{on } \Gamma_N, \\ {\mathbf g}_\partial & \text{on } \Gamma_D. \end{cases} \label{eq:ustarflux} \end{equation}

\begin{figure}[t]
\centering
\footnotesize
\sffamily
\scalebox{0.855}{\def\svgwidth{6in}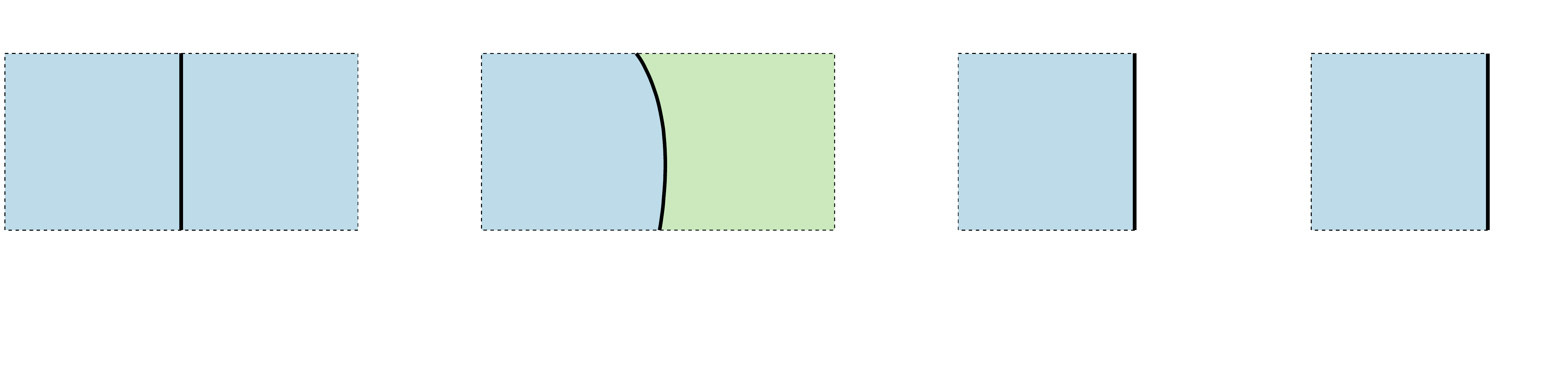}
\caption{Schematic of the numerical flux functions $\hat \vu$ and $\hat \vsigma$ defined by \eqref{eq:ustarflux} and \eqref{eq:qstarflux}. Except for interphase faces, the flux is single-valued; on interphase faces, the flux is multi-valued so as to incorporate the interfacial jump conditions $\jump{\vu} = {\mathbf g}_{ij}$ and $\jump{\vsigma \cdot \vn} = {\mathbf h}_{ij}$ on $\Gamma_{ij}$, $i < j$. A plus and minus sign denote the elemental values on the right and left of the face, respectively; e.g., for a point $x$ on the face, $\vu^\pm(x) = \lim_{\epsilon \to 0^+} \vu(x \pm \epsilon \vn)$.}
\label{fig:fluxdiagram}
\end{figure}

\noindent (See \cref{fig:fluxdiagram} for a schematic illustration.) Note that the flux is multivalued on interphase faces---on these faces, the interfacial jump condition $\jump{\vu} = {\mathbf g}_{ij}$ on $\Gamma_{ij}$ in \eqref{eq:appgovern2} is taken into account as follows: when an element ``reaches across'' the interface to evaluate the trace of $\vu$ on the other side, the trace value is compensated by the jump data to correctly account for the intended discontinuity in the solution. Note also that interfacial fluxes are weighted through a convex combination parameter $\lambda \in [0,1]$, which can vary from face to face. (For example, if $\lambda = 0$, then the numerical flux is sourced solely from the right element's trace $\vu^+$; if $\lambda = 1$, it is sourced solely from the left element's trace $\vu^-$.) The purpose of $\lambda$ is to implement the strategy of viscosity-upwinded numerical fluxes, as detailed in \S\ref{sec:viscosityupwind}. Upon summing \eqref{eq:geqn} over every element of the mesh, one has
\begin{equation} \label{eq:gsum} \begin{aligned} \int_{\Omega} \bmg : \vomega &= \sum_{E \in {\mathcal E}} \nabla \vu : \vomega + \int_{\Gamma_0} \bigl[ (\hat \vu - \vu^-) \cdot \vomega^- \cdot \vn - (\hat \vu - \vu^+) \cdot \vomega^+ \cdot \vn \bigr] \\ & \qquad + \sum_{j > i} \int_{\Gamma_{ij}} \bigl[ (\hat \vu_i - \vu^-) \cdot \vomega^- \cdot \vn - (\hat \vu_j - \vu^+) \cdot \vomega^+ \cdot \vn \bigr] \\
& \qquad + \int_{\Gamma_D} (\hat \vu - \vu^-) \cdot \vomega^- \cdot \vn + \int_{\Gamma_N} (\hat \vu - \vu^-) \cdot \vomega^- \cdot \vn \\
&= \sum_{E \in {\mathcal E}} \nabla \vu : \vomega - \int_{\Gamma_0} (\vu^- - \vu^+) \cdot \vomega^+ \cdot \vn \\ & \qquad - \sum_{j > i} \int_{\Gamma_{ij}} \bigl[ (1- \lambda) (\vu^-  - \vu^+) \cdot \vomega^- \cdot \vn + \lambda (\vu^-  - \vu^+ ) \cdot \vomega^+ \cdot \vn \bigr]  \\ & \qquad + \sum_{j > i} \int_{\Gamma_{ij}} (1 - \lambda) {\mathbf g}_{ij} \cdot \vomega^- \cdot \vn + \lambda\, {\mathbf g}_{ij} \cdot \vomega^+ \cdot \vn \\
& \qquad + \int_{\Gamma_D} ({\mathbf g}_\partial - \vu^-) \cdot \vomega^- \cdot \vn, 
\end{aligned} \end{equation}
which motivates the definition of the following operators:
\begin{itemize}
\item Let $\nabla_h : V_h \to V_h^d$ be the broken gradient operator and let $L: V_h \to V_h^d$ be the lifting operator, such that
\[ \int_\Omega (\nabla_h u) \cdot {\bm v} = \sum_{E \in \mathcal E} \int_E \nabla u \cdot {\bm v},\]
and
\begin{align*} \int_\Omega (Lu) \cdot {\bm v} &= \sum_{j > i} \int_{\Gamma_{ij}} (1 - \lambda) (u^+ - u^-) {\bm v}^- \cdot \vn + \lambda (u^+ - u^-) {\bm v}^+ \cdot \vn \\ & \qquad + \int_{\Gamma_0} (u^+ - u^-) {\bm v}^+ \cdot \vn - \int_{\Gamma_D} u^- {\bm v} ^- \cdot \vn,\end{align*}
holds for every ${\bm v} \in V_h^d$.
\item Define $J_g \in V_h^{d \times d}$ such that 
\[ \int_\Omega J_g : \vomega = \int_{\Gamma_D} {\mathbf g}_\partial \cdot \vomega^- \cdot \vn + \sum_{j > i} \int_{\Gamma_{ij}} (1 - \lambda) {\mathbf g}_{ij} \cdot \vomega^- \cdot \vn + \lambda\, {\mathbf g}_{ij} \cdot \vomega^+ \cdot \vn, \]
holds for every $\vomega \in V_h^{d \times d}$.
\end{itemize}
With these definitions, \eqref{eq:gsum} is equivalent to the statement that
\[ \bmgg_{ij} = G_j u_i + J_{g,ij}, \]
where $G : V_h \to V_h^d$ is the \textit{discrete gradient operator}, $G := \nabla_h + L$, having components $G = (G_1,\ldots,G_d)$, $\bmgg_{ij}$ denotes the $(i,j)$th component of $\bmg$, and $u_i$ denotes the $i$th component of $\vu$. To complete this step of the LDG construction, we define $\vtau$ as the natural discretisation of $\nabla \vu + \gamma\nabla \vu^\trans$,
\[ \vtau := {\bmg} + \gamma {\bmg}^\trans, \quad \text{i.e.,} \quad \tau_{ij} := G_j u_i + \gamma G_i u_j + J_{g,ij} + \gamma J_{g,ji}. \]

In the second step of the LDG formulation, a discrete approximation of $\vsigma = \mu \vtau - p {\mathbb I}$ is defined. In essence, this is implemented via an $L^2$ projection of $\mu \vtau$ onto $\smash{V_h^{d \times d}}$. We define $\smash{\vsigma \in V_h^{d \times d}}$ as the unique piecewise polynomial function such that
\begin{equation} \label{eq:sigmadef} \int_E \vsigma : \vomega = \int_E (\mu \vtau - p {\mathbb I}): \vomega \end{equation}
holds for every element $E \in \mathcal E$ and every test function $\vomega \in V_h^{d \times d}$. In the case that $\mu$ is piecewise constant, computing this $L^2$ projection is a simple matter of multiplying $\vtau$ by a scalar and subtracting the discrete pressure $p \in V_h$. In the general case, it is straightforward to show that a $\mu$-weighted $L^2$ projection is equivalent to multiplication by the block diagonal matrix $M^{-1} M_\mu$, where $M_\mu$ is the $\mu$-weighted mass matrix such that $u^\trans M_\mu v = \int_\Omega u\,\mu\,v$ holds for all $u, v \in V_h$. In this work, $M_\mu$ is computed with sufficiently high-order accurate quadrature schemes, typically Gaussian quadrature schemes. Using this relation, \eqref{eq:sigmadef} is equivalent to
\begin{equation} \label{eq:sigmaresult} \sigma_{ij} := M^{-1} M_\mu \tau_{ij} - p\,\delta_{ij} = M^{-1} M_\mu \bigl( G_j u_i + \gamma G_i u_j + J_{g,ij} + \gamma J_{g,ji} \bigr) - p\, \delta_{ij}. \end{equation}

In the third step, we consider a weak formulation for computing the divergence of $\vsigma$. This proceeds similarly to defining the discrete gradient of $\vu$, except numerical fluxes act in the opposite direction. For simplicity of presentation, the following numerical flux for $\vsigma$ is matrix-valued; however, only the normal component of the flux is used. Given $\vsigma \in V_h^{d \times d}$, define $\bm w \in V_h^d$ as the discrete divergence of $\vsigma$ such that
\begin{equation} \label{eq:wdef} \int_E \bm w \cdot \bm v = -\int_E \vsigma : \nabla \bm v + \int_{\partial E} \bm v \cdot {\hat \vsigma}_{\chi(E)} \cdot \vn \end{equation}
holds for every test function $\bm v \in V_h^d$ and every element $E \in \mathcal E$. Here, the numerical flux is defined by (see also \cref{fig:fluxdiagram})
\begin{equation} {\hat \vsigma}_\chi := \begin{cases} \vsigma^+ & \text{on any intraphase face}, \\ (1-\lambda) \vsigma^-  + \lambda (\vsigma^+ + {\mathbf h}_{\chi i} \otimes \vn) & \text{on $\Gamma_{\chi i}$ if $\chi < i$}, \\ (1 - \lambda) (\vsigma^- - {\mathbf h}_{i \chi} \otimes \vn) + \lambda \vsigma^+ & \text{on $\Gamma_{i \chi}$ if $\chi > i$}, \\ {\mathbf h}_\partial \otimes \vn & \text{on } \Gamma_N, \\ \vsigma^- & \text{on } \Gamma_D. \end{cases} \label{eq:qstarflux} \end{equation}
As in the numerical flux for $\hat \vu$, the interfacial jump condition $\jump{\vsigma \cdot \vn} = {\mathbf h}_{ij}$ on $\Gamma_{ij}$ in \eqref{eq:appgovern2} is taken into account via a multivalued interfacial flux, such that whenever an element reaches across the interface, the neighbouring element's trace is compensated by ${\mathbf h}_{ij}$ to correctly put it in the context of the source element. Summing \eqref{eq:wdef} over every element, one has, for every $\bm v \in V_h^d$,
\begin{equation} \label{eq:wworkings} \begin{aligned}
  \int_\Omega \bm w \cdot \bm v &= -\sum_{E \in \mathcal E} \int_E \vsigma : \nabla \bm v + \int_{\Gamma_0} (\bm v^- - \bm v^+) \cdot {\hat \vsigma} \cdot \vn \\ & \qquad + \sum_{j > i} \int_{\Gamma_{ij}} (\bm v^- \cdot {\hat \vsigma}_i - \bm v^+ \cdot {\hat \vsigma}_j) \cdot \vn \\ & \qquad + \int_{\Gamma_D} \bm v^- \cdot {\hat \vsigma} \cdot \vn + \int_{\Gamma_N} \bm v^- \cdot {\hat \vsigma} \cdot \vn \\
	&= -\sum_{E \in \mathcal E} \int_E \vsigma : \nabla \bm v + \int_{\Gamma_0} (\bm v^- - \bm v^+) \cdot \vsigma^+ \cdot \vn \\ & \qquad  + \sum_{j > i} \int_{\Gamma_{ij}} (1 - \lambda) (\bm v^- - \bm v^+) \cdot \vsigma^- \cdot \vn + \lambda (\bm v^- - \bm v^+) \cdot \vsigma^+ \cdot \vn  \\ & \qquad + \sum_{j > i} \int_{\Gamma_{ij}} \lambda \bm v^- \cdot {\mathbf h}_{ij} + (1-\lambda) \bm v^+ \cdot {\mathbf h}_{ij} \\
	& \qquad + \int_{\Gamma_D} \bm v^- \cdot \vsigma^- \cdot \vn + \int_{\Gamma_N} \bm v^- \cdot {\mathbf h}_\partial. 
\end{aligned} \end{equation}
Similar to the operator $J_g$ defined above, let $J_h \in V_h^d$ be such that
\[ \int_\Omega J_h \cdot \bm v = \int_{\Gamma_N} \bm v^- \cdot {\mathbf h}_\partial + \sum_{j > i} \int_{\Gamma_{ij}} \lambda \bm v^- \cdot {\mathbf h}_{ij} + (1-\lambda) \bm v^+ \cdot {\mathbf h}_{ij}\]
holds for every $\bm v \in V_h^d$. Then, upon using the lifting operator defined earlier, \eqref{eq:wworkings} is equivalent to the statement that, for every $\bm v \in V_h^d$,
\[ (\bm w, \bm v) = -\sum_{i=1}^d (\vsigma_i, \nabla_h v_i) - (\vsigma_i, L v_i) + (J_h, \bm v) = -\sum_{i=1}^d (\vsigma_i, G v_i) + (J_h, \bm v), \]
where $v_i$ denotes the $i$th component of $\bm v$ and $\vsigma_i$ denotes the $i$th row of the matrix $\vsigma$. Transferring the $G$ operator onto $\vsigma_i$ via the adjoint, it follows that $w_i = -\sum_{j=1}^d M^{-1} G_j^\trans M \sigma_{ij} + J_{h,i}$, where $w_i$ is the $i$th component of $\bm w$. Combining with \eqref{eq:sigmaresult}, we have
\begin{equation} \label{eq:wform} \begin{aligned} w_i &= - \sum_{j=1}^d  M^{-1} G_j^\trans M_\mu \bigl( G_j u_i + \gamma\,G_i u_j \bigr) + M^{-1} G_i^\trans M p \\ & \qquad -  \sum_{j=1}^d M^{-1} G_j^\trans (J_{g,ij} + \gamma J_{g,ji}) + J_{h,i}.  \end{aligned} \end{equation}
This is the weak statement that $\bm w$ is equal to the discrete divergence of $\vsigma$, itself a discrete approximation to $\mu(\nabla \vu + \gamma \nabla \vu^\trans) - p {\mathbb I}$, taking into account velocity Dirichlet boundary data, stress boundary data, and interfacial jump condition data (if any). One may recognise the first term of \eqref{eq:wform} as implementing $\nabla \cdot (\mu (\nabla \vu + \gamma \nabla \vu^\trans))$, the second term as implementing $-\nabla p$, while the remaining terms represent the contribution of any boundary or interfacial jump data. 

We now turn to the LDG discretisation of the divergence constraint of the Stokes equations. Given $\vu \in V_h^d$, define $w \in V_h$ as the discrete divergence of $\vu$ via the strong-weak form, such that
\begin{equation} \label{eq:ewv} \int_E w\, v = \int_E v\,\nabla \cdot \vu + \int_{\partial E} v\, ({\hat \vu}_{\chi(E)} - \vu) \cdot \vn \end{equation}
holds for every element $E \in \mathcal E$ and every $v \in V_h$. Here, the same numerical flux for $\vu$ as was used to define its gradient (see \eqref{eq:ustarflux}) is employed, but only the normal component is seen by \eqref{eq:ewv}. A similar derivation as before reveals that $w$ is essentially equal to the trace of $\bmg$, i.e.,
\[ w = \sum_{i=1}^d G_i u_i + J_{g\cdot n}, \]
where $J_{g\cdot n} \in V_h$ is such that $\int_\Omega J_{g\cdot n} v = \int_{\Gamma_D} v^- {\mathbf g}_\partial \cdot \vn + \sum_{j > i} \int_{\Gamma_{ij}} ((1 - \lambda) v^- + \lambda v^+) {\mathbf g}_{ij} \cdot \vn$ holds for every $v \in V_h$.

Last, in what is essentially the final step of the LDG formulation for the Stokes equations \eqref{eq:appgovern1}--\eqref{eq:appgovern3}, penalty stabilisation terms are added to ensure the well-posedness of the discrete problem \cite{CockburnKanschatSchotzauSchwab2002,ArnoldBrezziCockburnMarini,HesthavenWarburton}. These terms weakly impose solution continuity between neighbouring element polynomials in the same phase, weakly impose Dirichlet boundary conditions, and weakly enforce interfacial jump conditions (if any). Regarding the stabilisation parameters for velocity, we classify them according to three types: boundary ($\tau_{\vu,\partial}$), intraphase ($\tau_{\vu,0}$), and interphase ($\tau_{\vu,ij}$). Let ${\mathcal E}_{\vu,g} : V_h^d \to V_h^d$ be the affine operator such that, for each $\vu \in V_h$,
\begin{equation} \label{eq:Epenalty} \begin{aligned} \int_\Omega {\mathcal E}_{\vu,g}(\vu) \cdot \bm v &= \int_{\Gamma_0} \tau_{\vu,0} \jump{\vu} \cdot \jump{\bm v} + \sum_{j > i} \int_{\Gamma_{ij}} \tau_{\vu,ij} (\jump{\vu} - {\mathbf g}_{ij}) \cdot \jump{\bm v} \\ & \qquad + \int_{\Gamma_D} \tau_{\vu,\partial} (\vu^- - {\mathbf g}_\partial) \cdot \bm v^- \end{aligned} \end{equation}
holds for every test function $\bm v \in V_h^d$. Note that ${\mathcal E}_{\vu,g}(\vu) = ({\mathcal E}_{\vu,0} u_1,\ldots,{\mathcal E}_{\vu,0} u_d) + {\mathcal E}_{\vu,g}(0)$, where ${\mathcal E}_{\vu,0}$ represents the linear part of the operator acting on the vector field components, defined analogous to \eqref{eq:Epenalty} with homogeneous jump and boundary data. Concerning the stabilisation parameter $\tau_p$ for pressure, define the linear operator ${\mathcal E}_p : V_h \to V_h$ such that, for each $p \in V_h$,
\[ \int_\Omega {\mathcal E}_p(p) v = \int_{\Gamma_0} \tau_p \jump{p} \jump{v} \]
holds for every test function $v \in V_h$. Following the formulation presented in ref.~\cite{CockburnKanschatSchotzauSchwab2002}, subject to a suitable specification of the parameter values to be discussed shortly, these operators are added and subtracted to the discretisation of the Stokes momentum equations and divergence constraint, respectively. Specifically, the discretised multi-phase Stokes problem \eqref{eq:appgovern1}--\eqref{eq:appgovern3} consists of finding a velocity field $\vu \in V_h^d$ and pressure field $p \in V_h$ such that
\begin{equation} \label{eq:appstokes1} \begin{aligned} & \sum_{j=1}^d M^{-1} G_j^\trans M_\mu (G_j u_i + \gamma\, G_i u_j) + {\mathcal E}_{\vu,0} u_i - M^{-1} G_i^\trans M p \\ & \qquad = {\mathbb P}_{V_h^d}(f_i) - \sum_{j=1}^d M^{-1} G_j^\trans (J_{g,ij} + \gamma J_{g,ji}) + J_{h,i} - {\mathcal E}_{\vu,g,i}(0) \end{aligned} \end{equation}
holds for each $i = 1,\ldots,d$, where $f_i$ denotes the $i$th component of $\mathbf f$, subject to the divergence constraint
\begin{equation} \label{eq:appstokes2} -\sum_{j=1}^d G_i u_i - {\mathcal E}_p p = {\mathbb P}_{V_h} (f) + J_{g \cdot n}. \end{equation}
Multiplying both \eqref{eq:appstokes1} and \eqref{eq:appstokes2} by $M$, taken together these equations may be succinctly written in block form as
\begin{equation} \label{eq:stokessystem} \begin{pmatrix} A & M {\mathcal G} \\ {\mathcal G}^\trans M & -E \end{pmatrix} \begin{pmatrix} \mathbf u \\ p \end{pmatrix} = \begin{pmatrix} {\mathbf b}_{\vu} \\ {b}_p \end{pmatrix} \end{equation}
where the block operators are as follows:
\begin{itemize}
\item $A$ implements the viscous part of the Stokes momentum equations, and can be written in $d \times d$ block form corresponding to its action on the $d$ components of $\vu$, where the $(i,j)$th block is given by
\begin{equation} \label{eq:viscousop} A_{ij} = \delta_{ij} \bigl( \textstyle{\sum_{k=1}^d} G_k^\trans M_\mu G_k \bigr) + \gamma\, G_j^\trans M_\mu G_i + \delta_{ij} \tilde{E}, \end{equation}
where $\tilde E := M {\mathcal E}_{\vu,0}$ is the penalty stabilisation matrix associated with velocity stabilisation.
\item $\mathcal G = ({\mathcal G}_1,\ldots,{\mathcal G}_d)$ is a discrete gradient operator, closely related to the adjoint of $G$, whose components are given by
\[ {\mathcal G}_i = -M^{-1} G_i^\trans M. \]
\item $E$ is the penalty stabilisation matrix associated with pressure, defined by $E := M {\mathcal E}_p$.
\item Last, $({\mathbf b}_{\vu}, b_p)$ collects the entire influence of the right hand side data, $({\mathbf f}, f)$, together with Dirichlet, stress, and interfacial jump source data, and corresponds to the multiplication by $M$ of the right hand sides of \eqref{eq:appstokes1} and \eqref{eq:appstokes2}.
\end{itemize}

\subsection{Specification of penalty parameters}

As described, four different kinds of penalty parameters need specification for the LDG formulation of the Stokes equations---three for velocity on the boundary, intraphase, and interphase faces of the mesh, and one for pressure stabilisation. Remarks concerning their general specification and particular choices made in this work are provided here.
\begin{itemize}
\item In general, strictly positive parameters are sufficient to ensure well-posedness of the final linear system, i.e., to ensure it has the expected trivial kernel of the continuum Stokes operator, and to ensure the inf-sup conditions hold \cite{CockburnKanschatSchotzauSchwab2002,CockburnKanschatSchotzau2003,CockburnKanschatSchotzau2004}.\footnote{Subtleties may arise in the multi-phase case owing to the non-penalisation of pressure jumps across the interface, depending on mesh topology and the choice of numerical fluxes $\hat\vu$ and $\hat\vsigma$. Discussion is deferred to future work.} However, this is not a necessary condition. For example, on a regular Cartesian mesh, with purely one-sided intraphase numerical fluxes for $\hat \vu$ and $\hat \vsigma$ that ``upwind'' in the opposite direction, as used in this work, one can set the intraphase penalty parameter for velocity to zero, $\tau_{\vu,i} = 0$, see, e.g., ref.~\cite{Cockburn0701} in the case of LDG for scalar Poisson problems. On the other hand, a choice of penalty parameter which is too large can impact discretisation accuracy and overall conditioning of the final linear operator as well as multigrid solver efficiency.
\item If $\Gamma_D$ is nonempty, then $\tau_{\vu,\partial}$ must be positive to ensure well-posedness.
\item Although no specific lower bounds on parameter values are required for the LDG system to be well-posed, for consistent discretisation behaviour (including invariance with respect to stretching the mesh as well as preserving spectral characteristics as the mesh is refined), velocity penalty parameters should scale proportional to $h^{-1}$ and pressure penalty parameters should scale proportional to $h$. 
For anisotropic meshes, one can be more precise and require that penalty parameters on a particular mesh face scale appropriately with the measure of the face divided by the measure of the elements on either side \cite{CockburnKanschatSchotzauSchwab2002}, but this extra kind of precision is not pursued here.
\item To ensure correct scaling with respect to ellipticity coefficient, penalty parameters should also scale with viscosity. For velocity penalty parameters, this implies $\tau_{\vu,\partial} \sim \mu^-$ and $\tau_{\vu,0} \sim \mu$, where $\mu$ is the local value of the viscosity on the mesh face in question, and for interphase velocity penalty parameters, $\tau_{ij}$ should scale with the smaller of the two viscosities on either side of the interface \cite{fluxx}. Meanwhile, the pressure penalty parameter should scale inversely-proportional to the local viscosity to preserve the spectral characteristics of the Stokes operator; see, e.g., \cite{CockburnKanschatSchotzau2003}.
\item One can also choose to scale $\tau$ with the polynomial degree, and this can be important for the study of DG methods utilising very high degree polynomials; e.g., one could scale according to $\tau_{\vu} \sim p^2$. However, in this work, only moderate-order polynomials are used and a linear scaling in $p$ is applied, as defined next.
\end{itemize}
The specification of the pressure penalty stabilisation parameter is one of the main subjects of this work and is discussed in the main article. Regarding the remaining penalty parameters, unless otherwise specified, the following values have been used throughout this work:
\begin{itemize}
\item On faces associated with the imposition of velocity Dirichlet boundary conditions, $\tau_{\vu,\partial} = 10 p \mu^- / h$, where $p$ is the (one-dimensional) polynomial degree, $\mu^-$ is the local viscosity of the face, and $h$ is the typical element size.
\item For intraphase penalty parameters, $\tau_{\vu,0} = 0$ on all Cartesian mesh examples, and $\tau_{\vu,0} = 0.5 p \mu / h$ on all examples using semi-unstructured meshes, where $\mu$ is the local viscosity of the face.
\item For interphase penalty parameters, $\tau_{\vu,ij} = C p \min(\mu^-,\mu^+) / h$, where $\mu^\pm$ is the local viscosity of the two phases on either side of the interfacial mesh face. Here, $C = 3$ in the case of Cartesian grid meshes, whereas $C = 8$ in the case of implicitly defined meshes, which benefit from slightly increased penalty stabilisation owing to their marginally less uniform variety of element shapes next to the interface.
\item Last, regarding the time-dependent Stokes equations and its generalised form of the pressure penalty parameter $\tau_p$ (see \eqref{eq:taupgeneral}), $\tau_0 = 0.5 p$ is the prefactor used for the limiting case of vanishing viscosity.
\end{itemize}
These values have been chosen based on the the dual goals of obtaining high-order discretisation accuracy as well as good multigrid performance, and follow typical values employed in prior work \cite{fluxx}.

\subsection{Viscosity-upwinded numerical fluxes}
\label{sec:viscosityupwind}

As discussed in the main article, viscosity-upwinded numerical fluxes \cite{fluxx} are utilised for the multi-phase Stokes problems. This corresponds to defining $\lambda$ in the numerical flux functions $\hat\vu$, \eqref{eq:ustarflux}, and $\hat\vsigma$, \eqref{eq:qstarflux}, for interfacial mesh faces as follows (see also \cref{fig:fluxdiagram}):
\begin{equation} \label{eq:viscosityupwind} \lambda = \begin{cases} 0 & \text{if } \mu^- < \mu^+, \\ 0.5 & \text{if } \mu^- = \mu^+, \\ 1 & \text{if } \mu^- > \mu^+. \end{cases} \end{equation}
Accordingly, $\hat\vu$ is biased to the more viscous phase, and $\hat\vsigma$ is biased to the less viscous phase, with interfacial jump data ${\mathbf g}_{ij}$ and ${\mathbf h}_{ij}$ incorporated appropriately.

\section{Operator coarsening geometric multigrid methods}
\label{sec:op}

In this section, we describe an operator coarsening strategy which allows one to compute the discrete Stokes operator ${\mathcal A}_h$ on every level of the mesh hierarchy, without having to explicitly build the mesh. This approach is equivalent in function to a purely-geometric multigrid method, but provides a variety of convenient benefits, including: (i) elements and faces do not need to be enumerated on coarse meshes; (ii) construction of quadrature schemes for coarse mesh elements or faces can be avoided; (iii) the viscosity coefficient $\mu$ is automatically coarsened down the mesh hierarchy in a manner consistent with performing repeated $L^2$ projections; and (iv) LDG operators, such as the discrete gradient and penalty stabilisation operators, are built automatically such that the chosen numerical flux functions of the finest mesh are inherited consistently by the coarse meshes. The technique is also amenable to simple block-sparse linear algebra routines, providing an opportunity to optimise the implementation. These schemes were originally derived for LDG methods applied to elliptic problems in ref.~\cite{dgmg} and then extended to variable-coefficient elliptic interface problems with large viscosity jumps in ref.~\cite{fluxx}; here, they are extended to variable-viscosity multi-phase Stokes problems. 

The operator coarsening approach is similar to the ``RAT'' paradigm often seen in multigrid methods, where the coarse mesh operator is defined by the fine mesh operator, pre- and post-multiplied by the restriction (``R'') and interpolation (``T'') operators, respectively. Such an approach is often applied to the primary elliptic operator itself, e.g., the discrete Laplacian operator. However, as shown in ref.~\cite{dgmg}, this approach leads to breakdown of multigrid performance for LDG methods. Instead, a more appropriate strategy is to apply ``RAT'' to the individual discrete gradient and divergence operators underlying the LDG method. In ref.~\cite{dgmg}, it is proven this is equivalent to a purely-geometric method, and it is straightforward to extend the methods therein to prove the same is true for the LDG formulation of the Stokes equations described above in \S\ref{sec:ldg}. We summarise this construction here and omit the proof.

Given a general operator $A: V_h \to V_h$ defined on a fine mesh, its coarsened counterpart ${\mathcal C}(A) : V_{2h} \to V_{2h}$ on a coarse mesh is defined variationally. Specifically, ${\mathcal C}(A)$ is defined such that
\[ \bigl( {\mathcal C}(A) u, v \bigr)_{V_{2h}} = (A I_{2h}^h u, I_{2h}^h v)_{V_h} \]
for all $u, v \in V_{2h}$; here $(\cdot,\cdot)$ denotes the standard inner product, and $V_{2h}$ denotes the piecewise polynomial space associated with the coarse mesh. Equivalently, as a matrix acting on coefficient vectors in the chosen basis, 
\[ {\mathcal C}(A) = R_h^{2h} A I_{2h}^h = M_{2h}^{-1} (I_{2h}^h)^\trans M_h A I_{2h}^h \]
where $M_h$ and $M_{2h}$ are the mass matrices of the two meshes. 

\subsection{Time-independent Stokes equations}

We derive the operator coarsening strategy for the Stokes system \eqref{eq:stokessystem}, as follows:
\begin{enumerate}
\item The mass matrix of the coarse mesh is given by $M_{2h} = (I_{2h}^h)^\trans M_h I_{2h}^h$ \cite{dgmg}. If $M_{\mu,h}$ is a $\mu$-weighted mass matrix on the fine mesh, its coarsened counterpart is given by $M_{\mu,2h} = (I_{2h}^h)^\trans M_{\mu,h} I_{2h}^h$ \cite{fluxx}.
\item The coarsened operators making up the Stokes operator are given by ${\mathcal C}(G_h)$, ${\mathcal C}({\mathcal G}_h)$, $\tfrac12 M_{2h} {\mathcal C}(M_h^{-1} {\tilde E}_h)$, and $2\,M_{2h} {\mathcal C}(M_h^{-1} E_h)$ for the discrete gradient operator, the adjoint form of the discrete gradient operator, the velocity penalty stabilisation operator, and the pressure penalty stabilisation operator, respectively. In particular, note the $\tfrac12$ and $2$ prefactors in the coarsened penalty operators---these correspond to the observation that velocity penalty parameters scale proportional to $h^{-1}$ and pressure penalty parameters scale proportional to $h$, and that one should preserve this scaling across the full multigrid hierarchy to attain ideal multigrid performance \cite{fluxx}.
\item Last, the discrete Stokes operator on the coarse mesh is formed by computing the viscous operator in \eqref{eq:viscousop} using the coarsened discrete gradient operators, and then building the overall operator using the form given in \eqref{eq:stokessystem}.
\end{enumerate}

\begin{algorithm}[t]
	\caption{\sffamily\small Construction of coarse mesh operators for the time-independent multi-phase Stokes equations, given fine mesh operators $M_h$, $M_{\mu,h}$, $G_h$, ${\mathcal G}_h$, $\tilde{E}_h$, and $E_h$.}
	\begin{algorithmic}[1]
		\State $M_{2h} := (I_{2h}^h)^\trans M_h I_{2h}^h$
		\State $M_{\mu,2h} := (I_{2h}^h)^\trans M_{\mu,h} I_{2h}^h$
		\State $G_{2h} := M_{2h}^{-1} (I_{2h}^h)^\trans M_h G_h I_{2h}^h$
		\State ${\mathcal G}_{2h} := M_{2h}^{-1} (I_{2h}^h)^\trans M_h {\mathcal G}_h I_{2h}^h$; (equivalently, ${\mathcal G}_{2h,i} := -M_{2h}^{-1} G_{2h,i}^\trans M_{2h}$)
		\State $\tilde{E}_{2h} := \tfrac12 (I_{2h}^h)^\trans \tilde{E}_h I_{2h}^h$
		\State $E_{2h} := 2 (I_{2h}^h)^\trans E_h I_{2h}^h$
		\State Build the coarsened viscous operator in $d \times d$ block form, with $(i,j)$th block given by
		\[ A_{2h,ij} := \delta_{ij} \bigl( \textstyle{\sum_{k=1}^d} G_{2h,k}^\trans M_{\mu,2h} G_{2h,k} \bigr) + \gamma\, G_{2h,j}^\trans M_{\mu,2h} G_{2h,i} + \delta_{ij} \tilde{E}_{2h} \]
		\State Form the coarsened Stokes operator,
		\[ {\mathcal A}_{2h} := \begin{pmatrix} A_{2h} & M_{2h} {\mathcal G}_{2h} \\ {\mathcal G}_{2h}^\trans M_{2h} & -E_{2h} \end{pmatrix} \]
	\end{algorithmic}
	\label{algo:build}
\end{algorithm}

\noindent \cref{algo:build} defines the overall operator coarsening scheme, to be applied recursively down the mesh hierarchy. We note that in this algorithm, both of the discrete gradient operators, $G$ and ${\mathcal G}$, are coarsened; this is to assist in overall clarity, however it also suffices to coarsen just one.

\subsection{Time-dependent Stokes equations}

As discussed in the main article, a modification to the pressure penalty parameter scaling is appropriate for the time-dependent Stokes equations. In this setting, the relation $\tau_p = \tau\, h/\mu$ is replaced with
\begin{equation} \label{eq:taupgeneral} \tau_p = \Bigl( \frac{h \rho}{\tau_0 \delta} + \frac{\mu}{\tau h} \Bigr)^{-1}. \end{equation}
To incorporate this scaling into the operator coarsening strategy described above, the prefactor of $2$ in the coarsening of the pressure stabilisation operator, $E_{2h} = 2\,M_{2h} {\mathcal C}(M_h^{-1} E_h)$, should be appropriately modified to reflect the property that the scaling in \eqref{eq:taupgeneral} may change from $\tau_p \sim h$ to $\tau_p \sim 1/h$ down the multigrid hierarchy. In this work, this is implemented through the coarsening of \textit{two} stabilisation operators for pressure, which are then combined into a net result via a harmonic weighting. Specifically, on the finest mesh we define $E_{\mu}$ and $E_{\rho}$ such that
\[ u^T E_\mu v = \int_{\Gamma_0} \tfrac{\tau h}{\mu} \jump{u} \jump{v}, \quad \text{and} \quad u^T E_\rho v = \int_{\Gamma_0} \tfrac{\tau_0 \delta}{h \rho} \jump{u} \jump{v} \]
holds for every $u, v \in V_h$. Then, $E_{\mu}$ is coarsened using a prefactor of $2$, and $E_\rho$ is coarsened using a prefactor of $\frac12$. The two block-sparse operators are then combined using a heuristic based on the Frobenius norm of each block. Specifically, if $E_{\mu,ij}$ denotes the $(i,j)$th block of $E_{\mu}$ (similarly for $E_\rho$), where individual blocks correspond to individual elements of the mesh, then we define $E$ such that
\begin{equation} \label{eq:Etimedep} E_{ij} = \Bigl(\frac{\|E_{\rho,ij}\|_F}{\|E_{\mu,ij}\|_F + \|E_{\rho,ij}\|_F} \Bigr)^2 E_{\mu,ij} + \Bigl(\frac{\|E_{\mu,ij}\|_F}{\|E_{\mu,ij}\|_F + \|E_{\rho,ij}\|_F} \Bigr)^2 E_{\rho,ij}. \end{equation}
This formula is based on the identity $(1/x + 1/y)^{-1} = (x^2 y + y^2 x)/(x+y)^2$ and is simply a heuristic means of recombining $E_\mu$ and $E_\rho$ into one operator based on the form of \eqref{eq:taupgeneral}; other methods of deriving an operator coarsening strategy for the pressure stabilisation operator are likely possible. \cref{algo:build2} summarises the operator coarsening strategy for this generalisation to the time-dependent Stokes equations, which, as before, is to be applied recursively to define all necessary operators down the mesh hierarchy.

\begin{algorithm}[t]
	\caption{\sffamily\small Construction of coarse mesh operators for the time-dependent multi-phase Stokes equations, given fine mesh operators $M_h$, $M_{\mu,h}$, $M_{\rho,h}$, $G_h$, ${\mathcal G}_h$, $\tilde{E}_h$, $E_{\mu,h}$ and $E_{\rho,h}$.}
	\begin{algorithmic}[1]
		\State $M_{2h} := (I_{2h}^h)^\trans M_h I_{2h}^h$
		\State $M_{\mu,2h} := (I_{2h}^h)^\trans M_{\mu,h} I_{2h}^h$
		\State $M_{\rho,2h} := (I_{2h}^h)^\trans M_{\rho,h} I_{2h}^h$
		\State $G_{2h} := M_{2h}^{-1} (I_{2h}^h)^\trans M_h G_h I_{2h}^h$
		\State ${\mathcal G}_{2h} := M_{2h}^{-1} (I_{2h}^h)^\trans M_h {\mathcal G}_h I_{2h}^h$; (equivalently, ${\mathcal G}_{2h,i} := -M_{2h}^{-1} G_{2h,i}^\trans M_{2h}$)
		\State $\tilde{E}_{2h} := \tfrac12 (I_{2h}^h)^\trans \tilde{E}_h I_{2h}^h$
		\State $E_{\mu,2h} := 2 (I_{2h}^h)^\trans E_{\mu,h} I_{2h}^h$
		\State $E_{\rho,2h} := \tfrac12 (I_{2h}^h)^\trans E_{\rho,h} I_{2h}^h$
		\State Combine $E_{\mu,2h}$ and $E_{\rho,2h}$ into one effective pressure stabilisation operator, $E_{2h}$, via \eqref{eq:Etimedep}
		\State Build the coarsened viscous operator in $d \times d$ block form, with $(i,j)$th block given by
		\[ A_{2h,ij} := \delta_{ij} \bigl( \textstyle{\sum_{k=1}^d} G_{2h,k}^\trans M_{\mu,2h} G_{2h,k} \bigr) + \gamma\, G_{2h,j}^\trans M_{\mu,2h} G_{2h,i} + \delta_{ij} \tilde{E}_{2h} \]
		\State Form the coarsened Stokes operator,
		\[ {\mathcal A}_{2h} := \begin{pmatrix} \tfrac{1}{\delta} M_{\rho,2h} + A_{2h} & M_{2h} {\mathcal G}_{2h} \\ {\mathcal G}_{2h}^\trans M_{2h} & -E_{2h} \end{pmatrix} \]
	\end{algorithmic}
	\label{algo:build2}
\end{algorithm}

\section{Grid convergence analyses}
\label{sec:gridconv}

Throughout this work, grid convergence analyses have been used to measure the discretisation order of accuracy. In every test case, the same exact solution is used to define the source data $\mathbf f$, $f$, interfacial jump data ${\mathbf g}_{ij}$, ${\mathbf h}_{ij}$, and boundary data ${\mathbf g}_\partial$, ${\mathbf h}_\partial$, in the governing equations. This exact solution is based on a simple smooth sinusoidal wave, for both velocity and pressure, with a translation depending on the component (and phase in the multi-phase case) to avoid any coincidental alignment with the mesh or cancellation in jump data. It is defined by the velocity field $\vu : \Omega \to {\mathbb R}^d$, $\vu = (u_1,\ldots,u_d)$, and pressure field $p : \Omega \to {\mathbb R}$, where
\[ u_i(x) = \prod_{j=1}^d \sin 2 \pi \bigl( x_j - 0.2 i - 0.25 (\chi - 1) \bigr), \]
and
\[ p(x) = \mu_\chi \prod_{j=1}^d \sin 2 \pi \bigl(x_j + 0.2 - 0.25 (\chi - 1) \bigr), \]
where the term $0.2 i$ represents a component-dependent shift, $\chi$ is the phase (i.e., $x \in \Omega_\chi$), and $\mu_\chi$ is a typical viscosity coefficient of phase $\chi$.

\begin{figure}%[tbhp]
\centering\sffamily\footnotesize%
\includegraphics[scale=0.91]{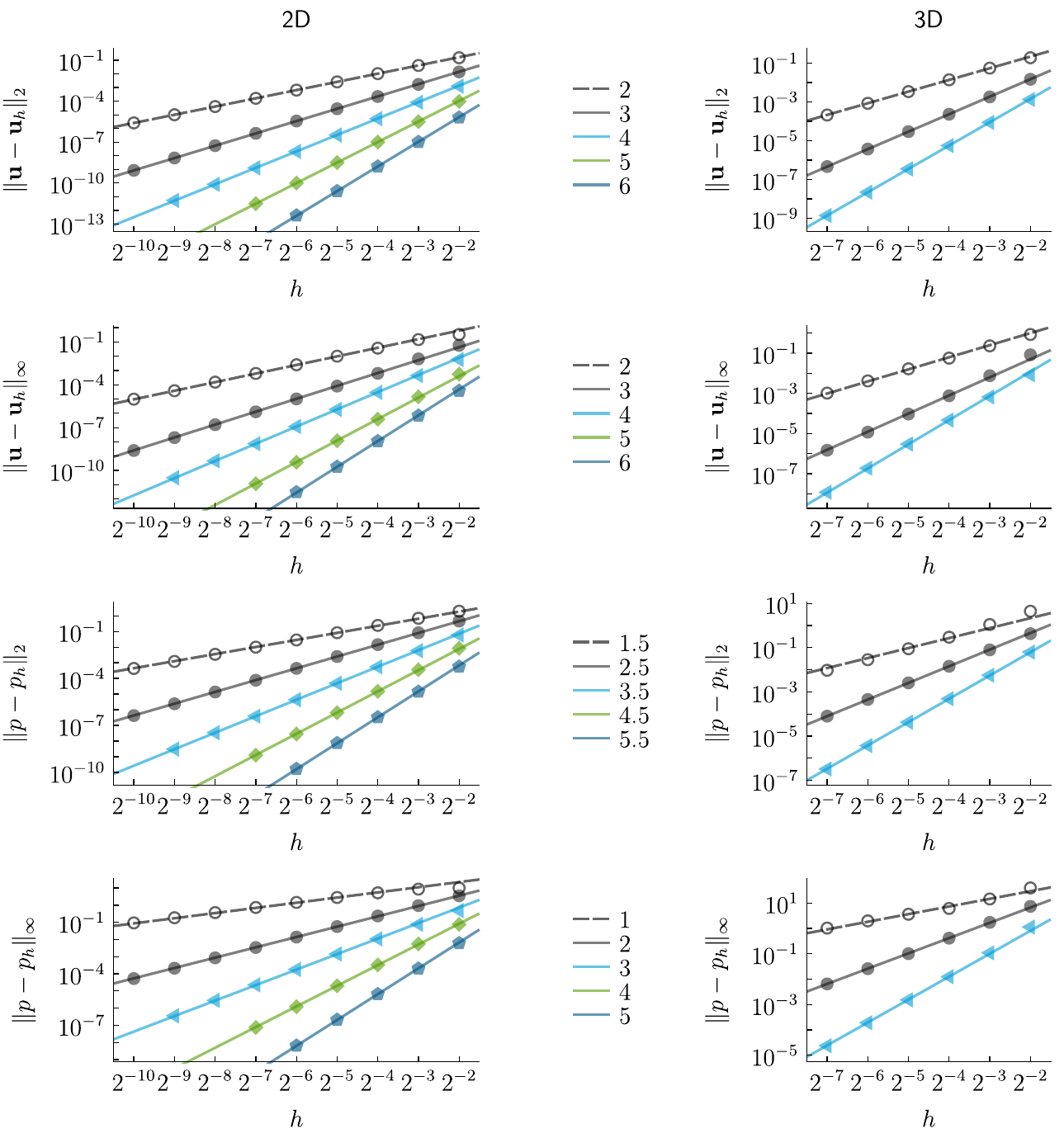}%
\caption[]{Grid convergence analysis for the single-phase Stokes problem in standard form, with $\mu = 1$ and velocity Dirichlet boundary conditions (see \cref{prob:dirichlet}). (Results for the test problem considered in \cref{prob:stress}, i.e., with stress boundary conditions, have similar characteristics.) Here, $h$ denotes the mesh element size and the lines of indicated slope illustrate the asymptotic convergence rate in the corresponding error norm, e.g., a slope of 4 indicates 4th order accuracy. Polynomial degrees are symbolised by $\circ$, $\bullet$, $\trianglesymbol$, $\squaresymbol$, and $\pentagonsymbol$ for $p = 1,2,3,4,$ and $5$, respectively.}%
\label{fig:si-grid-dirichlet-stokes}%
\end{figure}

\begin{figure}%[tbhp]
\centering\sffamily\footnotesize%
\includegraphics[scale=0.91]{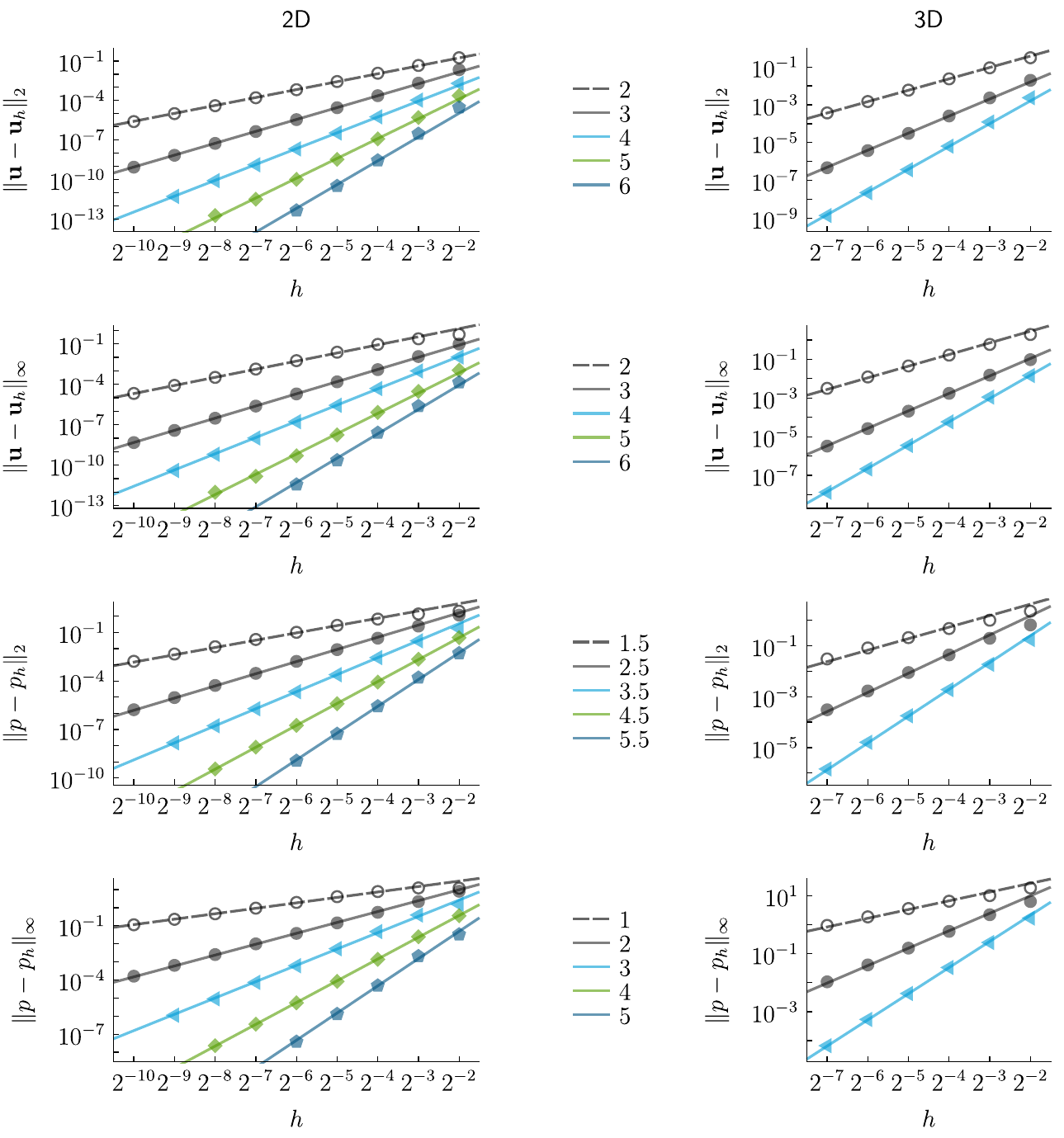}%
\caption[]{Grid convergence analysis for the single-phase Stokes problem in viscous-stress form, with non-constant variable viscosity and stress boundary conditions (see \cref{prob:variableviscosity}). Here, $h$ denotes the mesh element size and the lines of indicated slope illustrate the asymptotic convergence rate in the corresponding error norm, e.g., a slope of 4 indicates 4th order accuracy. Polynomial degrees are symbolised by $\circ$, $\bullet$, $\trianglesymbol$, $\squaresymbol$, and $\pentagonsymbol$ for $p = 1,2,3,4,$ and $5$, respectively.}%
\label{fig:si-grid-variable-mu-stress}%
\end{figure}

\begin{figure}%[tbhp]
\centering\sffamily\footnotesize%
\includegraphics[scale=0.91]{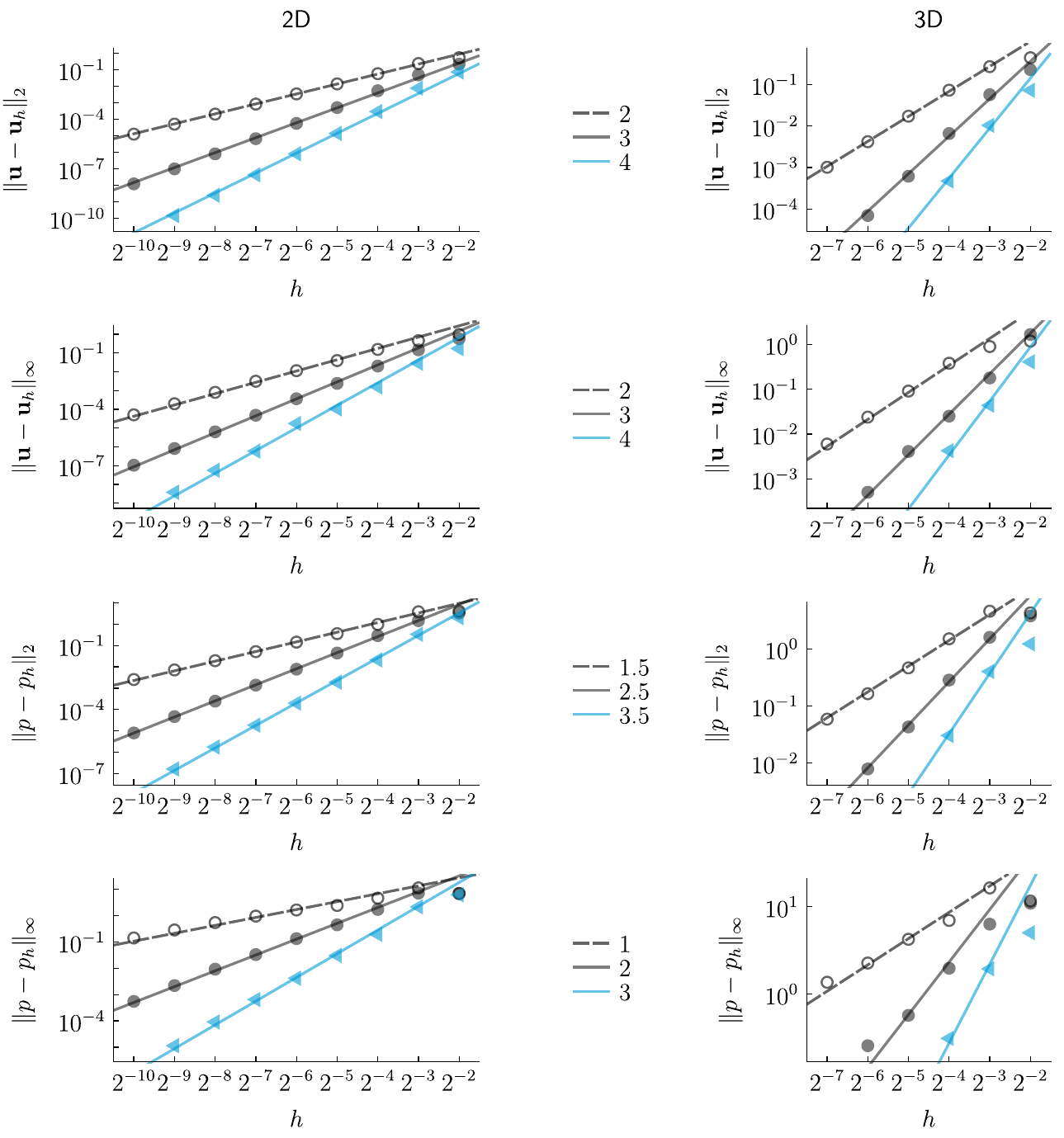}%
\caption[]{Grid convergence analysis for the single-phase Stokes problem in standard form, with $\mu = 1$ and velocity Dirichlet boundary conditions, in a unit diameter spherical domain using implicitly defined meshes (see \cref{prob:sphere}). Here, $h$ denotes the typical mesh element size and the lines of indicated slope illustrate the asymptotic convergence rate in the corresponding error norm, e.g., a slope of 4 indicates 4th order accuracy. Polynomial degrees are symbolised by $\circ$, $\bullet$, and $\trianglesymbol$ for $p = 1,2,$ and $3$, respectively.}%
\label{fig:si-grid-sphere}%
\end{figure}

\begin{figure}%[tbhp]
\centering\sffamily\footnotesize%
\includegraphics[scale=0.91]{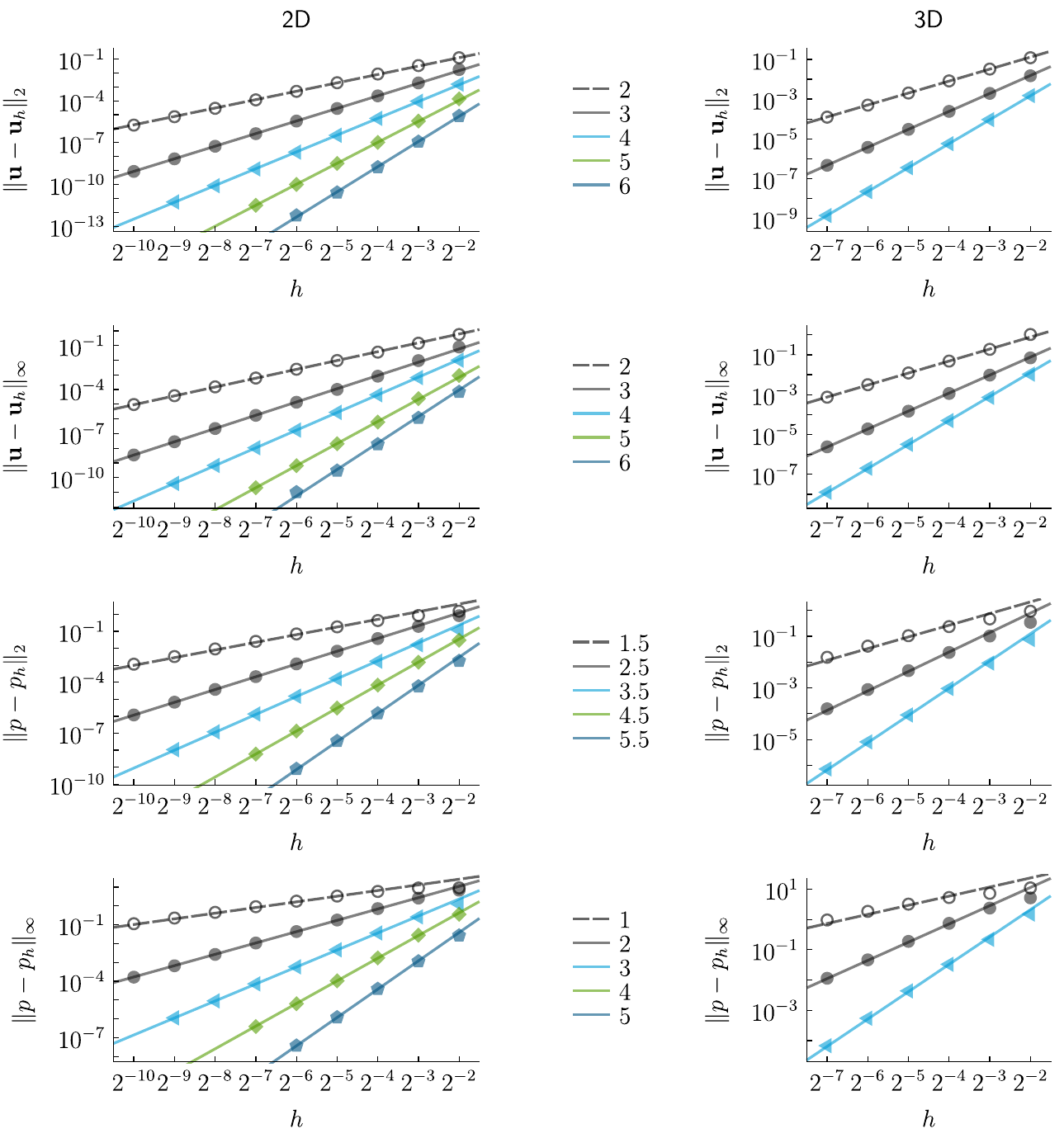}%
\caption[]{Grid convergence analysis for the multi-phase Stokes problem in viscous-stress form, in which $\Omega_1 = (\tfrac14, \tfrac34)^d$ has viscosity $\mu_1 = 10^{-6}$, and $\Omega_2 = (0,1)^d \setminus \overline{\Omega_1}$ has viscosity $\mu_2 = 1$, with stress boundary conditions (see \cref{prob:square}). (Results for the other viscosity ratios considered, i.e., $\mu_1 \in \{10^{-3},10^{+3},10^{+6}\}$, have similar characteristics.) Here, $h$ denotes the mesh element size and the lines of indicated slope illustrate the asymptotic convergence rate in the corresponding error norm, e.g., a slope of 4 indicates 4th order accuracy. Polynomial degrees are symbolised by $\circ$, $\bullet$, $\trianglesymbol$, $\squaresymbol$, and $\pentagonsymbol$ for $p = 1,2,3,4,$ and $5$, respectively.}%
\label{fig:si-grid-square-6}%
\end{figure}

\begin{figure}%[tbhp]
\centering\sffamily\footnotesize%
\includegraphics[scale=0.91]{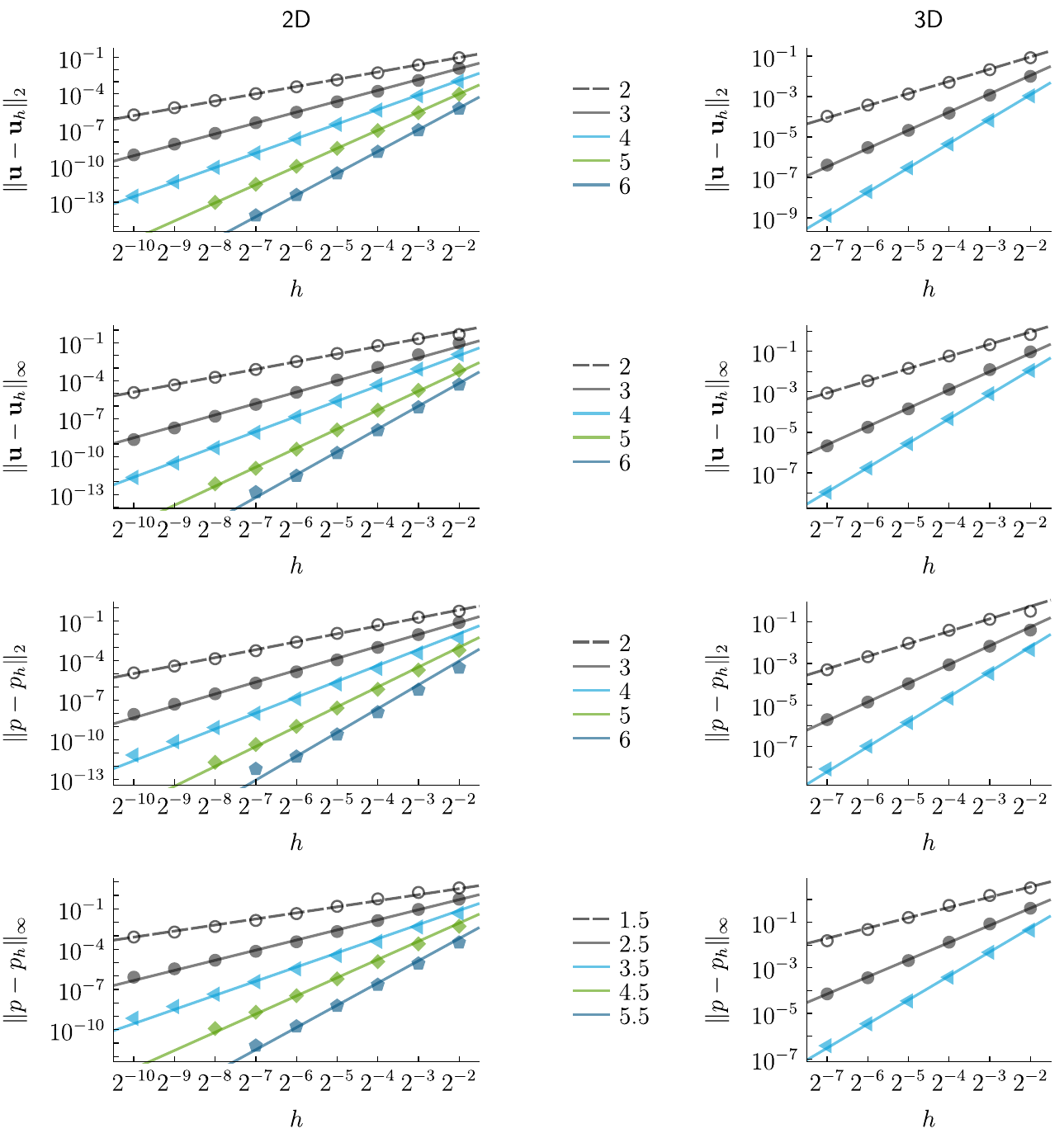}%
\caption[]{Grid convergence analysis for the time-dependent single-phase Stokes problem in standard form, with $\mu = 10^{-2}$ and $\delta = 0.1 h$, corresponding to $\textsf{Re} \approx 100$ (see \cref{prob:timedep1}). Here, $h$ denotes the mesh element size and the lines of indicated slope illustrate the asymptotic convergence rate in the corresponding error norm, e.g., a slope of 4 indicates 4th order accuracy. Polynomial degrees are symbolised by $\circ$, $\bullet$, $\trianglesymbol$, $\squaresymbol$, and $\pentagonsymbol$ for $p = 1,2,3,4,$ and $5$, respectively.}%
\label{fig:si-grid-time-single-phase-2}%
\end{figure}

\begin{figure}%[tbhp]
\centering\sffamily\footnotesize%
\includegraphics[scale=0.91]{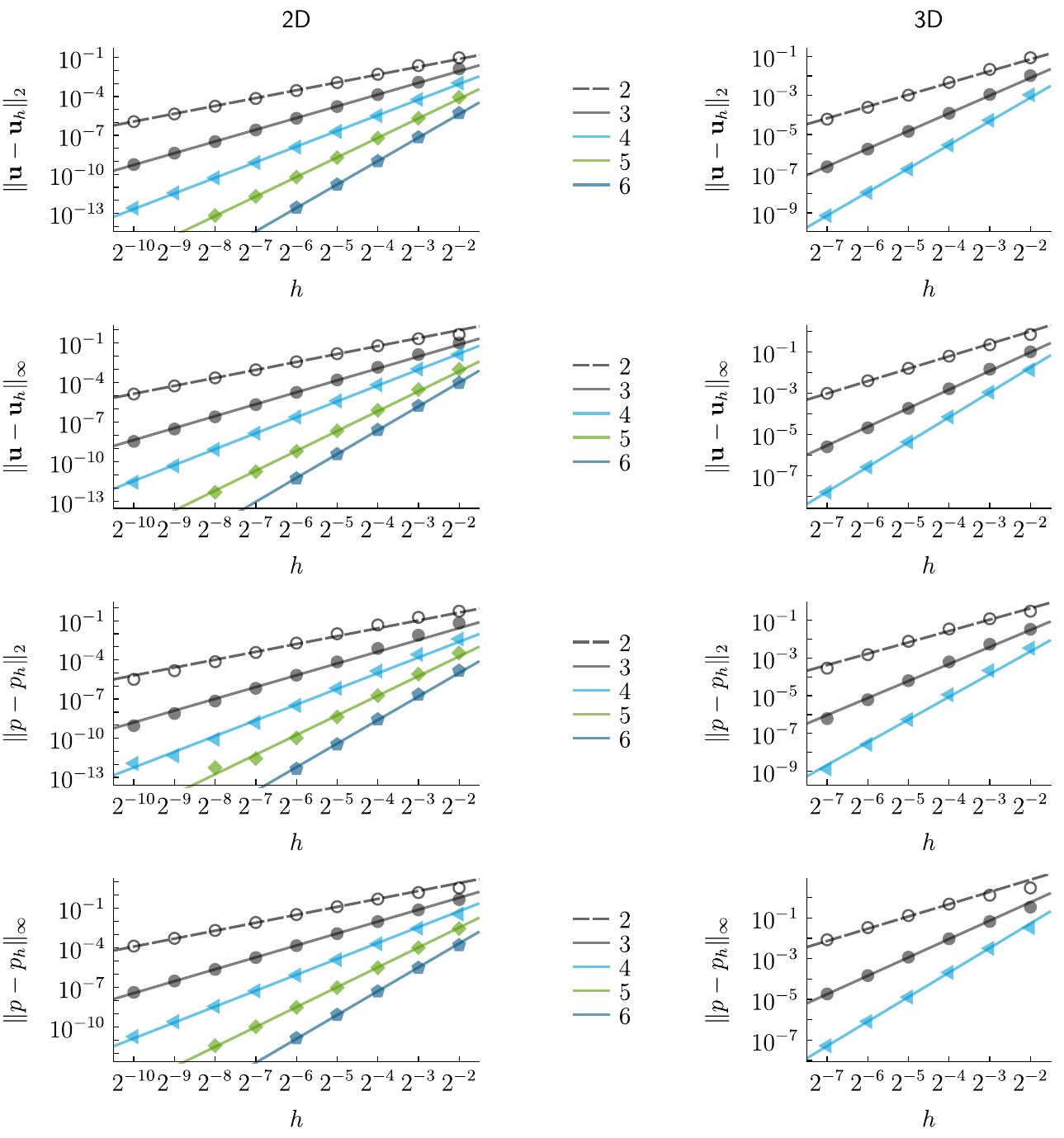}%
\caption[]{Grid convergence analysis for the time-dependent single-phase Stokes problem in standard form, with $\mu = 10^{-4}$ and $\delta = 0.1 h$, corresponding to $\textsf{Re} \approx 10,\!000$ (see \cref{prob:timedep1}). Here, $h$ denotes the mesh element size and the lines of indicated slope illustrate the asymptotic convergence rate in the corresponding error norm, e.g., a slope of 4 indicates 4th order accuracy. Polynomial degrees are symbolised by $\circ$, $\bullet$, $\trianglesymbol$, $\squaresymbol$, and $\pentagonsymbol$ for $p = 1,2,3,4,$ and $5$, respectively.}%
\label{fig:si-grid-time-single-phase-4}%
\end{figure}

\begin{figure}%[tbhp]
\centering\sffamily\footnotesize%
\includegraphics[scale=0.91]{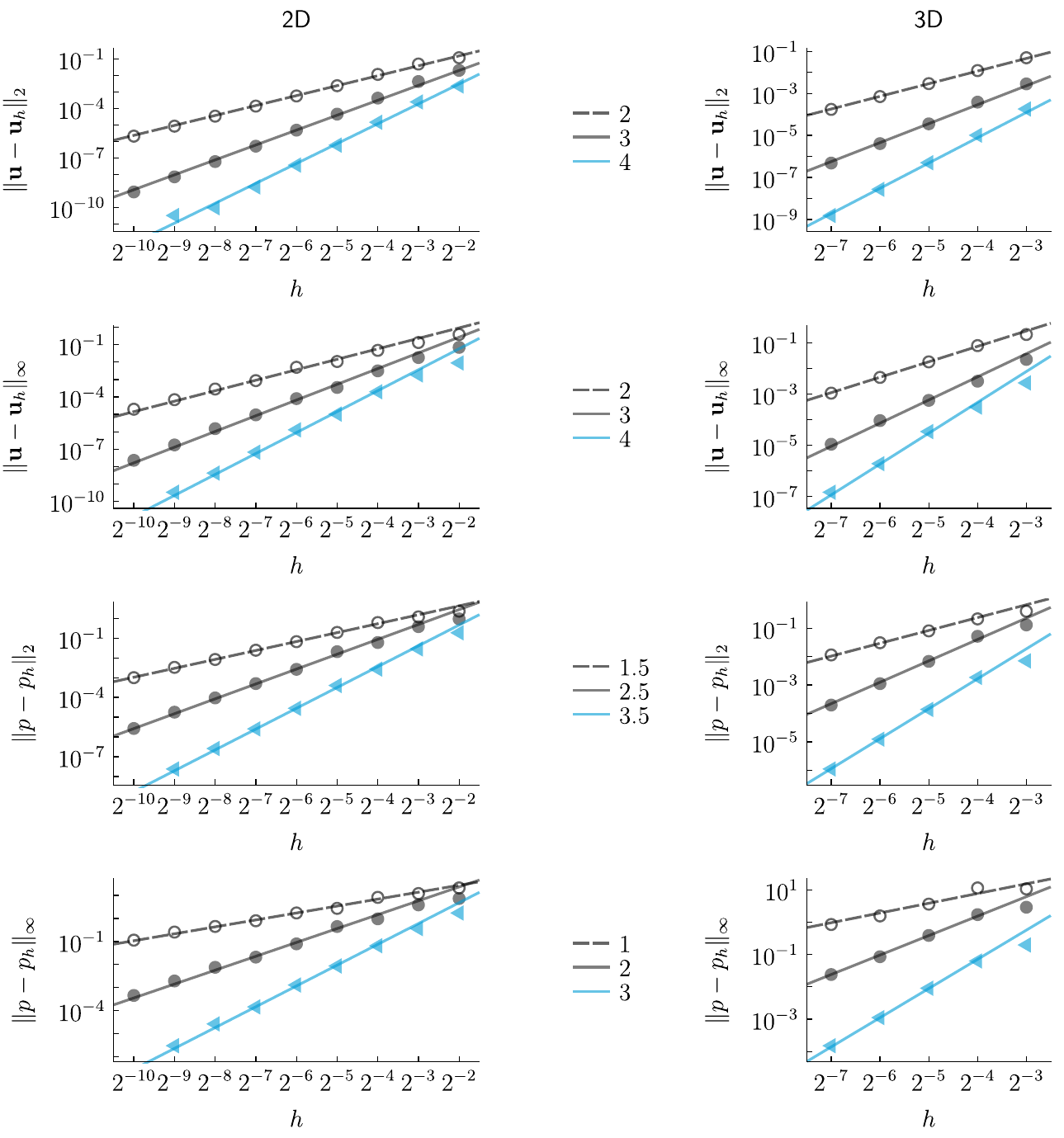}%
\caption[]{Grid convergence analysis for the time-dependent multi-phase Stokes problem in viscous-stress form, for a water bubble surrounded by gas, with $\rho_{\text{water}} = 1$, $\rho_{\text{gas}} = 0.001$, $\mu_{\text{water}} = 1$, $\mu_{\text{gas}} = 0.0002$, and $\delta = 0.1 h$, together with velocity Dirichlet boundary conditions (see \cref{prob:timedep2}). (Results for the opposite case, i.e., a gas bubble surrounded by water, have similar characteristics.) Here, $h$ denotes the typical mesh element size and the lines of indicated slope illustrate the asymptotic convergence rate in the corresponding error norm, e.g., a slope of 4 indicates 4th order accuracy. Polynomial degrees are symbolised by $\circ$, $\bullet$, and $\trianglesymbol$ for $p = 1,2,$ and $3$, respectively.}%
\label{fig:si-grid-water-bubble}%
\end{figure}

Whether the continuum (steady-state) Stokes problem is posed with velocity Dirichlet boundary conditions, or stress boundary conditions, or periodic boundary conditions, there is always an associated kernel of dimension at least one.\footnote{The dimension of the kernel may be smaller for time-dependent Stokes problems; for example, if Neumann or stress boundary conditions are applied, the kernel is zero-dimensional.} This kernel is referred to as the ``trivial kernel'', for it consists of constant-valued velocity fields and/or constant pressure fields, and possibly additional modes in the case of the viscous-stress form with stress boundary conditions (e.g., the velocity field $(x,y) \mapsto (-y,x)$, for which $\nabla \vu + \nabla \vu^\trans$ is zero). Since the Stokes operator is symmetric, it follows that the right hand side data $({\mathbf b}_\vu,b_p)$ must be orthogonal to the kernel; this is always the case for the method of manufactured solution applied here. However, the continuum solution is only unique up to modes in the kernel, and the discrete solution computed via the multigrid preconditioned GMRES method may contain arbitrary modes of the corresponding discrete kernel. To appropriately measure the discrete error, these modes are therefore disregarded. In particular, we compute the discrete error $(\vu - \vu_h, p - p_h)$ and nullify any components in the kernel through a simple Gram-Schmidt process applied to a basis of the kernel, known ahead of time. The resulting discrete error $(\vu - \vu_h, p - p_h)$ is then measured in the $L^2$ norm and the maximum norm, and is reported in the following collection of graphs, \crefrange{fig:si-grid-periodic-stokes}{fig:si-grid-water-bubble}. In each case, data points represent the measured error, and the lines of indicated slop are plotted to illustrate the asymptotic convergence rate. In some cases, the discrete error is saturated by numerical conditioning associated with double-precision arithmetic, forcing the cessation of high-order convergence; these data points are excluded from the graphs.

\subsection{Impact of pressure penalty parameter on discretisation error}
\label{app:tauaccuracy}

\begin{figure}%[tbhp]
\centering
\includegraphics[scale=0.91]{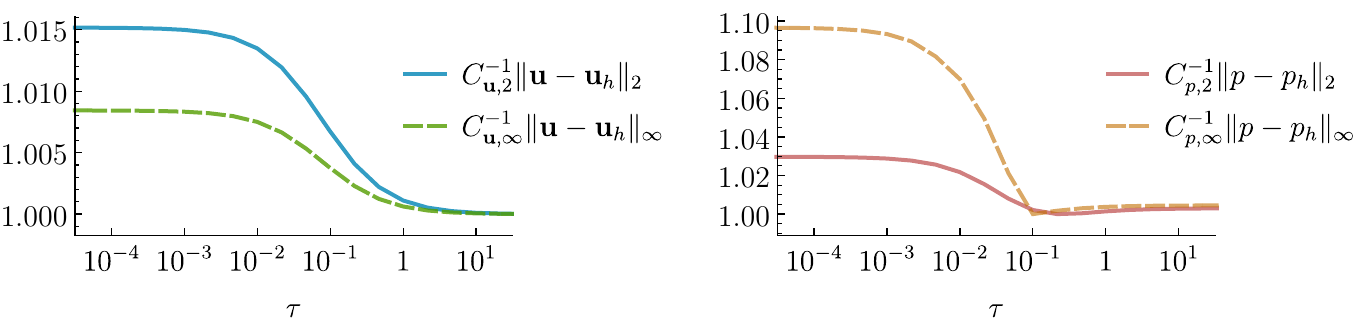}
\caption{Discretisation accuracy as a function of pressure penalty stabilisation parameter. Here, the errors in velocity (left) and pressure (right) are normalised by their minimum attained value.}
\label{fig:tauaccuracy}
\end{figure}

To supplement the discussion in \cref{sec:tauinfluence}, shown here is a test case examining the impact of the pressure penalty parameter $\tau$ on discretisation error. We consider a single-phase Stokes problem in standard form, with $\mu = 1$ and periodic boundary conditions, on a uniform Cartesian grid with $p = 2$ biquadratic elements. In the following tests, we consider extreme values of $\tau$ and utilise a direct solver instead of multigrid (thereby eliminating possible issues of multigrid non-convergence associated with extreme values of $\tau$); as such, a coarser grid of $16 \times 16$ is used. With the same range of $\tau$ as used in the discussion of \cref{sec:tauinfluence}, \cref{fig:tauaccuracy} shows the discrete error in the velocity and pressure in the $L^2$ and maximum norms; in each case the error is normalised by the minimum attained value, e.g., $C_{{\mathbf u},2} = \min_\tau \|\mathbf u - {\mathbf u}_h\|_2$, and similarly for the other quantities shown. (The order of magnitude of these errors can be inferred from the 2D results of \cref{fig:si-grid-periodic-stokes} below.) According to \cref{fig:tauaccuracy}, note that $\tau$ has relatively little influence on the error in both velocity and pressure, with maximal discretisation errors at most 10\% greater than optimal. Although only a minor improvement, note also that the best error in pressure is attained when $\tau$ is approximately equal or larger to the optimal $\tau$ values found in \cref{tab:optimaltau}. Overall, this kind of behaviour has been observed across many of the examples considered in this work, leading to the conclusion that, at least for relatively well-resolved Stokes problems, the discretisation error is largely insensitive to the value of $\tau$, thereby allowing this parameter to be tuned for excellent multigrid performance.

\begin{figure}%[tbhp]
\centering\sffamily\footnotesize%
\includegraphics[scale=0.91]{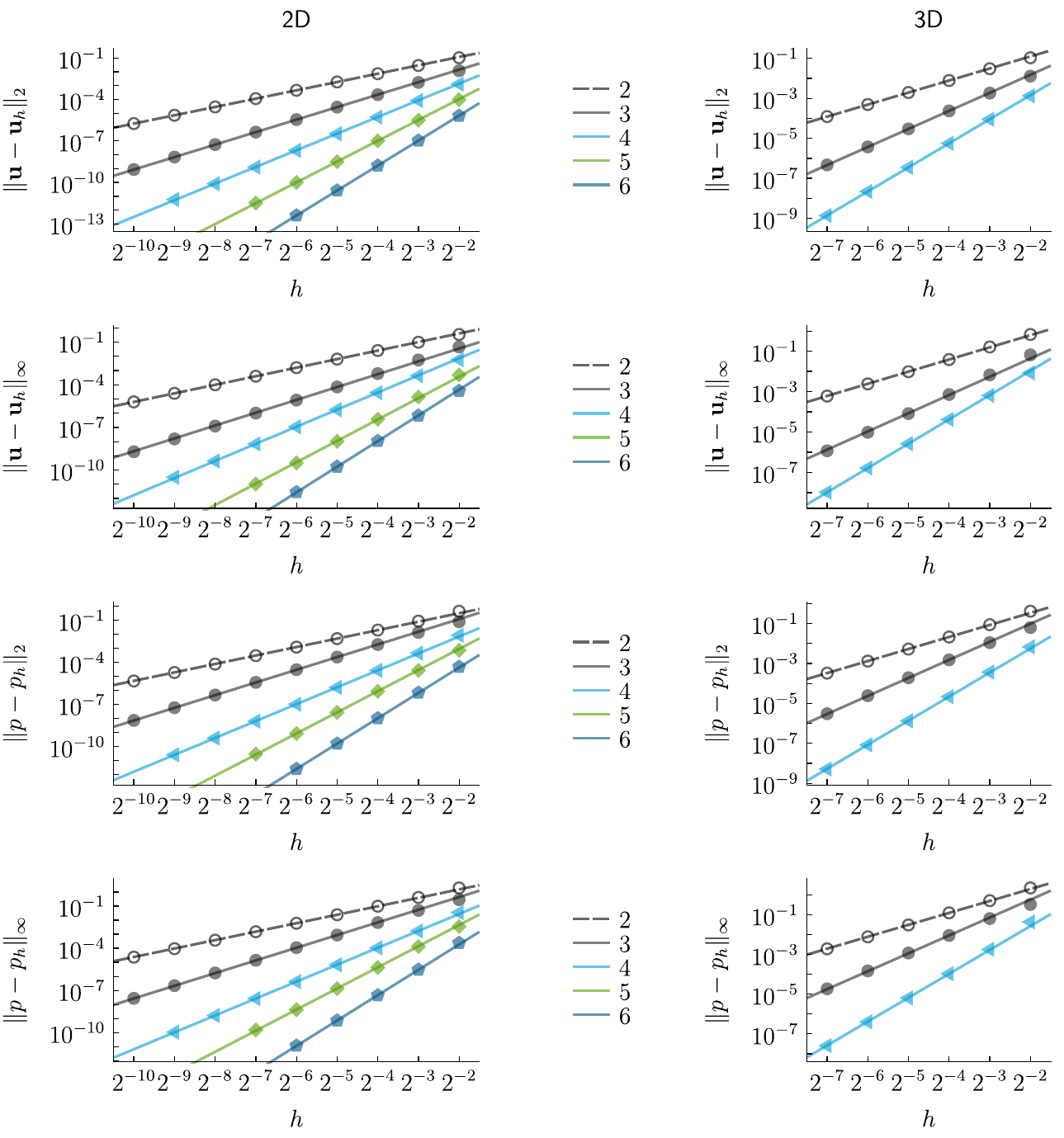}%
\caption[]{Grid convergence analysis for the single-phase Stokes problem in standard form, with $\mu = 1$ and periodic boundary conditions (see \cref{prob:periodic}). Here, $h$ denotes the mesh element size and the lines of indicated slope illustrate the asymptotic convergence rate in the corresponding error norm, e.g., a slope of 4 indicates 4th order accuracy. Polynomial degrees are symbolised by $\circ$, $\bullet$, $\trianglesymbol$, $\squaresymbol$, and $\pentagonsymbol$ for $p = 1,2,3,4,$ and $5$, respectively.}%
\label{fig:si-grid-periodic-stokes}%
\end{figure}

\clearpage

\bibliographystyle{siamplain}
\bibliography{references}

\end{document}

%% file: fluxdiagram.pdf_tex
%% Creator: Inkscape inkscape 0.92.2, www.inkscape.org
%% Accompanies image file 'fluxdiagram.pdf' (pdf, eps, ps)
\begingroup%
  \makeatletter%
  \providecommand\color[2][]{%
    \errmessage{(Inkscape) Color is used for the text in Inkscape, but the package 'color.sty' is not loaded}%
    \renewcommand\color[2][]{}%
  }%
  \providecommand\transparent[1]{%
    \errmessage{(Inkscape) Transparency is used (non-zero) for the text in Inkscape, but the package 'transparent.sty' is not loaded}%
    \renewcommand\transparent[1]{}%
  }%
  \providecommand\rotatebox[2]{#2}%
  \ifx\svgwidth\undefined%
    \setlength{\unitlength}{1332.13595581bp}%
    \ifx\svgscale\undefined%
      \relax%
    \else%
      \setlength{\unitlength}{\unitlength * \real{\svgscale}}%
    \fi%
  \else%
    \setlength{\unitlength}{\svgwidth}%
  \fi%
  \global\let\svgwidth\undefined%
  \global\let\svgscale\undefined%
  \makeatother%
  \begin{picture}(1,0.24265541)%
    \put(0,0){\includegraphics[width=\unitlength,page=1]{fluxdiagram.pdf}}%
    \put(0.11528163,0.22042016){\color[rgb]{0,0,0}\makebox(0,0)[b]{\smash{Intraphase face}}}%
    \put(0.41975153,0.22046413){\color[rgb]{0,0,0}\makebox(0,0)[b]{\smash{Interphase face on $\Gamma_{ij}$, $i<j$}}}%
    \put(0.66724325,0.22042016){\color[rgb]{0,0,0}\makebox(0,0)[b]{\smash{Stress boundary face}}}%
    \put(0.89244548,0.22042016){\color[rgb]{0,0,0}\makebox(0,0)[b]{\smash{Dirichlet boundary face}}}%
    \put(0,0){\includegraphics[width=\unitlength,page=2]{fluxdiagram.pdf}}%
    \put(0.05944506,0.14689979){\color[rgb]{0,0,0}\makebox(0,0)[b]{\smash{$-$}}}%
    \put(0.17204618,0.14689979){\color[rgb]{0,0,0}\makebox(0,0)[b]{\smash{$+$}}}%
    \put(0.36578114,0.14689979){\color[rgb]{0,0,0}\makebox(0,0)[b]{\smash{$-$}}}%
    \put(0.4781392,0.14689979){\color[rgb]{0,0,0}\makebox(0,0)[b]{\smash{$+$}}}%
    \put(0.66749107,0.14689979){\color[rgb]{0,0,0}\makebox(0,0)[b]{\smash{$-$}}}%
    \put(0.89269331,0.14689979){\color[rgb]{0,0,0}\makebox(0,0)[b]{\smash{$-$}}}%
    \put(0.13625725,0.17395976){\color[rgb]{0,0,0}\makebox(0,0)[b]{\smash{$\vn$}}}%
    \put(0.44123553,0.17384482){\color[rgb]{0,0,0}\makebox(0,0)[b]{\smash{$\vn$}}}%
    \put(0.74300083,0.13565215){\color[rgb]{0,0,0}\makebox(0,0)[b]{\smash{$\vn$}}}%
    \put(0.96887068,0.13565215){\color[rgb]{0,0,0}\makebox(0,0)[b]{\smash{$\vn$}}}%
    \put(0,0){\includegraphics[width=\unitlength,page=3]{fluxdiagram.pdf}}%
    \put(0.74673927,0.19017179){\color[rgb]{0,0,0}\makebox(0,0)[lb]{\smash{$\Gamma_N$}}}%
    \put(0.36352911,0.12875533){\color[rgb]{0,0,0}\makebox(0,0)[b]{\smash{$\chi=i\vphantom{j}$}}}%
    \put(0.47708756,0.12875533){\color[rgb]{0,0,0}\makebox(0,0)[b]{\smash{$\chi=j$}}}%
    \put(0,0){\includegraphics[width=\unitlength,page=4]{fluxdiagram.pdf}}%
    \put(0.9719415,0.19017179){\color[rgb]{0,0,0}\makebox(0,0)[lb]{\smash{$\Gamma_D$}}}%
    \put(0,0){\includegraphics[width=\unitlength,page=5]{fluxdiagram.pdf}}%
    \put(0.12227874,0.06146026){\color[rgb]{0,0,0}\makebox(0,0)[b]{\smash{$\begin{aligned}\hat \vu&=\vu^-\\[-2em]\\ {\hat\vsigma}&={\vsigma}^+\end{aligned}$}}}%
    \put(0.3399082,0.06146026){\color[rgb]{0,0,0}\makebox(0,0)[b]{\smash{$\begin{aligned}{\hat \vu}_i&=\lambda \vu^- + (1-\lambda)(\vu^+ + {\mathbf g}_{ij})\\[-2em]\\             {\hat\vsigma}_i&=(1-\lambda){\vsigma}^-+\lambda ({\vsigma}^+ + {\mathbf h}_{ij} \otimes \vn)\end{aligned}$}}}%
    \put(0.57910607,0.0063445){\color[rgb]{0,0,0}\makebox(0,0)[b]{\smash{$\begin{aligned}{\hat\vu}_j&=\lambda (\vu^--{\mathbf g}_{ij}) + (1-\lambda)\vu^+\\[-2em]\\ {\hat\vsigma}_j&=(1-\lambda)({\vsigma}^--{\mathbf h}_{ij} \otimes \vn)+\lambda {\vsigma}^+\end{aligned}$}}}%
    \put(0.74593819,0.06146026){\color[rgb]{0,0,0}\makebox(0,0)[b]{\smash{$\begin{aligned}\hat\vu&=\vu^-\\[-2em]\\ {\hat\vsigma}&={\mathbf h}_\partial \otimes \vn\end{aligned}$}}}%
    \put(0.95756885,0.06146026){\color[rgb]{0,0,0}\makebox(0,0)[b]{\smash{$\begin{aligned}\hat \vu&={\mathbf g}_\partial\\[-2em]\\ {\hat\vsigma}&={\vsigma}^-\end{aligned}$}}}%
  \end{picture}%
\endgroup%

%% file: paper.bbl
\begin{thebibliography}{10}

\bibitem{AdlerBensonMacLachlan2016}
{\sc J.~H. Adler, T.~R. Benson, and S.~P. MacLachlan}, {\em Preconditioning a
  mass-conserving discontinuous {G}alerkin discretization of the stokes
  equations}, Numerical Linear Algebra with Applications, 24 (2017), p.~e2047,
  \url{https://doi.org/10.1002/nla.2047}.

\bibitem{AksoyluUnlu2014}
{\sc B.~Aksoylu and Z.~Unlu}, {\em Robust preconditioners for the high-contrast
  {S}tokes equation}, Journal of Computational and Applied Mathematics, 259
  (2014), pp.~944--954, \url{https://doi.org/10.1016/j.cam.2013.10.016}.

\bibitem{ArnoldBrezziCockburnMarini}
{\sc D.~N. Arnold, F.~Brezzi, B.~Cockburn, and L.~D. Marini}, {\em Unified
  analysis of discontinuous {G}alerkin methods for elliptic problems}, SIAM
  Journal on Numerical Analysis, 39 (2002), pp.~1749--1779,
  \url{https://doi.org/10.1137/S0036142901384162}.

\bibitem{BauerKlementOberhuberZabka2016}
{\sc P.~Bauer, V.~Klement, T.~Oberhuber, and V.~{\v{Z}}abka}, {\em
  Implementation of the {V}anka-type multigrid solver for the finite element
  approximation of the {N}avier--{S}tokes equations on {GPU}}, Computer Physics
  Communications, 200 (2016), pp.~50--56,
  \url{https://doi.org/10.1016/j.cpc.2015.10.021}.

\bibitem{BenziGolubLiesen2005}
{\sc M.~Benzi, G.~H. Golub, and J.~Liesen}, {\em Numerical solution of saddle
  point problems}, Acta numerica, 14 (2005), pp.~1--137,
  \url{https://doi.org/10.1017/S0962492904000212}.

\bibitem{BorzacchielloLericheBlottiereGuillet2017}
{\sc D.~Borzacchiello, E.~Leriche, B.~Blotti{\`e}re, and J.~Guillet}, {\em
  Box-relaxation based multigrid solvers for the variable viscosity {S}tokes
  problem}, Computers \& Fluids, 156 (2017), pp.~515--525,
  \url{https://doi.org/10.1016/j.compfluid.2017.08.027}.

\bibitem{BrandtLivne}
{\sc A.~Brandt and O.~E. Livne}, {\em {Multigrid Techniques}}, Society for
  Industrial and Applied Mathematics, 2011,
  \url{https://doi.org/10.1137/1.9781611970753}.

\bibitem{Briggs_00_01}
{\sc W.~L. Briggs, V.~E. Henson, and S.~F. McCormick}, {\em {A Multigrid
  Tutorial, Second Edition}}, Society for Industrial and Applied Mathematics,
  2000, \url{https://doi.org/10.1137/1.9780898719505}.

\bibitem{Chen2015}
{\sc L.~Chen}, {\em Multigrid methods for saddle point systems using
  constrained smoothers}, Computers and Mathematics with Applications, 70
  (2015), pp.~2854--2866, \url{https://doi.org/10.1016/j.camwa.2015.09.020}.

\bibitem{Chorin}
{\sc A.~J. Chorin}, {\em Numerical solution of the {N}avier-{S}tokes
  equations}, Mathematics of Computation, 22 (1968), pp.~745--762,
  \url{https://doi.org/10.1090/S0025-5718-1968-0242392-2}.

\bibitem{Cockburn0701}
{\sc B.~Cockburn and B.~Dong}, {\em An analysis of the minimal dissipation
  local discontinuous {G}alerkin method for convection--diffusion problems}, J.
  Sci. Comput., 32 (2007), pp.~233--262,
  \url{https://doi.org/10.1007/s10915-007-9130-3}.

\bibitem{CockburnKanschatSchotzau2003}
{\sc B.~Cockburn, G.~Kanschat, and D.~Sch{\"o}tzau}, {\em {LDG} methods for
  {S}tokes flow problems}, in Numerical Mathematics and Advanced Applications,
  Springer, 2003, pp.~755--764,
  \url{https://doi.org/10.1007/978-88-470-2089-4_68}.

\bibitem{CockburnKanschatSchotzau2004}
{\sc B.~Cockburn, G.~Kanschat, and D.~Sch{\"o}tzau}, {\em A locally
  conservative {LDG} method for the incompressible {N}avier-{S}tokes
  equations}, Mathematics of Computation, 74 (2005), pp.~1067--1095,
  \url{https://doi.org/10.1090/S0025-5718-04-01718-1}.

\bibitem{CockburnKanschatSchotzauSchwab2002}
{\sc B.~Cockburn, G.~Kanschat, D.~Sch{\"o}tzau, and C.~Schwab}, {\em Local
  discontinuous {G}alerkin methods for the {S}tokes system}, SIAM Journal on
  Numerical Analysis, 40 (2002), pp.~319--343,
  \url{https://doi.org/10.1137/S0036142900380121}.

\bibitem{ldg}
{\sc B.~Cockburn and C.-W. Shu}, {\em The local discontinuous {G}alerkin method
  for time-dependent convection-diffusion systems}, SIAM Journal on Numerical
  Analysis, 35 (1998), pp.~2440--2463,
  \url{https://doi.org/10.1137/S0036142997316712}.

\bibitem{ColeyBenzakenEvans2017}
{\sc C.~Coley, J.~Benzaken, and J.~A. Evans}, {\em A geometric multigrid method
  for isogeometric compatible discretizations of the generalized {S}tokes and
  {O}seen problems}, Numerical Linear Algebra with Applications, 25 (2018),
  p.~e2145, \url{https://doi.org/10.1002/nla.2145}.

\bibitem{DrzisgaJohnRudeWohlmuthZulehner2018}
{\sc D.~Drzisga, L.~John, U.~Rude, B.~Wohlmuth, and W.~Zulehner}, {\em On the
  analysis of block smoothers for saddle point problems}, SIAM Journal on
  Matrix Analysis and Applications, 39 (2018), pp.~932--960,
  \url{https://doi.org/10.1137/16M1106304}.

\bibitem{ElsnerMhermann1991}
{\sc L.~Elsner and V.~Mehrmann}, {\em Convergence of block iterative methods
  for linear systems arising in the numerical solution of {E}uler equations},
  Numerische Mathematik, 59 (1991), pp.~541--559,
  \url{https://doi.org/10.1007/BF01385795}.

\bibitem{FarrellHeMacLachlan2019}
{\sc P.~E. Farrell, Y.~He, and S.~P. MacLachlan}, {\em A local {F}ourier
  analysis of additive {V}anka relaxation for the {S}tokes equations}, arXiv
  preprint arXiv:1908.09949,  (2019).

\bibitem{dgmg}
{\sc D.~Fortunato, C.~H. Rycroft, and R.~Saye}, {\em Efficient
  operator-coarsening multigrid schemes for local discontinuous {G}alerkin
  methods}, SIAM Journal on Scientific Computing, 41 (2019), pp.~A3913--A3937,
  \url{https://doi.org/10.1137/18M1206357}.

\bibitem{FuruichiMayTackley2011}
{\sc M.~Furuichi, D.~A. May, and P.~J. Tackley}, {\em Development of a {S}tokes
  flow solver robust to large viscosity jumps using a {S}chur complement
  approach with mixed precision arithmetic}, Journal of Computational Physics,
  230 (2011), pp.~8835--8851, \url{https://doi.org/10.1016/j.jcp.2011.09.007}.

\bibitem{GasparNotayOosterleeRodrigo2014}
{\sc F.~J. Gaspar, Y.~Notay, C.~W. Oosterlee, and C.~Rodrigo}, {\em A simple
  and efficient segregated smoother for the discrete {S}tokes equations}, SIAM
  journal on scientific computing, 36 (2014), pp.~A1187--A1206,
  \url{https://doi.org/10.1137/130920630}.

\bibitem{GmeinerHuberJohnRudeWohlmuth2016}
{\sc B.~Gmeiner, M.~Huber, L.~John, U.~R{\"u}de, and B.~Wohlmuth}, {\em A
  quantitative performance study for {S}tokes solvers at the extreme scale},
  Journal of Computational Science, 17 (2016), pp.~509--521,
  \url{https://doi.org/10.1016/j.jocs.2016.06.006}.

\bibitem{GuermondMinevShen}
{\sc J.~L. Guermond, P.~Minev, and J.~Shen}, {\em An overview of projection
  methods for incompressible flows}, Computer Methods in Applied Mechanics and
  Engineering, 195 (2006), pp.~6011--6045,
  \url{https://doi.org/10.1016/j.cma.2005.10.010}.

\bibitem{HeMacLachlan2017}
{\sc Y.~He and S.~P. MacLachlan}, {\em Local {F}ourier analysis of
  block-structured multigrid relaxation schemes for the {S}tokes equations},
  Numerical Linear Algebra with Applications, 25 (2018), p.~e2147,
  \url{https://doi.org/10.1002/nla.2147}.

\bibitem{HesthavenWarburton}
{\sc J.~S. Hesthaven and T.~Warburton}, {\em {Nodal Discontinuous Galerkin
  Methods: Algorithms, Analysis, and Applications}}, vol.~54 of Texts in
  Applied Mathematics, Springer, 2008,
  \url{https://doi.org/10.1007/978-0-387-72067-8}.

\bibitem{HongKrausXuZikatanov2016}
{\sc Q.~Hong, J.~Kraus, J.~Xu, and L.~Zikatanov}, {\em A robust multigrid
  method for discontinuous {G}alerkin discretizations of {S}tokes and linear
  elasticity equations}, Numerische Mathematik, 132 (2016), pp.~23--49,
  \url{https://doi.org/10.1007/s00211-015-0712-y}.

\bibitem{Kanschat2005}
{\sc G.~Kanschat}, {\em Block preconditioners for {LDG} discretizations of
  linear incompressible flow problems}, Journal of Scientific Computing, 22
  (2005), pp.~371--384, \url{https://doi.org/10.1007/s10915-004-4144-6}.

\bibitem{KanschatMao2015}
{\sc G.~Kanschat and Y.~Mao}, {\em Multigrid methods for {H}div-conforming
  discontinuous {G}alerkin methods for the {S}tokes equations}, Journal of
  Numerical Mathematics, 23 (2015), pp.~51--66,
  \url{https://doi.org/10.1515/jnma-2015-0005}.

\bibitem{Manservisi2006}
{\sc S.~Manservisi}, {\em Numerical analysis of {V}anka-type solvers for steady
  {S}tokes and {N}avier--{S}tokes flows}, SIAM Journal on Numerical Analysis,
  44 (2006), pp.~2025--2056, \url{https://doi.org/10.1137/060655407}.

\bibitem{sdpc}
{\sc M.~L. Minion and R.~Saye}, {\em {Higher-order temporal integration for the
  incompressible Navier--Stokes equations in bounded domains}}, Journal of
  Computational Physics, 375 (2018), pp.~797--822,
  \url{https://doi.org/10.1016/j.jcp.2018.08.054}.

\bibitem{Notay2017}
{\sc Y.~Notay}, {\em Algebraic multigrid for {S}tokes equations}, SIAM Journal
  on Scientific Computing, 39 (2017), pp.~S88--S111,
  \url{https://doi.org/10.1137/16M1071419}.

\bibitem{OlshanskiiTyrtyshnikov}
{\sc M.~Olshanskii and E.~Tyrtyshnikov}, {\em {Iterative Methods for Linear
  Systems: Theory and Applications}}, Society for Industrial and Applied
  Mathematics, 2014.

\bibitem{OosterleeGaspar2008}
{\sc C.~Oosterlee and F.~Gaspar}, {\em Multigrid relaxation methods for systems
  of saddle point type}, Applied Numerical Mathematics, 58 (2008),
  pp.~1933--1950, \url{https://doi.org/10.1016/j.apnum.2007.11.014}.

\bibitem{OosterleeWesseling1993}
{\sc C.~W. Oosterlee and P.~Wesseling}, {\em A robust multigrid method for a
  discretization of the incompressible {N}avier-{S}tokes equations in general
  coordinates}, IMPACT of Computing in Science and Engineering, 5 (1993),
  pp.~128--151, \url{https://doi.org/10.1006/icse.1993.1006}.

\bibitem{OosterleeWesseling1993b}
{\sc C.~W. Oosterlee and P.~Wesseling}, {\em Steady incompressible flow around
  objects in general coordinates with a multigrid solution method}, Numerical
  Methods for Partial Differential Equations, 10 (1994), pp.~295--308,
  \url{https://doi.org/10.1002/num.1690100304}.

\bibitem{ImplicitMeshPartOne}
{\sc R.~Saye}, {\em {Implicit mesh discontinuous Galerkin methods and
  interfacial gauge methods for high-order accurate interface dynamics, with
  applications to surface tension dynamics, rigid body fluid-structure
  interaction, and free surface flow: Part I}}, Journal of Computational
  Physics, 344 (2017), pp.~647--682,
  \url{https://doi.org/10.1016/j.jcp.2017.04.076}.

\bibitem{ImplicitMeshPartTwo}
{\sc R.~Saye}, {\em {Implicit mesh discontinuous Galerkin methods and
  interfacial gauge methods for high-order accurate interface dynamics, with
  applications to surface tension dynamics, rigid body fluid-structure
  interaction, and free surface flow: Part II}}, Journal of Computational
  Physics, 344 (2017), pp.~683--723,
  \url{https://doi.org/10.1016/j.jcp.2017.05.003}.

\bibitem{algoim}
{\sc R.~Saye}, {\em {Algoim -- Algorithms for implicitly defined geometry,
  level set methods, and Voronoi implicit interface methods}}.
\newblock \url{https://algoim.github.io/}, 2019.

\bibitem{HighorderImplicitQuad}
{\sc R.~I. Saye}, {\em {High-Order Quadrature Methods for Implicitly Defined
  Surfaces and Volumes in Hyperrectangles}}, SIAM Journal on Scientic
  Computing, 37 (2015), pp.~A993--A1019,
  \url{https://doi.org/10.1137/140966290}.

\bibitem{fluxx}
{\sc R.~I. Saye}, {\em Efficient multigrid solution of elliptic interface
  problems using viscosity-upwinded local discontinuous {G}alerkin methods},
  Communications in Applied Mathematics and Computational Science, 14 (2019),
  pp.~247--283, \url{https://doi.org/10.2140/camcos.2019.14.247}.

\bibitem{SchoberlZulehner2003}
{\sc J.~Sch{\"o}berl and W.~Zulehner}, {\em On {S}chwarz-type smoothers for
  saddle point problems}, Numerische Mathematik, 95 (2003), pp.~377--399,
  \url{https://doi.org/10.1007/s00211-002-0448-3}.

\bibitem{Sivaloganathan1991}
{\sc S.~Sivaloganathan}, {\em The use of local mode analysis in the design and
  comparison of multigrid methods}, Computer Physics Communications, 65 (1991),
  pp.~246--252, \url{https://doi.org/10.1016/0010-4655(91)90178-N}.

\bibitem{ThompsonFerziger1989}
{\sc M.~Thompson and J.~H. Ferziger}, {\em An adaptive multigrid technique for
  the incompressible {N}avier-{S}tokes equations}, Journal of computational
  Physics, 82 (1989), pp.~94--121,
  \url{https://doi.org/10.1016/0021-9991(89)90037-5}.

\bibitem{Vanka1986}
{\sc S.~P. Vanka}, {\em Block-implicit multigrid solution of {N}avier-{S}tokes
  equations in primitive variables}, Journal of Computational Physics, 65
  (1986), pp.~138--158, \url{https://doi.org/10.1016/0021-9991(86)90008-2}.

\bibitem{Vanka1986b}
{\sc S.~P. Vanka}, {\em A calculation procedure for three-dimensional steady
  recirculating flows using multigrid methods}, Computer Methods in Applied
  Mechanics and Engineering, 55 (1986), pp.~321--338,
  \url{https://doi.org/10.1016/0045-7825(86)90058-7}.

\bibitem{Wesseling2001}
{\sc P.~Wesseling}, {\em Principles of computational fluid dynamics}, vol.~29
  of Computational Mathematics, Springer, Berlin, Heidelberg, 2001,
  \url{https://doi.org/10.1007/978-3-642-05146-3}.

\bibitem{WesselingOosterlee2001}
{\sc P.~Wesseling and C.~W. Oosterlee}, {\em Geometric multigrid with
  applications to computational fluid dynamics}, Journal of Computational and
  Applied Mathematics, 128 (2001), pp.~311--334,
  \url{https://doi.org/10.1016/S0377-0427(00)00517-3}.

\bibitem{WobkerTurek2009}
{\sc H.~Wobker and S.~Turek}, {\em Numerical studies of {V}anka-type smoothers
  in computational solid mechanics}, Advances in Applied Mathematics and
  Mechanics, 1 (2009), pp.~29--55.

\end{thebibliography}
